\documentclass[11pt,reqno]{amsart}
\usepackage{amsmath,amssymb,amscd,verbatim,color}
\usepackage[backref]{hyperref}

\usepackage{amsfonts}
\usepackage[top=35mm, bottom=35mm, left=30mm, right=30mm]{geometry}

\usepackage{xcolor}

\theoremstyle{plain}
\newtheorem{thm}{Theorem}[section]
\newtheorem{cor}{Corollary}[section]
\newtheorem{lem}{Lemma}[section]
\newtheorem{prop}{Proposition}[section]

\newtheorem{defn}{Definition}[section]

\newtheorem{rem}{Remark}[section]

\numberwithin{equation}{section}

\newcommand{\mb}{\mathbb}
\newcommand{\mc}{\mathcal}
\def \a{\alpha} \def \b{\beta} \def \g{\gamma} \def \d{\delta}
\def \t{\theta}   \def \e{\epsilon} 
\def \s{\sigma} \def \l{\lambda}  \def \o{\omega}
\def \O{\Omega}   \def \r{\rho} \def \u{\upsilon}
  
 \def \spa{\mathrm{span}}
 
\def \L{\Lambda}

\allowdisplaybreaks

\renewcommand*{\backref}[1]{}
\renewcommand*{\backrefalt}[4]{\quad \tiny
  \ifcase #1 (\textbf{NOT CITED.})%
  \or    (Cited on page~#2.)%
  \else   (Cited on pages~#2.)%
  \fi}

\makeatletter

\def\MRbibitem{\@ifnextchar[\my@lbibitem\my@bibitem}

\def\mybiblabel#1#2{\@biblabel{{\hyperref{http://www.ams.org/mathscinet-getitem?mr=#1}{}{}{#2}}}}

\def\myhyperanchor#1{\Hy@raisedlink{\hyper@anchorstart{cite.#1}\hyper@anchorend}}

\def\my@lbibitem[#1]#2#3#4\par{%
  \item[\mybiblabel{#2}{#1}\myhyperanchor{#3}\hfill]#4%
  \@ifundefined{ifbackrefparscan}{}{\BR@backref{#3}}%
  \if@filesw{\let\protect\noexpand\immediate
    \write\@auxout{\string\bibcite{#3}{#1}}}\fi\ignorespaces%
}

\def\my@bibitem#1#2#3\par{%
  \refstepcounter\@listctr
  \item[\mybiblabel{#1}{\the\value\@listctr}\myhyperanchor{#2}\hfill]#3%
  \@ifundefined{ifbackrefparscan}{}{\BR@backref{#2}}%
  \if@filesw\immediate\write\@auxout
    {\string\bibcite{#2}{\the\value\@listctr}}\fi\ignorespaces%
}

\makeatother

\begin{document}



\author{Wen Huang} \address[Wen Huang] {  Department of Mathematics\\
 University of Science and Technology of China\\
  Hefei, Anhui, China}
\email[Wen.H]{wenh@ustc.edu.cn}

 \author{Zeng Lian} \address[Zeng Lian] {College of Mathematics\\ Sichuan University\\
    Chengdu, Sichuan, 610016, China} \email[Z.~Lian]{ZengLian@gmail.com}

\author{Kening Lu} \address[Kening Lu] {
    Department of Mathematics\\ Brigham Young University\\ Provo, Utah, 84602, United States} \email[K. Lu]{klu@math.byu.edu}

\title[Random Anosov systems mixing on fibers]{  Ergodic theory of Random Anosov systems \\ mixing on fibers}

\pagestyle{plain}

\thanks{This work is partially supported by NSF of China (11225105, 11371339, 11431012, 11671279, 11541003, 11725105).}


\begin{abstract}
In this paper,  we study the complicated dynamics of Anosov systems driven by an external force in the context of geometric theory (an abundance of random periodic points and random horseshoes) and smooth ergodic theory (random periodic measures and random Liv\v sic Theorem).
\end{abstract}


\maketitle

\setcounter{page}{1} \tableofcontents

\newpage

\parskip 5pt


\section{Introduction}\label{S:Introduction}

In this paper,  we study the complicated dynamics of Anosov systems driven by an external force in the context of geometric theory  and smooth ergodic theory.


\subsection{Description of Main Results.}
Let $(\O,d_\O)$ be a compact metric space and $\t:\O\to\O$ be a homeomorphism. Let $M$ be a connected compact smooth Riemannian manifold without boundary. A dynamical system driven by an external force $\theta$ or a cocycle  is a family of continuous maps
\[F(n,\cdot,\cdot):\Omega \times M \to M, \quad (\omega, x) \mapsto F(n, \omega, x), \quad\text{for } n\in\mathbb{Z}
\]
such that the map $F(n,\omega):=F(n,\omega,\cdot)$ forms a
cocycle over $\theta$:
\[
F(0, \omega)=Id, \quad \hbox{ for all }\; \omega \in \Omega,
\]
\[
F(n+m,\omega)=F(n,\theta^{m}\omega)F(m,\omega), \quad
\hbox{ for all }\; m, n \in \mathbb{Z}, \quad\omega \in \Omega.
\] When space $\Omega$ is endowed an invariant probability measure such as the Haar measure for $\theta$, $F$ is also called a random dynamical system, see Arnold \cite{A}.

Let $f(\omega)$ be the time-one map of the system, i.e., $f(\omega) x= F(1, \omega, x)$, which is assumed to be a diffeomorphism from $M$ to $M$. Conversely, such a family of  diffeomorphisms also generates a cocycle (or a random dynamical system):
\[
F(n, \omega, x)=
\begin{cases} f(\theta^{n-1}\omega)\circ \cdots \circ
f(\omega) x, & n >0, \\
x, & n=0, \\
f^{-1}(\theta^n\omega)\circ \cdots \circ
f^{-1}(\theta^{-1}\omega) x, & n<0.
\end{cases}
\]
Equivalently, putting $f$ and $\t$ together forms a skew product system $\phi:M\times \O\to M\times \O$  given by
\[\phi(x,\o)=(f(\o)x,\t\o)=(f_\o x,\t\o),\ \forall x\in M,\o\in\O
\]
 where we rewrite $f(\o)$ as $f_\o$ for the sake of convenience.

We consider a system $\phi(x,\o)=(f_\o x,\t\o)$ which we call {\bf Anosov on fibers}, i.e., the following hold:
for every $(x,\o) \in M\times \O$ there is a splitting of the tangential fiber of $M_\o:=M\times\{\o\}$ at $x$
\[
T_xM_\o=E_{(x,\o)}^u \oplus E_{(x,\o)}^{s}
\]
which depends continuously on $(x,\o) \in M\times \O$ with $\dim E_{(x,\o)}^u, \dim E_{(x,\o)}^s >0$ and satisfies that
\[
Df_{\o}(x) E_{(x,\o)}^u=E_{\phi(x,\o)}^u, \ \ \ Df_{\o}(x) E_{(x,\o)}^{s}= E_{\phi(x,\o)}^{s},
\]
and
\begin{equation*}
\left\{\begin{array}{ll}
|Df_{\o}(x) \xi| \geq e^{\l_0} |\xi|, \ \ \ &\forall\ \xi \in  E_{(x,\o)}^u, \\
|Df_{\o}(x)\eta| \leq  e^{-\l_0}|\eta|, \ \ \ &\forall\ \eta \in  E_{(x,\o)}^{s},
\end{array}
\right.
\end{equation*}
where $\lambda_0>0$ is a constant.

We note that if $\O$ is a differentiable manifold and $\theta$ is a diffeomorphims with expansion weaker than $e^{\lambda_0}$ and contraction weaker than $e^{-\lambda_0}$, then the system  is  a partially hyperbolic system.

We assume that $\phi:M\times \O\to M\times \O$ is   {\bf \em topological mixing on fibers}, that is,  that for any nonempty open sets $U,V\subset M$, there exists $N>0$ such that for any $n\ge N$ and $\o\in \O$
\[\phi^n(U \times \{\o\})\bigcap V\times \{\t^n\o\}\neq \emptyset.\]

In the following, we first introduce two ingredients of dynamical complexity: {\bf \em random periodic orbit} and {\bf \em random horseshoe.}

Let $L^\infty(\O, M)$ be the space of Borel measurable maps from $\O$ to $M$ endowed with the    metric $d_{L^\infty(\O, M)}$:
$$d_{L^\infty(\O, M)}(g_1,g_2)=\sup_{\o\in\O}d_M(g_1(\o),g_2(\o)),\ \forall g_1,g_2\in L^\infty(\O, M).$$
The set $\{graph(g_i)|\ g_i\in L^\infty(\O, M)\}_{1\le i\le n}$ is called  a {\bf \em  random periodic orbit} of $\phi$ with period $n$ if
 for all $1\le i\le n$
 \[\phi(\text{graph}(g_i))=\text{graph}(g_{i+1\mod n}),\ \forall 1\le i\le n.
 \]
 Equivalently,
 \begin{equation*}
 g_i \text{ is a periodic point of }\tilde \phi\text{ i.e., }\tilde \phi (g_i)=g_{i+1\mod n}\ \forall 1\le i\le n,
 \end{equation*}
 where $\tilde \phi$ is the induced map of $\phi$ on $L^\infty(\O, M)$ given by for each $g\in L^\infty(\O, M)$
 \[
 \tilde \phi(g)(\omega)=f_{\theta^{-1}\omega} g(\theta^{-1}\omega).
 \]
Each $graph(g_i)$ is called a  random periodic point of $\phi$. Without causing any confusion and for the sake of convenience, we  call  $g_i$ a random periodic point for short.
 Moreover, if $g_i$ is a continuous map then  we call   $g_i$ a continuous random periodic point.

Next, we introduce the concept of  random horseshoes on different levels, which are distinguished by the following separation functions on $L^\infty(\O, M)\times  L^\infty(\O, M)$  measuring the separations of elements of $L^\infty(\O, M)$ on different levels.
\begin{align*}\begin{split}
\overline d_{L^\infty(\O, M)}(g_1,g_2)&=\sup_{\O'\subset \O\text{ is open and nonempty}} \inf_{\o\in \O'}d_M(g_1(\o),g_2(\o)),\ \forall g_1,g_2\in L^\infty(\O, M);\\
\underline d_{L^\infty(\O, M)}(g_1,g_2)&=\inf_{\o\in\O}d_M(g_1(\o),g_2(\o)),\ \forall g_1,g_2\in L^\infty(\O, M).
\end{split}
\end{align*}
Let $\mc S_k=\{1,\ldots,k\}^{\mb Z}$ be the space of two sides sequences of $k(\geq 2))$ symbols endowed with the standard metric, and $\s$ be the (left-)shift map.

Using the separation funcation $\overline d_{L^\infty(\O, M)}$, we define a {\bf \em  full random horseshoe} with $k$-symbols of $\phi$  as a continuous embedding $\Psi:\mc S_k\to L^\infty(\O, M)$ satisfying
{\it \begin{itemize}
     \item[i)] There exists $\Delta>0$ such that for any $l\in\mb Z$ and $\hat a_1,\hat a_2\in \mc S_k\text{ with }\hat a_1(l)\neq \hat a_2(l)$
     \begin{equation*}\label{E:RandomSepaImgSymb}
     \overline d^{\tilde \phi}_{l}\left(\Psi({\hat a_1}),\Psi({\hat a_2})\right)>\Delta,
     \end{equation*}
     where for $l\ge 0$
     $$\overline d^{\tilde \phi}_{l}(g_1,g_2):=\max_{0\le i\le l}\left\{\overline d_{L^\infty(\O, M)}\left(\tilde \phi^i(g_1),\tilde \phi^i(g_2)\right)\right\},\ g_1,g_2\in L^\infty(\O, M)$$  and for $l<0$ we let $\overline d^{\tilde\phi}_{l}=\overline d^{\tilde\phi^{-1}}_{|l|};$
\item[ii)]  For all $\hat a\in\mc S_k$ $$\tilde \phi\left(\Psi({\hat a})\right)=\Psi({\s \hat a}).$$
\end{itemize}}
Then, using the stronger separation function $\underline d_{L^\infty(\O, M)}$, we define a {\bf \em  strong full random horseshoe} as a continuous embedding $\Psi:\mc S_k\to L^\infty(\O, M)$ satisfying
{\it \begin{itemize}
     \item[i)] There exists $\Delta>0$ such that for any $l\in\mb Z$ and $\hat a_1,\hat a_2\in \mc S_k\text{ with }\hat a_1(l)\neq \hat a_2(l)$
     \begin{equation*}\label{E:SRandomSepaImgSymb}
     \underline d^{\tilde \phi}_{l}\left(\Psi({\hat a_1}),\Psi({\hat a_2})\right)>\Delta,
     \end{equation*}
     where for $l\ge 0$
     $$\underline d^{\tilde \phi}_{l}(g_1,g_2):=\max_{0\le i\le l}\left\{\underline d_{L^\infty(\O, M)}\left(\tilde \phi^i(g_1),\tilde \phi^i(g_2)\right)\right\},\ g_1,g_2\in L^\infty(\O, M)$$  and for $l<0$ we let $\underline d^{\tilde\phi}_{l}=\underline d^{\tilde\phi^{-1}}_{|l|};$
\item[ii)]  For all $\hat a\in\mc S_k$ $$\tilde \phi\left(\Psi({\hat a})\right)=\Psi({\s \hat a}).$$
\end{itemize}
} A full random horseshoe we have here is an embedding of a two-sides symbolic dynamical system into $L^\infty(\O, M)$, which is an extension of the standard topological horseshoe.

Let $\mc H=\text{diff}^2(M)$ be the space of $C^2$ diffeomorphisms on $M$ with $C^2$-topology. Throughout this paper, we assume $f:\O\to \mc H$ is a continuous map.

Our first result states that there is an abundance of  random periodic points.

\noindent
{\bf Theorem A}  ({\bf Density of Random Periodic Points}).  {\it Assume that $\phi$ is Anosov and topological mixing on fibers. Then, for any $\e>0$, there exists $N\in\mb N$ such that for any $g\in L^{\infty}(\O)$ and $n\ge N$, there exists a random periodic point $\tilde g$ with period $n$ such that
 \begin{equation*}
  d_{L^{\infty}}(g,\tilde g)\le \e.
 \end{equation*}
 }

 \noindent
 {\bf Remark:} The study of hyperbolic dynamics  goes back to Poinc\'are \cite{Poincare} on $3$-body problem and Hadamard \cite{Hadamard} on geodesic flow. The modern theory of uniformly hyperbolic dynamical systems was initiated by Anosov \cite{Anosov} and Smale \cite{Smale} where Anosov and Axiom A diffeomeophism/flows were introduced respectively. The core ingredient in these systems is the uniform hyperbolicity which is an invariant geometric structure describing the exponential divergence of nearby orbits. This exponential divergence together with the compactness of phase space produces rich and complicated dynamical structures. The dynamics of these systems has been understood very well. One of results is that for a transitive Anosov diffeomorphism, the set of its periodic points (the number of the periodic orbits with same period is finite) is countable and dense in $M$. However, Theorem A tells that for the systems we study here, the set of random periodic points  is uncountable because of the non-separability of $L^\infty(\O, M)$,  where we exclude the case that $\O$ has trivial topology. In fact, Theorem A indicates that  the cardinality of random periodic orbits with sufficient large period is infinite (actually uncountable), thus the increasing rate of numbers of random periodic orbits with same period as periods tends to infinity is always infinite. Therefore such an increasing rate fails in representing the topological entropy of the system, which is one of the main tasks that periodic orbits perform for Anosov diffeomorphisms in deterministic scenario. Nevertheless, there are other candidates which are able to carry out the task of representing complexity of the system such as the random horseshoes, which will be introduced in the next.

 Another interesting observation based on an example given in Section \ref{Examples} is that the function corresponding to random periodic points are {\bf NOT} necessarily continuous. Actually, we  show that under certain conditions,  the graph corresponding to {\bf each} random periodic point is {\bf discontinuous}.

Our next result is on the existence of random horseshoes with entropy closed to the entropy of the system.

 \noindent
{\bf Theorem B} ({\bf Entropy and Random Horseshoe}).  {\it Assume that $\phi$ is  Anosov and topological mixing on fibers. Then the following statements hold
\begin{itemize}
\item[(a)] Strong Random Horseshoe: For any $\g>0$, there exist $N,k\in \mb N$ such that the following hold
   \begin{enumerate}
     \item[i)] $\frac1N\ln k\ge \underline{h}(M\times\O|\O)-\g$, where $\underline{h}(M\times\O|\O)>0$ and \[\underline{h}(M\times\O|\O):=\inf_{\o\in\O} h_{top}(\phi|_{M_\o});\]
     \item[ii)] $\phi^N$ has a $k$-symbol strong full random horseshoe.
   \end{enumerate}
\item[(b)] Random Horseshoe: For any $\g>0$,  there exist $N,k\in \mb N$ such that the following hold
   \begin{enumerate}
     \item[i)] $\frac1N\ln k\ge \overline{h}(M\times\O|\O)-\g$, where $\overline{h}(M\times\O|\O)>0$ and
      \[\overline{h}(M\times\O|\O):=\sup_{\o\in\O} h_{top}(\phi|_{M_\o});\]
     \item[ii)] $\phi^N$ has a $k$-symbol full random horseshoe.
   \end{enumerate}
 \end{itemize}
 }
 \noindent
 {\bf Remark:} The measure-theoretic entropy was introduced in 1950's by Kolmogorov \cite{Kol} and Sinai \cite{Sinai0} to measure the rate of increase in dynamical complexity as the  system evolves with time. Sinai \cite{Sinai64} studied an ergodic measure preserving automorphism  $f$ of a Lebesgue space $(X, \mu)$  and proved that if the measure-theoretic entropy of $f$ is positive, then $f$ contains factor automorphisms  which are isomorphic to Bernoulli shifts. The topological entropy was first introduced in 1965 by Adler, Konheim and McAndrew for studying dynamical systems in topological spaces. In his remarkable paper \cite{K}, A. Katok proved that for a nonuniformly hyperbolic $C^2$ diffeomorphism on a compact Riemannian manifold, the positive topological  entropy implies the existence of a Smale horseshoe. Recently, Lian and Young \cite{LY1}  extended Katok's results to  $C^2$ differentiable maps with a nonuniformly hyperbolic compact invariant set in a separable Hilbert space and to a $C^2$ semiflow  in a Hilbert space which has a nonuniformly hyperbolic compact invariant set. Without assuming any hyperbolicity, Huang and Lu \cite{HL} proved if a continuous random map in a Polish space has a positive topological entropy on a random compact invariant set ,  then it contains a weak horseshoe (see also \cite{HY} for topological dynamical systems in  compact metric spaces and \cite{KL} for $C^*$-dynamics). The horseshoe we have here is for hyperbolic systems driven by an external force, which cannot be obtained by using Katok's result or Lian and Young's result for flow. When space $\Omega$ is endowed an invariant probability measure such as the Haar measure for $\theta$, it follows from \cite{GK} that  $\phi$ has a so called random SRB measure. In particular, for the Anosov systems driven by a quasiperiodic force which appears in Section \ref{S:Example1},  the systems have a unique random SRB measure since the driven system is uniquely ergodic.

\vskip0.05in
Next, we consider random periodic $\phi$-invariant measures. We fix  a $\t$-invariant probability measure $\mb P$ on $\O$. For $g\in L^{\infty}(\O, M)$, we define a Borel probability measure $\mu_g$ in $\mc M(M\times\O)$ (the set of Borel probability measures) as follows: for any Borel set $A\subset M\times \O$
\begin{equation*}
\mu_g(A)=\mb P\left(\left\{\o\in\O\big|\ (g(\o),\o)\in A\cap M\times\{\o\}\right\}\right).
\end{equation*}
Let $ \mc I_{\mb P}(M\times\O)$ be the collection of $\phi$-invariant measures whose marginal is $\mb P$.
An invariant measure $\mu\in \mc I_{\mb P}(M\times\O)$ is called a random periodic measure if  there is a random periodic orbit $\{g_i\}_{i=1,\cdots,n}$ such that
$$\mu=\frac1n\sum_{i=1}^n\mu_{g_i}=\frac1n\sum_{i=1}^n\mu_{\tilde \phi^i(g_1)}.$$

Then, we have the following theorem which states that any $\phi$-invariant measure with marginal $\mb P$ can be approximated by random periodic measures.

\noindent{\bf Theorem C} ({\bf Density of Random Periodic Invariant Measures}).  {\it Assume that $\phi$ is  Anosov and topological mixing on fibers. Then the set of all periodic $\phi$-invariant measures is dense in $\mc I_{\mb P}(M\times\O)$ in the narrow topology. }

Unlike the standard periodic measure supported by one periodic orbit, the random periodic measure supported by one random periodic orbit is {\bf not} necessarily ergodic.

Finally, we present a weak version of the classic Liv\v sic Theorem for random dynamical systems. For a $\t$-invariant probability measure $\mb P$ and $\a\in(0,1)$, denote $L^1_{\mb P}(\O,C^{0,\a}(M))$ the collection of real functions $\Phi$ on $M\times\O$, which is H\"older on $M$ and measurable on $\O$, and satisfies the following condition
\begin{equation*}
\|\Phi\|_{L^1_{\mb P}(\O,C^{0,\a}(M))}:=\int_{\O}\|\Phi\|_{C^{0,\a}(M_\o)}d\mb P<\infty.
\end{equation*}
For each $g\in L^\infty(\O, M)$ and $\Phi\in L^1_{\mb P}(\O,C^{0,\a}(M))$, we define a functional $\tilde \Phi$ on $L^\infty(\O, M)$ as  follows
\begin{equation*}
\tilde \Phi(g)=\int_{\O}\Phi(g(\o),\o)d\mb P.
\end{equation*}

\noindent{\bf Theorem D} ({\bf Random Liv\v sic Theorem}).  {\it Assume that $\phi$ is  Anosov and topological mixing on fibers.
 Then for any $\Phi\in L^1_{\mb P}(\O,C^{0,\a}(M))$ satisfying
\begin{equation}\label{E:PeriodicNull}
\int_{M\times\O}\Phi d\mu=0\; \text{ for any random periodic } \phi\text{-invariant measure $\mu$ with marginal } \mb P,
\end{equation}
there exists a functional $\tilde\Psi:L^\infty(\O, M)\to \mb R$ which is $\a$-H\"older continuous and satisfies the following cohomological equation
\begin{equation}\label{E:CoBoundary}
\tilde\Phi(g)=\tilde\Psi(\tilde\phi(g))-\tilde \Psi(g),\ \forall g\in L^\infty(\O, M).
\end{equation}
Moreover, $\tilde \Psi$ is uniquely determined up to a constant functional.
}

\noindent
{\bf Remark:} The celebrated theorem of Liv\v sic  \cite{Liv71, Liv72a,Liv72b} was established by  Liv\v sic for hyperbolic diffeomorphisms. Since then, there is a large amount of literature on generalizations of this theorem, see, for example, \cite{AKL18,dlLW10,dlLW11,Kal11,KK96, KP16,NP13,NT95, NT96, Wil13}.
Kifer and Liu \cite{KifLiu} conjectured that for a random Anosov diffeomorphism $f_\omega$ on $\mathbb{T}^d$, and a given function $\Phi(x,\omega)$ which is  measurable in $(x,\omega)$  and $\a$-H\"older continuous  in $x$, and satisfies $\int \Phi d\mu=0$ for any $\phi$-invariant measure $\mu$ with marginal  $\mb P$, the following cohomological equation
\[
\Phi=\Psi \circ \phi - \Psi+h,
\]
has a solution $(\Psi(x, \omega), h(\omega))$,  where $\Psi(x, \omega)$ is  a measurable function satisfying the random
H\"older condition $|\Psi(x, \omega)-\Psi(y,\omega)|\leq  K_{\Psi}(\omega)(d(x, y))^\alpha$ with  $\int \log  K_{\Psi}(\omega) \;d {\mb P}< \infty$, and $h(\omega)$ is a random variable with $\int h \;d {\mb P}=0$.

Theorem D is a weak version of Liv\v sic theorem for random dynamical systems,
in the proof of which  a significant difference from the deterministic case shows up: Because of the absence of the separability the phase space $L^{\infty}(\O,M)$, there is no  dense orbit any more, which plays a key role in the proof for the deterministic case.  The proof for random case  involves a different and more complicated procedure.

\subsection{Examples.} We study following two examples: system S1) and system S2) in  section \ref{Examples}﹝
\begin{enumerate}
\item[S1)] {\bf Anosov systems diven by a quasiperiodic forcing:} $(\t,\O)$ is a minimal irrational rotation on a compact torus of dimenison $d$, and $\phi$ is Anosov on fibers and topological transitive on $M\times \O$;
\item[S2)] {\bf Measure preserving Anosov systems driven by a random forcing:} $\phi$ is Anosov on fibers, and there exists an $f$-invariant Borel probability measure $\nu$ with full support  (i.e. $supp\nu=M$).
\end{enumerate}
We show that both system $S1$ and system $S2$ are topological mixing on fibers. Thus, {\bf Theorem A-D} hold for both systems. More specific examples are also given.  For system S1) we consider a $2$-$d$ tori automorphism driven by a quasiperiodic force and for system S2) we look at the random composition of $2\times 2$ area-preserving positive matrices with integer entries.

\noindent
{\bf Remark on Horseshoe:}
 System S1) is an Anosov system driven by a quasi-periodic forcing and is also a partially hyperbolic system. In fact, with each invariant Borel probability measure $\mu$, system S1) has zero Lyapunov exponent of multiplicity $d$ (the dimension of $\Omega=T^d$). Thus, the Katok's result (\cite{K}) does not apply to this system. Lian and Young's results (\cite{LY}, \cite{LY1}) for semiflows or flows, which allow only one-dimensional center associated to the flow direction, are also not applicable to this case because of high dimensional center. Futhermore, the quasi-periodic forced systems  has no periodic points, and so has no standard horseshoe. However, because of the uniform hyperbolicity of the map between fibers ($f(\o):M_\o\to M_{\t\o}$) and the rigidity of the evolution of quasiperiodic forcing, one expects the co-existence of certain kinds of periodicity and chaotic behavior over some random invariant structures.

 The main difficulty appearing in these systems is that the pseudo orbits derived from  Poincar\'e Recurrence Theorem only provides returns in $M\times \Omega$, thus does not have the shadowing property since the hyperbolicity only happens on fibers which are identified by $M$. To overcome this difficulty, instead of looking at one single pseudo orbit, we construct a group of globally returned pseudo orbits for all $\o\in \Omega$ at the same time.

 For non-uniformly hyperbolic systems, one needs to look at the return map on Pesin blocks (as being done in \cite{K}). If there is a globally returning time, the projection of such  Pesin blocks to $\O$ should be a $\mb P$-full-measure-set. This is not satisfied even in the setting of \cite{GK} in which hyperbolicity depends on the sample ``$\o$". Such observation reveals that for non-uniform hyperbolic system, accurate information about ``shapes" of the Pesin blocks is required, which makes the situation much more complicated than autonomous case. Second, even dealing with systems which is Anosov on fibers, constructing a global return is not trivial, since one needs infinite many (usually uncountable many) orbits return at a same time. Our immediate feeling is that for non-uniformly hyperbolic systems on fibers, one may need a new approach.

Finally, we mention some existing results which are closely relevant  to the results derived in this paper, which may shed a light on further research directions, for example, chaotic behavior of Anosov flows driven by quasi-periodic noises. Of course, this is far from being a complete picture. In \cite{G95}, Gundlach studied the random homoclinic points for smooth invertible random dynamical systems, and derived a random version of {\em Birkhoff-Smale Theorem}. This result is applicable to systems which are provided by small random perturbations of systems with homoclinic points.  In the recent papers \cite{LW1} and \cite{LW2}, Lu and Wang studied the chaos behavior of a type of differential equations driven by a nonautononous force deduced from  stochastic processes which are the truncated and classical Brownian motions  in \cite{LW1}  and in \cite{LW2} respectively. The authors extended the concept of topological horseshoe and showed that such structure (thus chaotic behavior) exists almost surely. Moreover, they further applied the results to randomly forced Duffing equations and pendulum equations. There are also many works on the existence of random periodic solutions for stochastic differential equations, see \cite{ZhZh}



 \section{Setting, Definitions,  and Statement of the Main Results}\label{S:Setting}
In this section, we set up the systems, introduce notions, and give  the statements of the main results precisely. Although some of the definitions have been given in Section \ref{S:Introduction}, we  restate them formally for the sake of completeness and convenience.

Let $M$ be a connected compact smooth Riemannian manifold with no boundary,  $d_M$ be the induced Riemannian metric of $M$,  $(\O,d_\O)$ be a compact metric space, and $\t:\O\to\O$ be a homeomorphism.
 Let $\mc H=\text{diff}^2(M)$ be the space of $C^2$ diffeomorphisms on $M$ with $C^2$-topology, and $f:\O\to \mc H$ be a continuous map. Then the skew product of $f$ and $\t$ induces a map $\phi:M\times \O\to M\times \O$  in the following way:
 $$\phi(x,\o)=(f(\o)x,\t\o)=(f_\o x,\t\o),\ \forall x\in M,\o\in\O$$
 where we rewrite $f(\o)$ as $f_\o$ for the sake of convenience. Clearly, $\phi$ induces a cocycle over $(\O,\t)$, which can be inductively defined by
$$\begin{cases}
\phi^0(x,\o)&=(x,\o)\\
\phi^n(x,\o)&=(f_{\t^{n-1}\o}\phi^{n-1}(x,\o),\t^n\o),\ n\ge 1\\
\phi^n(x,\o)&=(f^{-1}_{\t^{n+1}\o}\phi^{n+1}(x,\o),\t^n\o),\ n\le -1
\end{cases}.
$$
Note that $M\times \O$ is a compact metric space when we introduce a metric $d$ by defining that
 $$d((x_1,\o_2),(x_2,\o_2))=d_M(x_1,x_2)+d_\O(\o_1,\o_2), \forall (x_1,\o_1),(x_2,\o_2)\in M\times \O.$$
 Then $\phi$  and $D_M\phi(\text{where }D_M\phi(x,\o):=Df_\o(x))$ are continuous on $M\times \O$.
\begin{defn}\label{D:Anosov}
The system generated by $\phi$ (or simply $\phi$) is called  {\bf Anosov on fibers} if the following holds:
for every $(x,\o) \in M\times \O$ there is a splitting of the tangential fiber of $M_\o:=M\times\{\o\}$ at $x$
\[
T_xM_\o=E_{(x,\o)}^u \oplus E_{(x,\o)}^{s}
\]
which depends continuously on $(x,\o) \in M\times \O$ with $\dim E_{(x,\o)}^u,\dim E_{(x,\o)}^s >0$ and satisfies that
\[
Df_{\o}(x) E_{(x,\o)}^u=E_{\phi(x,\o)}^u, \ \ \ Df_{\o}(x) E_{(x,\o)}^{s}= E_{\phi(x,\o)}^{s},
\]
and
\begin{equation}\label{E:UniformHyperbolic}
\left\{\begin{array}{ll}
|Df_{\o}(x) \xi| \geq e^{\l_0} |\xi|, \ \ \ &\forall\ \xi \in  E_{(x,\o)}^u, \\
|Df_{\o}(x)\eta| \leq  e^{-\l_0}|\eta|, \ \ \ &\forall\ \eta \in  E_{(x,\o)}^{s},
\end{array}
\right.
\end{equation}
where $\lambda_0>0$ is a constant.
\end{defn}
\begin{rem}\label{R:ContInvDistr}
It is clear that the splitting $T_xM_\o=E_{(x,\o)}^u \oplus E_{(x,\o)}^{s}$ is uniformly continuous on $M\times\O$. Actually, one has that the above splitting is H\"older continuous on $x$, which follows from the standard arguments. Moreover, the H\"older constant and exponent can be chosen independent  on $(x,\o)$. For more details, we refer the readers to \cite{Brin}, Section 5.3 of \cite{BaPe}, and Section  5.4 of \cite{LLL}.
\end{rem}

Denote $\pi^u_{(x,\o)},\pi^s_{(x,\o)}$ the projections  associated to the splitting $T_xM_\o=E_{(x,\o)}^u \oplus E_{(x,\o)}^{s}$, which map $T_xM_\o$  onto $E_{(x,\o)}^u$ and  $E_{(x,\o)}^{s}$ respectively. Since the splitting $T_xM_\o=E_{(x,\o)}^u \oplus E_{(x,\o)}^{s}$ is uniformly continuous on $M\times\O$, there exists a real number $\mc P>1$ such that
\begin{equation}\label{E:PrjBound}
\max\left\{\|\pi_{(x,\o)}^{\tau}\|\big|\ (x,\o)\in M\times\O, \tau=u,s\right\}<\mc P.
\end{equation}

\medskip

The general  systems we consider is topological mixing on fibers which is defined as follows:
\begin{defn}\label{D:Mixing}
$\phi:M\times \O\to M\times \O$ is called  {\bf topological mixing on fibers} if for any nonempty open sets $U,V\subset M$, there exists $N>0$ such that for any $n\ge N$ and $\o\in \O$
$$\phi^n(U \times \{\o\})\bigcap V\times \{\t^n\o\}\neq \emptyset.$$
\end{defn}
In the following, we first introduce two ingredients of dynamical complexity: {\bf \em random periodic orbit} and {\bf \em random horseshoe.}

Let $L^\infty(\O, M)$ be the space of Borel measurable maps from $\O$ to $M$ endowed with the    metric $d_{L^\infty(\O, M)}$:
$$d_{L^\infty(\O, M)}(g_1,g_2)=\sup_{\o\in\O}d_M(g_1(\o),g_2(\o)),\ \forall g_1,g_2\in L^\infty(\O, M).$$
For later use, we also define two functions on $L^\infty(\O, M)\times  L^\infty(\O, M)$ to measure the separations of elements of $L^\infty(\O, M)$ on different levels:
\begin{align}\begin{split}\label{E:SeparationFunc}
\overline d_{L^\infty(\O, M)}(g_1,g_2)&=\sup_{\O'\subset \O\text{ is open and nonempty}} \inf_{\o\in \O'}d_M(g_1(\o),g_2(\o)),\ \forall g_1,g_2\in L^\infty(\O, M).\\
\underline d_{L^\infty(\O, M)}(g_1,g_2)&=\inf_{\o\in\O}d_M(g_1(\o),g_2(\o)),\ \forall g_1,g_2\in L^\infty(\O, M).
\end{split}
\end{align}

Since $\phi$  and $\phi^{-1}$ are both continuous, $\phi$ induces a map $\tilde \phi$ on $L^\infty(\O, M)$ given by
 \[
 \tilde \phi(g)(\omega)=f_{\theta^{-1}\omega} g(\theta^{-1}\omega).
 \]


\begin{defn}\label{D:RandomPeriodic}
$\{graph(g_i)|\ g_i\in L^\infty(\O, M)\}_{1\le i\le n}$ is called  a {\em random periodic orbit} of $\phi$ with period $n$ if
 for all $1\le i\le n$ $$\phi(\text{graph}(g_i))=\text{graph}(g_{i+1\mod n})\ \forall 1\le i\le n,$$
 equivalently,
 \begin{equation}\label{E:PeriodicPointTPh}
 g_i \text{ is a periodic point of }\tilde \phi\text{ i.e. }\tilde \phi (g_i)=g_{i+1\mod n}\ \forall 1\le i\le n.
 \end{equation}
\end{defn}
Each $graph(g_i)$ is called a random periodic point of $\phi$. Without causing any confusion and for the sake of convenience, we  call  $g_i$ a random periodic point for short.
 Moreover, if $g_i$ is a continuous map then  we call   $g_i$ a continuous random periodic point.

\begin{rem}\label{R:RandomDynSys}
The system generated by $\phi$ is usually called a skew-product system rather than a random dynamical system, which makes the name " random periodic orbit" sound a little bit more in name than reality. Nevertheless, once a $\t$-invariant measure $\mb P$ on $\O$ is introduced,  the resulting system induced by $\phi$ becomes a random dynamical systems over the  metric dynamical system $(\O,\mb P,\t)$ (see \cite{A} for more details).  On the other hand, since the name "random periodic orbit" has been introduced in exiting literature (see \cite{Kl, ZhZh}) and has the essentially equivalent meaning,  we prefer not to introduce different name or notions.
\end{rem}

Next, we introduce the concept of  random horseshoes in different levels, which are distinguished by the separation functions given by (\ref{E:SeparationFunc}).

 Let $\mc S_k=\{1,\ldots,k\}^{\mb Z}$ be the space of two sides sequences of $k(\geq 2))$ symbols endowed with the standard metric, and $\s$ be the (left-)shift map. The following is the standard definition of Smale horseshoe for the deterministic map $\tilde \phi$.
 \begin{defn}\label{D:Horseshoe}
 We  call a continuous embedding $\Psi:\mc S_k\to L^\infty(\O, M)$ a {\bf \em  horseshoe} with $k$-symbols of $\tilde \phi$ if
the following hold:
\begin{itemize}
     \item[i)] There exists $\Delta>0$ such that for any $l\in\mb Z$ and $\hat a_1,\hat a_2\in \mc S_k\text{ with }\hat a_1(l)\neq \hat a_2(l)$,
     \begin{equation}\label{E:SeparationImgSymbols}
     d^{\tilde \phi}_{l}\left(\Psi({\hat a_1}),\Psi({\hat a_2})\right)>\Delta,
     \end{equation}
     where for $l\ge 0$
     $$d^{\tilde \phi}_{l}(g_1,g_2):=\max_{0\le i\le l}\left\{d_{L^\infty(\O, M)}\left(\tilde \phi^i(g_1),\tilde \phi^i(g_2)\right)\right\},\ g_1,g_2\in L^\infty(\O, M)$$ is the Bowen metric induced by $\tilde \phi$ and for $l<0$ we let $d^{\tilde\phi}_{l}=d^{\tilde\phi^{-1}}_{|l|};$
\item[ii)]  For all $\hat a\in\mc S_k$ $$\tilde \phi\left(\Psi({\hat a})\right)=\Psi({\s \hat a}).$$
\end{itemize}
\end{defn}

 The next definition we introduce here is the random horseshoe for a skew-product system, which has  a slightly finer structure by employing the separation function $\overline d_{L^\infty(\O, M)}$ in (\ref{E:SeparationFunc}) to replace $d_{L^\infty(\O, M)}$.
 \begin{defn}\label{D:RandomHorseshoe}
 We  call a continuous embedding $\Psi:\mc S_k\to L^\infty(\O, M)$ a {\bf \em  full random horseshoe} with $k$-symbols of $\phi$ if
the following hold:
\begin{itemize}
     \item[i)] There exists $\Delta>0$ such that for any $l\in\mb Z$ and $\hat a_1,\hat a_2\in \mc S_k\text{ with }\hat a_1(l)\neq \hat a_2(l)$
     \begin{equation}\label{E:RandomSepaImgSymb}
     \overline d^{\tilde \phi}_{l}\left(\Psi({\hat a_1}),\Psi({\hat a_2})\right)>\Delta,
     \end{equation}
     where for $l\ge 0$
     $$\overline d^{\tilde \phi}_{l}(g_1,g_2):=\max_{0\le i\le l}\left\{\overline d_{L^\infty(\O, M)}\left(\tilde \phi^i(g_1),\tilde \phi^i(g_2)\right)\right\},\ g_1,g_2\in L^\infty(\O, M)$$  and for $l<0$ we let $\overline d^{\tilde\phi}_{l}=\overline d^{\tilde\phi^{-1}}_{|l|};$
\item[ii)]  For all $\hat a\in\mc S_k$ $$\tilde \phi\left(\Psi({\hat a})\right)=\Psi({\s \hat a}).$$
\end{itemize}
\end{defn}

\noindent
Finally, by using the stronger separation function $\underline d_{L^\infty(\O, M)}$ as in (\ref{E:SeparationFunc}), we define a stronger random horseshoe as follows:
 \begin{defn}\label{D:SRandomHorseshoe}
 We  call a continuous embedding $\Psi:\mc S_k\to L^\infty(\O, M)$ a {\bf \em  strong full random horseshoe} with $k$-symbols of $\phi$ if
the following hold:
\begin{itemize}
     \item[i)] There exists $\Delta>0$ such that for any $l\in\mb Z$ and $\hat a_1,\hat a_2\in \mc S_k\text{ with }\hat a_1(l)\neq \hat a_2(l)$
     \begin{equation}\label{E:SRandomSepaImgSymb}
     \underline d^{\tilde \phi}_{l}\left(\Psi({\hat a_1}),\Psi({\hat a_2})\right)>\Delta,
     \end{equation}
     where for $l\ge 0$
     $$\underline d^{\tilde \phi}_{l}(g_1,g_2):=\max_{0\le i\le l}\left\{\underline d_{L^\infty(\O, M)}\left(\tilde \phi^i(g_1),\tilde \phi^i(g_2)\right)\right\},\ g_1,g_2\in L^\infty(\O, M)$$  and for $l<0$ we let $\underline d^{\tilde\phi}_{l}=\underline d^{\tilde\phi^{-1}}_{|l|};$
\item[ii)]  For all $\hat a\in\mc S_k$ $$\tilde \phi\left(\Psi({\hat a})\right)=\Psi({\s \hat a}).$$
\end{itemize}
\end{defn}


\begin{rem}\label{R:CompareHorseshoes}
Clearly, we have that
$d_{L^\infty(\O, M)}\ge\overline d_{L^\infty(\O, M)}\ge\underline d_{L^\infty(\O, M)}.$ Therefore, we have that  a strong full random horseshoe is automatically a full random horseshoe, and a full random horseshoe of $\phi$ automatically induces a horseshoe of $\tilde \phi$.
\end{rem}

Before stating the main results,  we  introduce several notations. Let
 \begin{align*}
 &N(\o,\e,n):=\\&\quad\quad\quad\quad\quad\max\left\{card(E)\big|\ E\subset M_{\o}, E \text{ is an }\e\text{ separated set with respect to }d^{\phi}_{M\times \O,n}\right\},
 \end{align*}
where $d^\phi_{M\times \O,n}$ is the Bowen metric. Let
$$h_{top}(\phi|_{M_\o})=\lim_{\e\to 0^+}\limsup_{n\to\infty}\frac1n\log N(\o,\e,n),$$
$$\underline{h}(M\times\O|\O):=\inf_{\o\in\O} h_{top}(\phi|_{M_\o}),$$
$$ \overline{h}(M\times\O|\O):=\sup_{\o\in\O} h_{top}(\phi|_{M_\o}).$$
Here $h_{top}(\phi|_{M_\o})$ is called the topological fiber  entropy on the fiber $M_\o$, $\underline{h}(M\times\O|\O)$ and $\overline{h}(M\times\O|\O)$ are called lower and upper topological fiber entropy respectively.

The systems we consider here satisfy the following two hypotheses:
\begin{itemize}
\item[{\bf H1)}] $\phi$ is Anosov on fibers.
\item[{\bf H2)}] $\phi$ is topological mixing on fibers.
\end{itemize}
\begin{thm}[Density of Random Periodic Points]\label{T:TheoryAnosovMix}
 Let $\phi$ satisfy condition  H1) and H2). Then, for any $\e>0$, there exists $N\in\mb N$ such that for any $g\in L^{\infty}(\O)$ and $n\ge N$, there exists a random periodic point $\tilde g$ with period $n$ such that
 \begin{equation*}\label{E:LInfiClose2}
  d_{L^{\infty}}(g,\tilde g)\le \e.
 \end{equation*}
 \end{thm}

\begin{thm}[Entropy and Random Horseshoe]\label{T:TheoryAnosovMix2}
 Let $\phi$ satisfy condition H1)and H2). Then the following assertions hold:
 \begin{itemize}
  \item[{\bf A.}]  For any $\g>0$, there exist $N,k\in \mb N$ such that the following hold
   \begin{enumerate}
     \item[i)] $\frac1N\ln k\ge \underline{h}(M\times\O|\O)-\g$, where $\underline{h}(M\times\O|\O)>0$;
     \item[ii)] $\phi^N$ has a $k$-symbol strong full random horseshoe.
   \end{enumerate}
   \item[{\bf B.}] For any $\g>0$,  there exist $N,k\in \mb N$ such that the following hold
   \begin{enumerate}
     \item[i)] $\frac1N\ln k\ge \overline{h}(M\times\O|\O)-\g$;
     \item[ii)] $\phi^N$ has a $k$-symbol full random horseshoe.
   \end{enumerate}
 \end{itemize}
 \end{thm}

Next, we fix  a $\t$-invariant probability measure $\mb P$ on $\O$, and consider the random dynamical systems over $\theta$.

For $g\in L^{\infty}(\O)$, we define a Borel probability measure, $\mu_g$, in $\mc M(M\times\O)$ as follows: for any Borel set $A\subset M\times \O$
\begin{equation}\label{E:Mug}
\mu_g(A)=\mb P\left(\left\{\o\in\O\big|\ (g(\o),\o)\in A\cap M\times\{\o\}\right\}\right).
\end{equation}
Let $ \mc I_{\mb P}(M\times\O)$ be the collection of $\phi$-invariant measures whose marginal are $\mb P$. Now, we define the periodic measures.
\begin{defn}\label{D:PeriodicMeasure}
$\mu\in \mc I_{\mb P}(M\times\O)$ is called a random periodic measure if and only if there is a random periodic orbit $\{g_i\}_{i=1,\cdots,n}$ such that
$$\mu=\frac1n\sum_{i=1}^n\mu_{g_i}=\frac1n\sum_{i=1}^n\mu_{\tilde \phi^i(g_1)}.$$
\end{defn}
The following theorem states that any $\phi$-invariant measure with marginal $\mb P$ can be approximated by random periodic measures.
\begin{thm}[Density of Random Periodic Invariant Measures]\label{T:PeriodicApprox}
Let $\phi$ satisfy  H1) and  H2). Then, the set of all periodic measures is dense in $\mc I_{\mb P}(M\times\O)$ in the narrow topology.
\end{thm}
\begin{rem}\label{R:NonErgPer}
We remark here that, unlike the standard periodic measure supported by one periodic orbit, the random periodic measure supported by one random periodic orbit is {\bf not} necessarily ergodic.
\end{rem}

Finally,  we present a weak version of the classic Liv\v sic Theorem for random dynamical systems. For the $\t$-invariant probability measure $\mb P$ and $\a\in(0,1)$, denote $L^1_{\mb P}(\O,C^{0,\a}(M))$ the collection of functions $\Phi$ on $M\times\O$, which are H\"older on $M$ and measurable on $\O$, and satisfy the following condition
\begin{equation}\label{E:L1Holder}
\|\Phi\|_{L^1_{\mb P}(\O,C^{0,\a}(M))}:=\int_{\O}\|\Phi\|_{C^{0,\a}(M_\o)}d\mb P<\infty,
\end{equation}
where
$$\|\Phi\|_{C^{0,\a}(M_\o)}:=\sup_{x\in M_\o}|\Phi(x,\o)|+\sup_{x,y\in M_\o,\ x\neq y}\frac{|\Phi(x,\o)-\Phi(y,\o)|}{\left(d_M(x,y)\right)^\a}.$$
For any $g\in L^\infty(\O, M)$ and $\Phi\in L^1_{\mb P}(\O,C^{0,\a}(M))$, define a functional $\tilde \Phi$ on $L^\infty(\O, M)$ as follows
\begin{equation}\label{E:Phi(g)}
\tilde \Phi(g)=\int_{\O}\Phi(g(\o),\o)d\mb P.
\end{equation}
\begin{rem}\label{R:WellDefined}
First, one needs to show that the above integration is well defined and is finite. Since $\Phi(x,\o)$ is continuous on $x$ and Borel measurable on $\o$, we have that $\Phi$ is Borel measurable on $M\times\O$.  Then, for any real number  $a\in \mb R$,
$U:=\Phi^{-1}\left((-\infty,a)\right)$ is a Borel subset of $M\times \O$.
Note that $graph(g)$ is also a Borel subset of $M\times \O$ since $g:\O\to M$ is Borel. Therefore, $U\cap graph(g)$ is Borel in $M\times\O$. By applying Theorem \ref{T:BMST}, we have that $$\left\{\o|\ \Phi(g(\o),\o)<a\right\}= \pi_\O\left(U\cap graph(g)\right)\text{ is Borel in }\O.$$
Hence $\Phi(g(\cdot),\cdot):\O\to \mb R$ is Borel measurable. As $\Phi$ satisfying (\ref{E:L1Holder}), we have that (\ref{E:Phi(g)}) is well defined and is finite.


\end{rem}

We call the following theorem a weak version of random Liv\v sic Theorem. In the case that there is no confusion caused, we simply use "periodic measures with marginal $\mb P$" short for  "periodic $\phi$-invariant Borel probability measure with marginal $\mb P$"  for the sake of convenience.
\begin{thm}[Random Liv\v sic Theorem]\label{T:WLivTh}
Let $\phi$ satisfy condition  H1) and H2). Then for any $\Phi\in L^1_{\mb P}(\O,C^{0,\a}(M))$ satisfying
\begin{equation}\label{E:PeriodicNull}
\int_{M\times\O}\Phi d\mu=0\text{ when }\mu\text{ is periodic with marginal } \mb P,
\end{equation}
there exists a functional $\tilde\Psi:L^\infty(\O, M)\to \mb R$ which is $\a$-H\"older continuous and satisfies the following cohomological equation
\begin{equation}\label{E:CoBoundary}
\tilde\Phi(g)=\tilde\Psi(\tilde\phi(g))-\tilde \Psi(g),\ \forall g\in L^\infty(\O, M).
\end{equation}
Moreover, $\tilde \Psi$ is uniquely determined up to a constant functional.
\end{thm}
\section{Preliminary  Lemmas}\label{S:Preliminary}

This section is devoted to the preparation on technical tools needed for the proof of main theorems.  These fundamental tools stated in this section only require the systems satisfy H1).
\subsection{Invariant Manifolds}\label{S:InvMan}
We define the local stable manifolds and unstable manifolds as follows:
\begin{align*}
&W^s_\e(x,\o)=\{(y,\o)\in M_\o|\ d_M(\pi_M\phi^n(y,\o),\pi_M\phi^n(x,\o))\le \e\text{ for all }n\ge 0\}\\
&W^u_\e(x,\o)=\{(y,\o)\in M_\o|\ d_M(\pi_M\phi^n(y,\o),\pi_M\phi^n(x,\o))\le \e\text{ for all }n\le 0\},
\end{align*}
where $\pi_M:M\times \O\to M$ is the nature coordinate projection.
The following lemma can mainly be viewed as a special version of Theorem 3.1 from \cite{GK}. The only difference is that these local manifolds depends on $\o$ measurably in \cite{GK} while depends on $\o$ continuously in this paper. The reason is that the invariant splitting $E^u_{(x,\o)}\oplus E^s_{(x,\o)}$ varies continuously in both $x$ and $\o$, and $\phi$ and $D_M\phi$ are also continuous in $\o$. The proof of iii) follows from the standard procedure for which we refer to the proof of Lemma 3.1 in \cite{LLL}. Thus, we omit the proof  of the following Lemma \ref{L:InvMani}.
\begin{lem}\label{L:InvMani}
For any $\l\in (0,\l_0)$, there exists $\e_0>0$ such that for any $\e\in(0,\e_0]$, the following hold:
\begin{itemize}
\item[i)] $W^s_\e(x,\o), W^u_\e(x,\o)$ are $C^2$-embedded discs for all $(x,\o)\in M\times \O$ which satisfy \[T_xW_\e^\tau(x,\o)=E^\tau_{(x,\o)}, \text{for}
    ;\tau=u,s.\]
Moreover,  there exist a constant $L>1$ and $C^2$ maps
$$g^u_{(x,\o)}:B^u_{\mc P\e}(x,\o)\to E^s_{(x,\o)} \text{ and }g^s_{(x,\o)}:B^s_{\mc P\e}(x,\o)\to E^u_{(x,\o)} $$
such that $$W^\tau_\e(x,\o)\subset Exp_{(x,\o)}\left(graph(g^\tau_{(x,\o)})\right),$$ and $\|Dg^\tau_{(x,\o)}\|<\frac1{10\mc P}$,  $Lip Dg^\tau_{(x,\o)}<L$ for $\tau=u,s$, where $\mc P$ is as in (\ref{E:PrjBound}) and $B^\tau_{\mc P \e}(x,\o)$ is the $\mc P\e$-ball of $E^\tau_{(x,\o)}$ centered at  origin;
\item[ii)] $d_M(\pi_M\phi^n(x,\o),\pi_M\phi^n(y,\o))\le e^{-n\l}d_M(x,y)$ for $(y,\o)\in W^s_\e(x,\o)$, $n\ge 0$, and\\
$d_M(\pi_M\phi^{-n}(x,\o),\pi_M\phi^{-n}(y,\o))\le e^{-n\l}d_M(x,y)$ for $(y,\o)\in W^u_\e(x,\o)$, $n\ge 0$;
\item[iii)] $W^s_\e(x,\o), W^u_\e(x,\o)$ vary continuously on $(x,\o)$ (in $C^1$ topology).
\end{itemize}
\end{lem}
The local stable and unstable manifolds can be used to construct the global stable and unstable manifolds respectively,
\begin{align*}
&W^s(x,\o)=\{(y,\o)\in M_\o|\ |d_M(\pi_M\phi^n(y,\o),\pi_M\phi^n(x,\o))\to 0\text{ as }n\to \infty\}\\
&W^u(x,\o)=\{(y,\o)\in M_\o|\ |d_M(\pi_M\phi^n(y,\o),\pi_M\phi^n(x,\o))\to 0\text{ as }n\to -\infty\},
\end{align*}
as
\begin{align*}
&W^s(x,\o)=\bigcup_{n=0}^\infty \phi^{-n}(W^s_\e(\phi^n(x,\o))),\\
&W^u(x,\o)=\bigcup_{n=0}^\infty \phi^{n}(W^u_\e(\phi^{-n}(x,\o))),
\end{align*}
where $\e$ is an arbitrarily fixed small positive number as in \ref{L:InvMani}.
The following lemma provides local canonical coordinates on $M_\o$.
\begin{lem}\label{L:LocCoor}
For any $\e\in(0,\e_0)$ there is a $\d\in (0,\e)$ such that  for any $x,y\in M$ with $d_M(x,y)<\d$,  $W^s_\e(x,\o)\bigcap W^u_\e(y,\o)$ consists of a single point, which is denoted by $[x,y]_\o$. Furthermore
$$[\cdot,\cdot]_{\cdot}:\{(x,y,\o)\in M\times M\times \O|\ d_M(x,y)<\d\}\to M\text{ is continuous}.$$
We say the system $(M\times \O,\phi)$ has {\bf local product structure} with size $\d$.
\end{lem}
\begin{proof}
This simply follows from  iii) of Lemma \ref{L:InvMani}, the continuous (thus uniformly continuous) dependence of $W^s_\e(x,\o),W^u_\e(x,\o)$ on $(x,\o)$, and the uniform continuity of the invariant splitting $E^u_{(x,\o)}\oplus E^s_{(x,\o)}$. For details, we refer to \cite{Bow4}.
\end{proof}
\begin{lem}\label{L:Expansive}
There is an $\e>0$ such that for any $\o\in\O$ and $x,y\in M_\o$ with $y\neq x$, the following holds:
$$d_M(\pi_M\phi^k(x,\o),\pi_M\phi^k(y,\o))>\e\text{ for some }k\in\mb Z.$$
\end{lem}
\begin{proof}
Take $\e$ smaller than $\e_0$ as in Lemma \ref{L:InvMani}. Otherwise, $(y,\o)\in W^s_\e(x,\o)\cap W^u_\e(x,\o)$, therefore $y=x$.
\end{proof}
\vskip0.1in

\subsection{Shadowing Lemma}\label{S:Shadow}
We first introduce some relevant concepts. For any $\a>0$, an orbit $\{(x_i,\t^i\o)\}_{i\in \mb Z}\subset M\times \O$ is called {\bf $(\o,\a)$-pseudo orbit}  if for any $i\in \mb Z$
$$d_M(f_{\t^i\o}(x_i),x_{i+1})<\a.$$
\begin{lem}\label{L:Shadowing}
For any $\b>0$, there exists an $\a>0$ such that for any $(\o,\a)$-pseudo orbit $\{(x_i,\t^i\o)\}_{i\in \mb Z}$, there is a true orbit $\{(y_i,\t^i\o)\}_{i\in \mb Z}$ of $\phi$ such that
$$\sup_{i\in \mb Z}d_M(y_i,x_i)<\b.$$
$\{(y_i,\t^i\o)\}_{i\in \mb Z}$ is called {\bf $(\o,\b)$-shadowing} orbit of $\{(x_i,\t^i\o)\}_{i\in \mb Z}$.\\
 Further more, there is a $\b_0>0$ such that for any $0<\b<\b_0$, the above true orbit is unique.
\end{lem}
\begin{proof}
The proof follows from the proof of Proposition 3.6 in \cite{GK}, thus is omitted here.
\end{proof}
\begin{rem}\label{R:Beta0}
For the sake of convenience, we remark here that $\b_0$ is always chosen to be less than $\frac13\d$, where $\d$ is the one in Lemma \ref{L:LocCoor} for $\e=\frac12\e_0$.
\end{rem}
The next lemma gives the detailed relation between "$\a$" and "$\b$" in Lemma \ref{L:Shadowing} above.
\begin{lem}\label{L:ClosingLemma}
There exist $\a_0>0$ and $C\ge 1$ so that for any $\a\in (0,\a_0)$ and any $(\o,\a)$-pseudo orbit $\{(x_i,\t^i\o)\}_{i\in\mb Z}$, there exists a true orbit $\{(y_i,\t^i\o)\}_{i\in\mb Z}$ of $\phi$ satisfying the following
$$\sup_{i\in \mb Z} d_M(y_i,x_i)<C\a.$$
\end{lem}
\begin{proof}
By i) and iii) of Lemma \ref{L:InvMani}, for any $\e_1>0$, there exist $\d_1, L_1>0$ such that the following properties hold
\begin{itemize}
\item[P1)] For any $x,y\in M$ with $d_M(x,y)<\d_1$,  $[x,y]_\o=W^s_{\e}(x,\o)\cap W^u_{\e}(y,\o)$ is well defined, where $\e$ is as in Lemma \ref{L:InvMani};

\item[P2)] For any $x,y\in M$ with $d_M(x,y)<\d_1$
$$\max\left\{d_M([x,y]_\o,x),d_M([x,y]_\o,y)\right\}\le L_1d_M(x,y);$$

\item[P3)] For any $x,y\in M$ with $d_M(x,y)<\d_1$
$$\max\left\{\left|d_M([x,y]_\o,x)-d_M([y,x]_\o,y)\right|,\ \left|d_M([y,x]_\o,x)-d_M([x,y]_\o,y)\right|\right\}\le \e_1d_M(x,y).$$
\end{itemize}
It is obvious that P1) holds when taking $\d_1<\d$ where $\d$ is as in Lemma \ref{L:LocCoor}.
The proof of the  existence of such $\d_1$ and $L_1$ to fit P2) and P3) may be standard, nevertheless, we  sketch the proof briefly for the sake of completeness. Note that  the above properties only appear in a local scenario, for the sake of simplicity,  we  fix an atlas of the manifold and always keep our discussion in one single local chart in the rest of the proof. Without causing any confusion, we identify $T_xM_\o$ and $T_yM_\o$ to $\mb R^n$ for $d_M(x,y)$ small, thus such $x$ and $y$ are considered to be two points in $\mb R^n$ thus $x-y$ are well defined, and $\exp{(E^u_{(x,\o)})},\ \exp{(E^s_{(x,\o)})}$ are considered to be  hyperplanes in $\mb R^n$, which contain the point $x\in \mb R^n$. We also note that the distance in the local chart $\mb R^n$ is equivalent to the metric $d_M$ and the comparison constant is some system constant, therefore we will not distinguish  them in the local argument.\\

Now we propose the condition of $\d_1$ and $L_1$ which make P2) valid. As discussed in Remark \ref{R:ContInvDistr}, the splitting $T_xM_\o=E^u_{(x,\o)}\oplus E^s_{(x,\o)}$ is uniformly continuous on $(x,\o)$, thus one can make $E^u_{(x,\o)}\oplus E^s_{(x,\o)}$ closed enough to $E^u_{(y,\o)}\oplus E^s_{(y,\o)}$ by taking $d_M(x,y)$ small enough. Precisely, consider the splittings
$$\mb R^n=E^u_{(x,\o)}\oplus E^s_{(y,\o)},\text{ with associated projections } \pi^u_{x,y}, \pi^s_{x,y},$$
and
$$\mb R^n=E^u_{(y,\o)}\oplus E^s_{(x,\o)},\text{ with associated projections } \pi^u_{y,x}, \pi^s_{y,x},$$
then we request that for any $d_M(x,y)<\d_1$, the following hold
\begin{equation}\label{E:ErrProj}
\max\left\{\|\pi^\tau_{x,y}-\pi^\tau_{(z,\o)}\|,\|\pi^\tau_{y,x}-\pi^\tau_{(z,\o)}\||\right\}<\e',  \tau=u,s,\ z=x,y,
\end{equation}
where $\e'>0$ is a prefixed small constant.

Let $z_1=[x,y]_\o$ and $z_2=\exp(E^u_{(y,\o)})\cap \exp(E^s_{(x,\o)})$. Note that, locally, $z_1-z_2$, $\pi^s_{y,x}(z_1-z_2)$ and $\pi^u_{y,x}(z_1-z_2)$ are well defined and satisfy the following:
$$d_M\left(z_1,\exp(E^u_{(y,\o)})\right)\ge\frac{|\pi^s_{y,x}(z_1-z_2)|}{\|\pi^s_{y,x}\|} \ge\frac{|\pi^s_{y,x}(z_1-z_2)|}{\mc P+\e'},$$
and
$$d_M\left(z_1,\exp(E^s_{(x,\o)})\right)\ge\frac{|\pi^u_{y,x}(z_1-z_2)|}{\|\pi^u_{y,x}\|} \ge\frac{|\pi^u_{y,x}(z_1-z_2)|}{\mc P+\e'}.$$
Since $|z_1-z_2|\le |\pi^u_{y,x}(z_1-z_2)|+|\pi^s_{y,x}(z_1-z_2)|$, without losing generality, we assume that $|\pi^s_{y,x}(z_1-z_2)|>\frac12|z_1-z_2|$. By i) of Lemma \ref{L:InvMani}, we have that
\begin{equation}\label{E:P2)1}
|z_1-z_2|\le 2(P+\e')d_M\left(z_1,\exp(E^u_{(y,\o)})\right)\le \frac{2(P+\e')}{10\mc P}d_M(y,z_1).
\end{equation}
Also note that
\begin{align*}
&d_M(y,z_2)\le \|\pi^s_{y,x}\|d_M(x,y)\le (\mc P+\e')d_M(x,y),\\
&d_M(y,z_1)\le d_M(y,z_2)+|z_1-z_2|.
\end{align*}
Then we have
$$|z_1-z_2|\le \frac{2(P+\e')^2}{10\mc P}d_M(x,y)+\frac{2(P+\e')}{10\mc P}|z_1-z_2|,$$
which implies that
\begin{equation}\label{E:P2)2}
|z_1-z_2|\le L_1'd_M(x,y), \text{ where }L_1'=\frac{\frac{2(P+\e')^2}{10\mc P}}{1-\frac{2(P+\e')}{10\mc P}}.
\end{equation}
Thus, by taking $L_1=L_1'+P+\e'$, we have that
$$d_M([x,y]_\o,y)\le |z_1-z_2|+d_M(y,z_2)\le L_1d(x,y).$$
The other part follows exactly from the same argument. Proof for P2) is completed.\\

Next, we sketch the proof for P3). Let $z_2=\exp(E^u_{(y,\o)})\cap \exp(E^s_{(x,\o)})$ which is the same as above and  $z_3=\exp(E^s_{(y,\o)})\cap \exp(E^u_{(x,\o)})$. Then we have
\begin{align}\begin{split}\label{E:P3)1}
&\left| |y-z_3|-|x-z_2|\right|\\
\le&\left||y-z_3|-|\pi^u_{(x,\o)}(x-y)|\right|+\left||x-z_2|-|\pi^u_{(x,\o)}(x-y)|\right|\\
=&\left||\pi^s_{x,y}(y-x)|-|\pi^u_{(x,\o)}(x-y)|\right|+\left||\pi^s_{y,x}(y-x)|-|\pi^u_{(x,\o)}(x-y)|\right|\\
\le &2\e'|x-y|,
\end{split}
\end{align}
where $\e'>0$ is a prefixed small number from (\ref{E:ErrProj}. To complete the proof of P3), we need to estimate the distance between $z_1=[x,y]_\o$ and $z_2$ (similarly the distance betwwen $z_4:=[y,x]_\o$ and $z_3$), which has been given by (\ref{E:P2)1}) and (\ref{E:P2)2}). The last step is to take $\d_1>0$ small enough such that the constant "$10\mc P$" in (\ref{E:P2)1}) and (\ref{E:P2)2}) can be replaced by a large enough constant $K>0$ so that $L'_1<\frac14\e_1$. This can be done, because of "$LipDg^\tau_{(x,\o)}<L$" from i) and iii) of Lemma \ref{L:InvMani}. Combining with (\ref{E:P3)1}), we have that  P3) holds when we take $\e'<\frac18\e_1$.


 We now are ready to prove the  lemma.

 By the setting of $\phi$, there exists $L_2\ge 1$ such that the Lipschitz constants of $f_\o$ and $f^{-1}_\o$ are less that $L_2$ for all $\o\in\O$.

Set $C_0:=\frac{(\e_1+L_1)(L_2+1)}{1-\e_1-(1+\e_1)e^{-\l}}$, where we require the smallness of $\e_1$ to ensure the positivity of $C_0$.
Now take $\a_0>0$ with $4C_0\a_0<\d_1$, we will show that there exists $C\ge 1$ satisfying the required properties.

We call a subset $R_\o\subset M_\o$ a {\em rectangle} if its diameter is small and $[x,y]_\o\in R_\o$ for all $x,y\in R_\o$. Let $A,B\subset M$ be Borel sets with diameter less than $2\e$ and their Hausdorff distance are smaller than $\d_1$. It is not hard to see that $[A,B]_\o:=\{[x,y]_\o|\ x\in A, y\in B\}$ is a rectangle in $M_\o$. For a rectangle $R_\o\subset M_\o$, define that
$$\partial ^sR_\o=\{x\in R_\o:x\notin int(W^u_{\e}(x,\o)\cap R_\o)\}$$
$$\partial ^uR_\o=\{x\in R_\o:x\notin int(W^s_{\e}(x,\o)\cap R_\o)\}$$

For a given $\a\in(0,\a_0)$, $x,y\in M_\o$ and $z\in M_{\t\o}$ satisfying that $d_M(x,y)<\a$ and $d_M(f_\o(x),z)<\a$,  we have the following:
\begin{itemize}
\item[p1)] $d_M(f_\o(y),z)\le d_M(f_\o(y),f_\o(x))+d_M(f_\o(x),z) \le (L_2+1)\a$.
\item[p2)]
$\max\{d_M([z,f_\o(y)]_{\t\o},z),d_M([f_\o(y),z]_{\t\o},y)\}\le L_1(L_2+1)\a$ by P2).
\item[p3)] For any $y'\in f_\o(W^s_{C_0\a}(y,\o))$, by P2), p1) and p2),
\begin{align*}
d_M([z,y']_{\t\o},z)
&\le d_M([z,y']_{\t\o},[z,f(y)]_{\t\o})+d_M([z,f(y)]_{\t\o},z)\\
&\le d_M(y',f_\o(y))+\e_1d_M(y',[z,f(y)]_{\t\o})+d_M([z,f(y)]_{\t\o},z)\\
&\le C_0\a e^{-\l}+\e_1(d_M(y',f_\o(y))+d_M(f_\o(y),z))+L_1(L_2+1)\a\\
&\le (1+\e_1)C_0\a e^{-\l}+(\e_1+L_1)(L_2+1)\a\\
&= (1-\e_1)C_0\a.
\end{align*}
\item[p4)] Similar to the arguments for p3), we have that for any $z'\in f_{\t\o}^{-1}(W^u_{C_0\a}(z,\t\o))$
$$ d_M([z',y]_\o,y)\le (1-\e_1)C_0\a.$$
\end{itemize}
Let $R_\o(y):=\overline {[W^u_{C_0\a}(y,\o),W^s_{C_0\a}(y,\o)]_\o}$. By P1) and P3), we have that for any $x\in R_\o(y)$
\begin{align}\begin{split}\label{E:DisEst1}
d_M(x,y)&\le d_M(x,[y,x]_\o)+d_M([y,x]_\o,y)\\
&\le L_1d_M([x,y]_\o,[y,x]_\o)+C_0\a\\
&\le (2L_1+1)C_0\a.
\end{split}\end{align}
 Next, define two complete metric spaces $\mc V^s_\o(y)$ and $\mc V^u_\o(y)$ as follows:
$$\mc V^\tau_\o(y):=\{W^\tau_{3C_0\a}(x,\o)\cap R_\o(y)|\ x\in R_\o(y)\},\ \tau=u,s,$$
where the metric in $\mc V^\tau_\o(y)$ is the Hausdorff distance.

Define a map $\mc F^u_\o(y,z):\mc V^u_\o(y)\to\mc V^u_{\t\o}(z)$ by letting
$$\mc \mc F^u_\o(y,z)W=f_\o(W)\cap R_{\t\o}(z),\text{ where }W\in \mc V^u_\o(y).$$
By the choice of $C_0$ and p3), we can see that the above map is well defined. Moreover, p3) also implies that $\mc \mc F^u_\o(y,z)$ is a contraction on $\mc V^u_\o(y)$ with contracting rate $1-\e_1$.
Similarly, we can define a contracting map $\mc F^s_{\t\o}(z,y):\mc V^s_{\t\o}(z)\to\mc V^s_{\o}(y)$ by p4).\\

 Now for a given $(\o,\a)$-pseudo orbit $\{(x_i,\t^i\o)\}_{i\in\mb Z}$ and $i,j\in\mb Z$ with $i<j$, define that
$$\mc F^u_{i,j}:=\mc F^u_{\t^{j-1}\o}(x_{j-1},x_{j})\circ\cdots\circ \mc F^u_{\t^i\o}(x_i,x_{i+1})(:\mc V^u_{\t^i\o}(x_i)\to \mc V^u_{\t^j\o}(x_j)),$$
$$\mc F^s_{j,i}:=\mc F^s_{\t^{i+1}\o}(x_{i+1},x_{i})\circ\cdots\circ \mc F^s_{\t^j\o}(x_j,x_{i-1})(:\mc V^s_{\t^j\o}(x_j)\to \mc V^s_{\t^i\o}(x_i)),$$
which are both contractions. Thus the following set consists of only one point which is denoted by $y_i$
  $$\lim_{k\to\infty}\left(\mc F^u_{i-k,i}\left(V^u_{\t^{i-k}\o}(x_{i-k})\right)\cap \mc F^s_{i+k,i}\left(V^s_{\t^{i+k}\o}(x_{i+k})\right)\right).$$
  It is clear that $\{(y_i,\t^i\o)\}_{i\in \mb Z}$ is a true orbit and $y_i\in R_{\t^i\o}(x_i)$. Therefore, (\ref{E:DisEst1}) tells that taking $C=(2L_1+1)C_0$ will fit the need.

\end{proof}

\begin{lem}\label{L:ContinousShadowing}
For a given $\b\in(0,\b_0)$, let $\a\in (0,\a_0)$ be as in   Lemma \ref{L:Shadowing} corresponding to $\b$ where $\a_0$ is as in Lemma  \ref{L:ClosingLemma}.  Then for any $\tau>0$, there exists $N_0\in\mb N$ such that given any  two $(\o,\a)$-pseudo orbits $\{(x_i,\t^i\o)\}_{i\in\mb Z}$ and $\{(x_i',\t^i\o)\}_{i\in\mb Z}$ which are $(\o,\b)$-shadowed by two uniquely defined true orbits $\{(y_i,\t^i\o)\}_{i\in\mb Z}$ and $\{(y_i',\t^i\o)\}_{i\in\mb Z}$ respectively, for any $N>N_0$ if $x_i=x_i',\ \forall i\in[-N,N]$, then $d_M(y_0,y_0')<\tau$.
\end{lem}
\begin{proof}
Note that if $x_i=x_i',\ \forall i\in[-N,N]$, then $d_M(y_i,y_i')<2\b,\ \forall i\in[-N,N]$.  Let $z=[y_0,y'_0]_\o$, where $[\cdot,\cdot]_\cdot$ is defined in \ref{L:LocCoor} for a fixed $\d$ corresponding to $\e=\frac12\e_0$. By the definition of local unstable manifold and ii) of Lemma \ref{L:InvMani}, we have that
\begin{align*}
d_M(z,y'_0)&\le e^{-N\l}d_M(\pi_M\phi^N(z,\o),y_N')\\
&\le e^{-N\l}(d_M(y_N,y'_N)+d_M(\pi_M\phi^N(z,\o),y_N))\\
&\le e^{-N\l}(2\b+d_M(z,y_0))\\
&\le e^{-N\l}(2\b+\frac12\e_0).
\end{align*}
The same argument is applicable to $d_M(z,y_0)$ if we reverse the time. Therefore,
$$d_M(y_0,y'_0)\le  2L_1e^{-N\l}(2\b+\frac12\e_0),$$
where $L_1$ is a constant given by P2) at beginning of the proof of Lemma  \ref{L:ClosingLemma}.
The proof is completed.
\end{proof}

Next, we consider the shadowing property of the map $\tilde \phi:L^\infty(\O, M)\to L^\infty(\O, M)$ induced by $\phi$. In this case, the system is a deterministic system, and the pseudo orbit is defined in the standard way:  for any $\a>0$, an orbit $\{g_i\}_{i\in \mb Z}\subset L^\infty(\O, M)$ is called an {\bf $\a$-pseudo orbit}  if for any $i\in \mb Z$
$$d_{L^\infty(\O, M)}(\tilde \phi(g_i),g_{i+1})<\a.$$
\begin{lem}\label{L:LShad}
For any $\b\in(0,\b_0)$, there exists an $\a>0$ such that for any $\a$-pseudo orbit $\{g'_i\}_{i\in \mb Z}$, there is a unique true orbit $\{g_i\}_{i\in \mb Z}$ of $\tilde\phi$ such that
$$\sup_{i\in \mb Z}d_{L^\infty(\O, M)}(g_i,g'_i)<\b,$$
where $\b_0$ is as in Lemma \ref{L:Shadowing}, and
$\{g_i\}_{i\in \mb Z}$ is called {\bf $\b$-shadowing} orbit of $\{g'_i\}_{i\in \mb Z}$.

Moreover,  there exist an $\a_0>0$ and $C\ge 1$ so that for any $\a\in (0,\a_0)$ and any $\a$-pseudo orbit $\{g'_i\}_{i\in\mb Z}$, there exists a true orbit $\{g_i\}_{i\in\mb Z}$ of $\tilde\phi$ satisfying the following
$$\sup_{i\in \mb Z} d_{L^\infty(\O, M)}(g_i,g'_i)<C\a.$$
\end{lem}
\begin{proof}
First, we note that each $\a$-pseudo orbit $\{g'_i\}_{i\in \mb Z}$ induces an $(\o,\a)$-pseudo orbit \[\{(g'_i(\t^i\o),\t^i\o)\}_{i\in \mb Z}\]
for each $\o\in \O$. Thus Lemma \ref{L:Shadowing} is applicable and produces a map $g:\O\to M$ such that $\{(g(\t^i\o),\t^i\o)\}_{i\in \mb Z}$ is $(\o,\b)$-shadowing $\{(g_i(\t^i\o),\t^i\o)\}_{i\in \mb Z}$, where $\a$ and $\b$ are as in Lemma \ref{L:Shadowing}. Once $g$ is Borel measurable, the proof follows from Lemma \ref{L:Shadowing} and \ref{L:ClosingLemma}.

Next, we prove that $g$ is Borel measurable. For $h\in L^\infty(\O, M)$ and $r>0$, denote $B_\O(h,r)$ the open tubular neighborhood of $graph(h)$, i.e.,
$$B_\O(h,r):=\left\{(x,\o)|\ \o\in\O, d_M(x,h(\o))<r\right\}.$$
It is not hard to see that $B_\O(h,r)$ is a Borel subset of $M\times\O$. In fact, since $M$ is compact, there exists a sequence of simple functions $\{h_i:\O\to M\}_{i\in \mb N}$ such that
$$d_{L^\infty(\O, M)}(h_i,h)<\frac{r}{i+2},\ i\in \mb N.$$
Note that each $B_\O(h_i,\frac{ir}{i+2})$ is a Borel subset of $M\times\O$ because it is a union of finite rectangle sets, thus $B_\O(h,r)=\cup_{i\in N}B_\O(h_i,\frac{ir}{i+2})$ implies the Borel measurability of  $B_\O(h,r)$.

Coming back to $g$ and $g_i'$, by the uniqueness of $g$, we have that
$$graph(g)=\bigcap_{i\in\mb Z}\phi^{-i}\left(B_\O(g_i,\b)\right).$$
For $k\in\mb N$, let $$S_k:=\bigcap_{i\in[-k,k]\cap \mb Z}\phi^{-i}\left(B_\O(g_i,\b)\right).$$
By the continuity of $\phi$ and $\phi^{-1}$, $S_k$ is a Borel subset in $M\times\O$. Also note that for each $\o\in\O$,
$$S_k\cap M_\o=\bigcap_{i\in[-k,k]\cap \mb Z}\pi_M\phi^{-i}\left(\left\{(x,\t^i\o)|\ d_M(x,g_i(\t^i\o))<\b\right\}\right),$$
which is a nonempty open subset of $M_\o$. By applying Theorem \ref{T:BMST} to each $S_k$, we have that for each $k\in\mb N$, there exists a Borel measurable map $g_k:\O\to M$ such that
$$graph(g_k)\subset S_k.$$
By the uniqueness of $g$ and Lemma \ref{L:Expansive}, we have that $$\lim_{k\to\infty}\sup_{\o\in\O}\{d_M(g_k(\o),g(\o)\}=0,$$
which implies that $g$ is Borel measurable. The proof is completed.
\end{proof}
\begin{rem}\label{R:Alpha0&C}
We remark here that $\a_0$, $\b_0$ and $C$ in Lemma \ref{L:Shadowing}, \ref{L:ClosingLemma}, and \ref{L:LShad} can take the same value. Thus we will not distinguish them when these Lemmas are used in the rest of the paper.
\end{rem}

The next lemma is a straightforward consequence of Lemma \ref{L:ContinousShadowing} and \ref{L:LShad}.
\begin{lem}\label{L:LContShad}
Let $\a_0$ and  $\b_0$ are as in Lemma  \ref{L:LShad}, and for a given $\b\in(0,\b_0)$, $\a\in (0,\a_0)$ be the one as in  Lemma \ref{L:LShad} corresponding to $\b$.   Then for any $\tau>0$, there exists $N_0\in\mb N$ such that the following holds: given any  two $\a$-pseudo orbits $\{g_i\}_{i\in\mb Z}$ and $\{g_i'\}_{i\in\mb Z}$ of $\tilde \phi$, which are $\b$-shadowed by two true orbits $\{\tilde g_i\}_{i\in\mb Z}$ and $\{\tilde g_i'\}_{i\in\mb Z}$ respectively, then for any $N>N_0$, $g_i=g_i',\ \forall i\in[-N,N]$, we have that
$d_{L^\infty(\O, M)}(\tilde g_0,\tilde g_0')<\tau.$
\end{lem}



\section{Density of Random Periodic Points}
In this section, we prove Theorem \ref{T:TheoryAnosovMix} based on the mixing property and the shadowing property.

\vskip0.05in
\noindent {\it Proof of Theorem \ref{T:TheoryAnosovMix}}. Firstly, fix a given Borel function $g\in L^\infty(\O, M)$ and $\e>0$. Without losing  generality, we assume that $\e<\frac12\b_0$, where $\b_0$ is as in Lemma \ref{L:Shadowing} and \ref{L:LShad}. Taking $\b=\frac13\e$ and applying Lemma \ref{L:LShad}, there exists $\a>0$ corresponding to $\b$ as in Lemma \ref{L:LShad}. Again, for the sake of convenience, we assume that $\a<\b$.
By the compactness of $M$, there exists  a finite open cover $\{U_i,1\le i\le n_0\}$ of which each $U_i=B(x_i,\frac13\a)$ is a ball of radius $\a$ in $M$.

Secondly, define a function $k:\O\to \{1,2,\cdots, n_0\}$ by the following
\begin{equation}\label{E:MeasurableIndex}
k(\o)=\min\left\{i\in \{1,\cdots,n_0\}\big|\ g(\o)\in U_i\right\},\ \forall \o\in\O.
\end{equation}
It is not hard to see that $k$ is a measurable  function since $g$ is measurable and the following holds
\begin{align*}
&\left\{\o|\ k(\o)=i\right\}\\
=&\begin{cases}
\{\o|\ g(\o)\in U_i\}&\text{ when } i=1\\
\{\o|\ g(\o)\in U_i\}\setminus\left(\bigcup_{1\le j\le i-1}\{\o|\ g(\o)\in U_j\}\right)&\text{ when } i\in\{2,\cdots, n_0\}
\end{cases}.
\end{align*}
This Borel measurable function $k(\o)$ induces a finite measurable partition of $\O$,
$$\left\{\O_i:=\{\o|\ k(\o)=i\}\right\}_{i\in\{1,\cdots,n_0\}}.$$
Note that, by condition H2), there exists $N\in\mb N$ such that
\begin{equation}\label{E:MixingTimeS2}
\pi_M\phi^n(U_i,\o)\cap U_j\neq \emptyset,\ \forall \o\in \O, i,j\in\{1,\cdots,n_0\},n\ge N.
\end{equation}
Then, for any $n\ge N$, $k$ and $n$ induces a finite measurable partition of $\O$, $\{\O^{n}_{ij}\}_{i,j\in\{1,\cdots,n_0\}}$, by letting
$$\O^n_{ij}:=\left\{\o\big|\ k(\o)=i, k(\t^n\o)=j\right\}=\t^{-n}(\O_j)\cap \O_i.$$
By the uniform continuity of $\phi$, there exists $\d_1>0$ such that the following holds
\begin{equation}\label{E:CondDelt1}
d_M(\pi_M\phi(x,\o_1),\pi_M\phi(x,\o_2))<\frac16\a,\ \forall x\in M\text{ and } d_\O(\o_1,\o_2)<\d_1.
\end{equation}
Additionally, by the uniform continuity of $\t$, for a fixed $n\ge N$, there exists $\d_2>0$ such that
\begin{equation}\label{E:CondDelt2}
d_\O(\t^m\o_1,\t^m\o_2)<\d_1\ \forall0\le m\le  n-1 \text{ and  }d_\O(w_1,w_2)<\d_2.
\end{equation}
Now take a finite measurable partition of $\O$, $\{\O'_s\}_{s\in S}$, which is taken to be a refinement of $\{\O^n_{ij}\}_{i,j\in\{1,\cdots,n_0\}}$ and has diameter less than $\d_2$. By (\ref{E:MixingTimeS2}), we have that for any $\o\in\O^n_{ij}$, there exists an $x\in U_i$ such that  $\phi^n(x,\o)\in (U_j,t^n\o)$.
Then, for each $\O'_s(\subset \O^n_{ij})$, arbitrarily fix a point $(y_s,\o_s)\in U_i\times \O'_s$ such that
 $$\pi_M\phi^n(y_s,\o_s)\in U_j.$$
Define a simple function $y=\sum_{s\in S}y_s\chi_{\O'_s}$, where $\chi_{\O'_s}$ is the characteristic function of $\O'_s$. \\

Now, we are ready to construct the $(\o,\a)$-pseudo-orbits which visit the tubular neighborhood of $graph(g)$ periodically (with period of $n$).  For a given $\o\in\O$, denote $q_\o:\mb Z\to S$ the function satisfying that for any $l\in \mb Z$, $\t^{ln}(\o)\in \O_{q(l)}'$. Then define
\begin{equation}\label{E:DefPseOrbS2)}
y'_i(\o)=
\begin{cases}
y(\t^{ln}\o),&\text{when } i=ln\\
\pi_M\left(\phi^{i-ln}(y(\t^{ln}\o),\o_{q_\o(l)})\right),&\text{when }i\in[ln+1,(l+1)n-1]
\end{cases},
\end{equation}
where the map $y'_i:\O\to M$ is clearly Borel measurable.
For $l\in\mb Z$ and $i\in  [ln,(l+1)n-1]$, by the choice of $\d_2$ and (\ref{E:CondDelt2}), we have that
$$d_\O(\t^{i}\o,\t^{i-ln}\o_{q_\o(l)})<\d_1.$$
Therefore, by (\ref{E:CondDelt1}), we have that, for $i\in  [ln,(l+1)n-2]$
\begin{align*}
&d_M\left(\pi_M\phi(y_i'(\o)),y_{i+1}'(\o)\right)\\
=&d_M\left(\pi_M\phi\left(\pi_M\phi^{i-ln}(y(\t^{ln}\o),\o_{q_\o(l)}),\t^i\o\right),\pi_M\phi\left(\pi_M\phi^{i-ln}(y(\t^{ln}\o),\o_{q_\o(l)}),\t^i\o_{q_\o(l)}\right)\right)\\
\le&\frac16\a,
\end{align*}
and for $i=(l+1)n-1$,
\begin{align*}
&d_M\left(\pi_M\phi(y'_{(l+1)n-1}(\o),\t^{(l+1)n-1}\o),y'_{(l+1)n}(\o)\right)\\
\le &d_M\left(\pi_M\phi\left(\pi_M\phi^{i-ln}(y(\t^{ln}\o),\o_{q_\o(l)}),\t^i\o\right),\pi_M\phi\left(\pi_M\phi^{i-ln}(y(\t^{ln}\o),\o_{q_\o(l)}),\t^i\o_{q_\o(l)}\right)\right)\\
&+d_M\left(\pi_M\phi^{N}(y(\t^{ln}\o),\o_{q_\o(l)}),y(\t^{(l+1)n}\o)\right)\\
<&\frac16\a+\frac23\a<\a.
\end{align*}
In the above estimate, we used the fact that $\phi^{n}(y(\t^{ln}\o),\o_{q_\o(l)})$ and $y(\t^{(l+1)n}\o)$  fall in a same element of  $\{U_i\}_{i\in\{1,\cdots,n_0\}}$. This conclude that $\{(y'_i\in L^\infty(\O, M)\}_{i\in\mb Z}$ forms an $\a$-pseudo orbit of $\tilde\phi$. Thus, by Lemma \ref{L:LShad}, there exists a unique true orbit $\{\tilde\phi^i(\tilde g), \tilde g\in L^\infty(\O, M)\}_{i\in\mb Z}$ which is $\b$-shadowing  $\{y'_i\}_{i\in\mb Z}$. Note that $\{\phi^i(\tilde g)\}_{i\in\mb Z}$ and $\{\phi^{i+n}(\tilde g)\}_{i\in\mb Z}$ are two true orbits and
$$d_{L^\infty(\O, M)}(\tilde\phi^{i+n}(\tilde g),\tilde\phi^i(\tilde g))\le 2\b<\frac13\b_0,\ \forall i\in \mb Z.$$
Then Lemma \ref{L:LContShad} implies that
$$\tilde\phi^n(\tilde g)=\tilde g\ i.e.\ \phi^{n} (graph(\tilde g))=graph(\tilde g).$$
This completes the proof of this theorem.  \qed

\medskip

From the above proof, we summarize a technical lemma as follows, which can be used to join two points in $L^\infty(\O, M)$ by using a segment of an orbit of $\tilde \phi$. Since the proof follows exactly the same idea as above, we omit it here.
\begin{lem}\label{L:JoinSegments}
For any $\e>0$, there exists a $N\in \mb N$ such that for any $g_1,g_2\in L^\infty(\O, M)$ and $n\in [N,\infty)\cap \mb N$, there exists $g\in L^\infty(\O, M)$ satisfies the following
$$\max\left\{d_{L^\infty(\O, M)}(g,g_1),d_{L^\infty(\O, M)}(\tilde\phi^n(g),g_2)\right\}<\e.$$
\end{lem}

\section{Random Horseshoe}

In this section, we first show that $\underline h(M\times\O|\O)>0$. Then, we prove Theorem \ref{T:TheoryAnosovMix2}, i.e., the existence of strong full random  horseshoes and full random and  horseshoes.

Before going to the proof, we need to introduce an alternative definition of $\underline h(M\times \O|\O)$ which is more convenient.
\begin{lem}\label{L:AlDefH}
Suppose  that $\phi$ satisfy Condition H1). Then the following holds:
$$\underline h(M\times \O|\O)=\lim_{\e\to0^+}\inf_{\o\in\O}\limsup_{n\to\infty}\frac1n\log N(\o,\e,n),$$ where
\begin{align*}
&N(\o,\e,n):=\\&\max\  \left\{card(E)\big|\ E\subset M_{\o}, E \text{ is an }\e\text{ separated set with respect to }d^{\phi}_{M\times \O,n}\right\}.
\end{align*}
\end{lem}
\begin{proof}
Recall that the original definition of $\underline h(M\times \O|\O)$ is given by
$$\underline h(M\times \O|\O)=\inf_{\o\in\O}\lim_{\e\to0^+}\limsup_{n\to\infty}\frac1n\log N(\o,\e,n).$$
This lemma simply states that the above order of $\inf$ and $\lim_{\e\to0^+}$ can be exchanged. First, we clearly have that
$$\underline h(M\times \O|\O)=\inf_{\o\in\O}\lim_{\e\to0^+}\limsup_{n\to\infty}\frac1n\log N(\o,\e,n)\ge\lim_{\e\to0^+}\inf_{\o\in\O}\limsup_{n\to\infty}\frac1n\log N(\o,\e,n).$$
Next, we show that the opposite inequality also holds. To see this, we claim

\noindent
{\bf Claim}:{\it  There exists $\eta >0$ such that
\begin{equation}\label{E:ConstantEta}
h_{top}(\phi|M_\o):=\lim_{\e\to0^+}\limsup_{n\to\infty}\frac1n\log N(\o,\e,n)=\limsup_{n\to\infty}\frac1n\log N(\o,\eta,n),\ \forall \o\in\O.
\end{equation}}
\noindent{\em Proof of the Claim}: It is obvious that $h_{top}(\phi|M_\o)\ge \limsup_{n\to\infty}\frac1n\log N(\o,\eta,n)$ by definition. Hence, it remains to prove the opposite inequality.

Firstly, we fix an $\eta>0$ satisfying the following property: for any $x,y\in M_\o$
$$\sup_{n\in\mb Z}d(\phi^n(x,\o),\phi^n(y,\o))\le \eta\;\text{ implies }\;x=y.$$
The existence of such $\eta$ follows from the expansive property of $\phi$, where one only requires $\eta<\e$ for the $\e$ from Lemma \ref{L:Expansive}.

Secondly, by Lemma \ref{L:ContinousShadowing}, we have that for any $\e'\in(0,\eta)$, there exists $L(\e')\in \mb N$ such that for any $x,y\in M_\o$
\begin{equation}\label{E:SepControl}
\max_{|n|\le L(\e')}d(\phi^n(x,\o),\phi^n(y,\o))\le \eta\text{ implies }d_M(x,y)<\e'.
\end{equation}
By (\ref{E:SepControl}), we have for any $\o\in\O$, $n\in\mb N$, and $x,y\in M_\o$,
$$d^\phi_{M\times \O,n+2L(\e')}((x,\o),(y,\o))\le \eta\text{ implies }d^\phi_{M\times \O,n}(\phi^{L(\e')}(x,\o),\phi^{L(\e')}(y,\o))<\e',$$
which yields
\begin{equation}\label{E:SepNEst1}
N(\t^{L(\e')}\o,\e',n)\le N(\o,\eta,n+2L(\e')),\ \forall \o\in\O, n\in \mb N, \e'\in(0,\eta).
\end{equation}
Note that
\begin{equation}\label{E:SepNEst2}
N(\o,\e',n+m)\le N(\o,\e',m)N(\t^m\o,\e',n).
\end{equation}
(\ref{E:SepNEst1}) and (\ref{E:SepNEst2}) imply that
\begin{align*}
h_{top}(\phi|M_\o)=&\lim_{\e'\to0^+}\limsup_{n\to\infty}\frac1n\log N(\o,\e',n)\\
=&\lim_{\e'\to0^+}\limsup_{n\to\infty}\frac1{n+L(\e')}\log N(\o,\e',n+L(\e'))\\
\le&\lim_{\e'\to0^+}\limsup_{n\to\infty}\frac1{n+L(\e')}\left(\log N(\o,\e',L(\e'))+\log N(\t^{L(\e')}\o,\e',n)\right)\\
\le&\lim_{\e'\to0^+}\limsup_{n\to\infty}\frac1{n+L(\e')}\left(\log N(\o,\e',L(\e'))+\log N(\o,\eta,n+2L(\e'))\right)\\
=&\lim_{\e'\to0^+}\limsup_{n\to\infty}\frac1{n+L(\e')}\log N(\o,\eta,n+2L(\e'))\\
=&\lim_{\e'\to0^+}\limsup_{n\to\infty}\frac1{n}\log N(\o,\eta,n)\\
=&\limsup_{n\to\infty}\frac1{n}\log N(\o,\eta,n).
\end{align*}
Therfore (\ref{E:ConstantEta}) holds, which yields
\begin{align*}
\underline h(M\times \O|\O)=&\inf_{\o\in\O}\lim_{\e\to0^+}\limsup_{n\to\infty}\frac1n\log N(\o,\e,n)\\
=&\inf_{\o\in \O} h_{top}(\phi|M_\o)\\
=&\inf_{\o\in\O}\limsup_{n\to\infty}\frac1n\log N(\o,\eta,n)\\
\le &\lim_{\eta\to0^+}\inf_{\o\in\O}\limsup_{n\to\infty}\frac1n\log N(\o,\eta,n),
\end{align*}
which completes the proof of this lemma.
\end{proof}

Next, we will  show that $\underline h(M\times \O|\O)>0$. It is not hard to see that once $\phi^N$ has a strong horseshoe with $k$ symbols for some $N,k,\e$ as in the Definition \ref{D:SRandomHorseshoe}, then by the definition of $\underline h(M\times\O|\O)$, we have that
$$\underline h(M\times\O|\O)\ge \frac{\log k}N>0.$$
Thus, the first assertion on the positivity of $\underline h(M\times\O|\O)$ in i) of Part B of Theorem \ref{T:TheoryAnosovMix} follows from the lemma below:
\begin{lem}\label{L:Horseshoe}
If $\phi$ satisfies condition H1) and H2), then for any $k\in \mb N\setminus\{1\}$ there exists $N\in \mb N$ such that $\phi^N$ has a $k$-symbol strong full random horseshoe.
\end{lem}
\begin{proof}

Let $\{x_i\}_{1\le i\le k}$ be an $3\e$-separated finite subset of $M$. Such a set exists for small enough $\e>0$, and we fix this $\e$.  Without loss of  generality, we assume that $\e<\frac12\b_0$, where $\b_0$ is as in Lemma \ref{L:LShad}. For $\b=\frac13\e$,  there exists an $\a>0$ corresponding to $\b$ as in Lemma \ref{L:LShad}. Again, for the sake  of convenience, we assume that $\a<\b$. For any $i\in \{1,\cdots,k\}$, let $U_i=B(x_i,\frac13\a)$.

By Condition H2), there exists an $N\in\mb N$ such that the following holds
\begin{equation}\label{E:MixingTimeS2S}
\pi_M\phi^n(U_i,\o)\cap U_j\neq \emptyset,\ \forall \o\in \O, i,j\in\{1,\cdots,k\},n\ge N.
\end{equation}
By the uniform continuity of $\phi$ and $\t$, there exist $\d_1,\d_2>0$ such that (\ref{E:CondDelt1}) and (\ref{E:CondDelt2}) hold. Let $\{\O_s'\}_{s\in S}$ be a finite Borel measurable partition of $\O$ with diameter less than $\d_2$. By (\ref{E:MixingTimeS2S}), for any $s\in S$ and $i,j\in\{1\cdots,k\}$, there exists $(y^{i,j}_s,\o^{i,j}_s)\in U_i\times \O'_s$ such that
$$\pi_M\phi^N(y^{i,j}_s,\o^{i,j}_s)\in U_j.$$
We define a simple function $y^{i,j}=\sum_{s\in S}y^{i,j}_s\chi_{\O'_s}$ for all $i,j\in\{1,\cdots,k\}$.

Now, we are ready to build the strong random horseshoe by constructing a proper $\a$-pseudo-orbits of $\tilde\phi$ associated to elements from $\mc S_k:=\{1,\cdots, k\}^{\mb Z}$.

For a given $\hat a=(a_i)_{i\in\mb Z}\in \mc S_k$, and a given $\o\in\O$, we denote $q_\o:\mb Z\to S$ the function satisfying that for any $l\in \mb Z$, $\t^{lN}(\o)\in \O_{q(l)}'$. Then, we define
\begin{equation}\label{E:DefPseOrbS2)}
y'_{\hat a,\o}(i)=
\begin{cases}
y^{a_l,a_{l+1}}(\t^{lN}\o),&\text{when } i=lN\\
\pi_M\left(\phi^{i-lN}(y^{a_l,a_{l+1}}(\t^{lN}\o),\o^{a_l,a_{l+1}}_{q_\o(l)})\right),&\text{when }i\in[lN+1,(l+1)N-1]
\end{cases}.
\end{equation}
Set $$g_i(\o):=y'_{\hat a,\o}(i),\ \forall \o\in\O,\ i\in \mb Z.$$
It is clear that each $g_i:\O\to M$ is Borel measurable since $g_i$ is either a simple function or the image of a simple function under $\phi$ or $\phi^{-1}$'s iterations.

For $l\in\mb Z$ and $i\in  [lN,(l+1)N-1]$, by the choice of $\d_2$ and (\ref{E:CondDelt2}), we have that
$$d_\O(\t^{i}\o,\t^{i-lN}\o^{a_l,a_{l+1}}_{q_\o(l)})<\d_1.$$
Therefore, by (\ref{E:CondDelt1}), we have that, for $i\in  [lN,(l+1)N-2]$, \\
(to save space, let $z=\pi_M\phi^{i-lN}\left(y^{a_l,a_{l+1}}(\t^{lN}\o),\o^{a_l,a_{l+1}}_{q_\o(l)}\right)$)
\begin{align*}
&d_M\left(\pi_M\phi(y'_{\hat a,\o}(i)),y'_{\hat a,\o}(i+1)\right)
=d_M\left(\pi_M\phi\left(z,\t^i\o\right),\pi_M\phi\left(z,\t^i\o^{a_l,a_{l+1}}_{q_\o(l)}\right)\right)
\le\frac16\a,
\end{align*}
while for $i=(l+1)N-1$,
\begin{align*}
&d_M\left(\pi_M\phi(y'_{\hat a,\o}(i),\t^{i}\o),y'_{\hat a,\o}(i+1)\right)\\
\le &d_M\left(\pi_M\phi\left(z,\t^i\o\right),\pi_M\phi\left(z,\t^i\o^{a_l,a_{l+1}}_{q_\o(l)}\right)\right)\\
&+d_M\left(\pi_M\phi^{N}(y^{a_l,a_{l+1}}(\t^{lN}\o),\o^{a_l,a_{l+1}}_{q_\o(l)}),y^{a_{l+1},a_{l+2}}(\t^{(l+1)N}\o)\right)\\
<&\frac16\a+\frac23\a<\a.
\end{align*}
 This concludes that for any $\o\in \O$, $\{(y'_{\hat,\o}(i),\t^i\o)\}_{i\in\mb Z}$ forms an $(\o,\a)$-pseudo orbit of $\phi$. Thus $\{g_i\}_{i\in\mb Z}$ forms an $\a$-pseudo orbit of $\tilde \phi$. Thereafter, Lemma \ref{L:LShad} implies the existence of a unique true orbit $\{\phi^i(\tilde )\}_{i\in\mb Z}$ which is $\b$-shadowing  $\{g_i\}_{i\in\mb Z}$.

 We define the map $\Psi:\mc S_k\to L^\infty(\O, M)$ by letting
$$\Psi(\hat a)=g_{\hat a},\ \forall \hat a\in \mc S_k.$$
 where the continuity of $\Psi$ follows from Lemma \ref{L:LContShad}.

 Note that $\{\tilde\phi^i(\tilde g_{\hat a})\}_{i\in\mb Z}$ and $\{\tilde\phi^i(\tilde g_{\s \hat a})\}_{i\in\mb Z}$ are two true orbits and
$$d_{L^\infty(\O, M)}(\tilde\phi^{i+N}(\tilde g),\tilde\phi^i(\tilde g_{\s\hat a})\le 2\b<\frac13\b_0,$$
therefore, by Lemma \ref{L:LContShad}, we have that
$$\tilde\phi^N(\tilde g_{\hat a})=(\tilde g_{\s\hat a} )\ i.e.\ \phi^{N} graph(\tilde g_{\hat a})=graph(\tilde g_{\s\hat a}),$$
which implies (ii) in Definition \ref{D:SRandomHorseshoe}.

 Now we prove i) in Definition \ref{D:SRandomHorseshoe}. For any $\hat a_1=(a_1(i))_{i\in \mb Z},\hat a_2=(a_2(i))_{i\in \mb Z}\in\mc S_k$ with $\hat a_1\neq \hat a_2$, let $s=\min\{|i||\ a_1(i)\neq a_2(i)\}$. Then,  by the choice of $x_i,U_i,\e,\a,\b$ at the beginning of this section,  we have that
 \begin{equation}\label{E:SepSym}
 d^{\phi}_{M\times\O,N}\left(\phi^{rN}(\Psi(\hat a_1)(\o),\o),\phi^{rN}(\Psi(\hat a_2)(\o),\o)\right)>\e,\ \forall\o\in\O,
 \end{equation}
 where $a_1(r)\neq a_2(r)$ with $|r|=s$ and $d^{\phi}_{M\times\O,N}$ is the Bowen metric. Note that there is a constant $L>1$ such that for any $\o\in\O$ and $x,y\in M$ with $d(x,y)>0$, we have
 $$d(\phi^{\pm1}(x,\o),\phi^{\pm1}(y,\o))\le Ld_M(x,y).$$
 Therefore
 \begin{equation}\label{E:SHorsSeparation}
 d_M\left(\pi_M\left(\Psi(\hat a_1)(\o),\o\right),\pi_M\phi^{rN}\left(\Psi(\hat a_2)(\o),\o\right)\right)>\e L^{-N},\ \forall\o\in\O,
\end{equation}
which implies Condition i) of Definition \ref{D:SRandomHorseshoe} by taking $\Delta =\e L^{-N}$.

 Finally, the continuity of $\Psi^{-1}$ can be obtained by the following estimate:
 \begin{align}\begin{split}\label{E:PsiInverseContinuity}
& d_M(\Psi(\hat a_1)(\o),\Psi(\hat a_2)(\o))\\
=&d\left((\Psi(\hat a_1)(\o),\o),(\Psi(\hat a_2)(\o),\o)\right)\\
\ge& L^{-|rN|}d\left(\phi^{rN}(\Psi(\hat a_1)(\o),\o),\phi^{rN}(\Psi(\hat a_2)(\o),\o)\right)\\
\ge& \e L^{-(s+1)N}.
\end{split}
 \end{align}
 The proof is completed.
\end{proof}

\begin{rem}\label{R:PositivLHVolEst}
Actually, the positivity of $\underline h(M\times\O|\O)$ can be derived more directly  by estimating the volume expanding rate on local unstable manifolds, which follows from the hyperbolicity of the systems on fibers. Nevertheless, we prefer to present the current proof  which is demonstrative for the ideas on constructing random horseshoes by mixing property.
\end{rem}

Next, we  construct a strong horseshoe which can capture enough entropy which is close to the lower topological  fiber entropy. The construction of the desired horseshoe follows from the same procedure as the one used in the proof of Lemma \ref {L:Horseshoe}.

\begin{proof}[Proof of Part A. of Theorem \ref{T:TheoryAnosovMix2}]
By Lemma \ref{L:AlDefH}, we have that there exists $\e'_0>0$ such that for any $\e'\in(0,\e'_0)$ the following holds
\begin{equation}\label{E:SepNEst}
\inf_{\o\in\O}\limsup_{n\to\infty}\frac1n\log N(\o,\e',n)>\underline h(M\times\O|\O)-\frac14\g.
\end{equation}

Let $\e'_1=\min\left\{\frac12\b_0,\e'_0\right\}$, where $\b_0$ is as in Lemma \ref{L:LShad}. Taking $\b=\frac16\e'_1$, there exists an $\a>0$ corresponding to $\b$ as in Lemma \ref{L:Shadowing}. Again, for the sake  of convenience, we assume that $\a<\b$. Now we fix an $\e\in(0,\b)$, and also fix a set of open balls in $M$ with radius equaling $\min\{\frac13\a,\frac16\e\}$, $\{U_i,0\le i\le p\}$, where $\{U_i\}_{1\le i\le p}$ forms an open cover of $M$.

By Condition H2), there exists $N_1\in\mb N$ such that the following holds
\begin{equation}\label{E:MixingLS}
\pi_M\phi^n(U_i,\o)\cap U_j\neq \emptyset,\ \forall \o\in \O, i,j\in\{0,1,\cdots,p\},n\ge N_1.
\end{equation}

By (\ref{E:SepNEst}) and $\b\le \frac16 \e_0'$, for any $\o\in\O$, there exist $l(\o),m(\o)\in \mb N$, and a $3\b$-separated set $\{z_i(\o)\}_{1\le i\le l(\o)}$ with respect to $d^\phi_{M\times \O,m(\o)}$ satisfies that
\begin{equation}\label{E:SepNEstO}
\frac{1}{m(\o)+2N_1}\log l(\o)>\underline h(M\times\O|\O)-\frac12\g.
\end{equation}
Since $\phi$ is continuous and $\O$ is compact, the $3\b$-separated set can be chosen properly to make $m(\o)$ being uniformly bounded on $\O$. The bound is denoted by $N_2$. Let $N=L(N_2+2N_1)$, where $L$ is a prefixed large positive integer.

Note that, by the choice of $\a$ and $\b$, if two $(\o,\a)$-pseudo orbits have a $3\b$-separation on time $n$, then the corresponding true orbits which are $(\o,\b)$-shadowing them respectively  have a $\b$-separation on time $n$. Therefore, here instead of constructing specific pseudo orbits, we would like to count the number of  $(\o,\frac23\a)$-pseudo orbits which have $3\b$-separations happened during time interval $[0,N]$.

 First, we describe the procedure defining the pseudo orbits. For a given $\o\in\O$, by (\ref{E:MixingLS}), we can use a segment of a true orbit to connect $(U_0,\o)$ to any one of $\{(U_i,\t^{N_1})\}_{1\le i\le p}$. This orbit segment of connecting $U_0$ to some $U_i$ is the first piece of the $(\o,\frac23\a)$-pseudo orbit. The second piece of the  $(\o,\frac23\a)$-pseudo orbit is one of the orbit segments \[\{\phi^j(z_i(\t^{N_1}\o),\t^{N_1}\o)\}_{0\le j\le m(\t^{N_1}\o)-1}
  \]which is a segment of a true orbit. So in order to ensure that these two pieces can be merged with variation $\le \frac23\a$, $U_i$ chosen in the first step should contains the $z_i(\t^{N_1}\o)$ chosen in the second step. The third piece of these $(\o,\frac23\a)$-pseudo orbits is to connect $\phi^{m(\t^{N_1}\o)}(z_i(\t^{N_1}\o),\t^{N_1}\o)$ to $U_0$ by a true orbit segment with length $N_1$, which can be done again by (\ref{E:MixingLS}). We repeat this procedure and  derive a sequence of integers in the following way
 \begin{align*}
 R_0&=0\\
 R_1&=2N_1+m(\t^{N_1}\o)\\
R_2&=2N_1+m(\t^{R_1+N_1}\o)\\
&\cdots\\
R_{n+1}&=2N_1+m(\t^{\sum_{i=1}^nR_i+N_1}\o)\\
&\cdots
 \end{align*}
Let $$s:=\max\left\{n\Big|\ \sum_{i=1}^nR_i\le (L-1)(N_2+2N_1)\right\}.$$
 The above procedure of constructing pseudo orbits repeats $s$ times in the time interval $[0,\sum_{i=1}^sR_i]$. For the time interval $[\sum_{i=1}^sR_i,N]$, we use only one orbit segment to connect $(U_0,\t^{\sum_{i=1}^sR_i}\o)$ to $(U_0,\t^{N}\o)$. It is not hard to see that the number of $(\o,\frac23\a)$-pseudo orbits which are $3\b$-separated in the time interval $[0,N]$ is
 $$K(\o)=\prod_{i=1}^sl\left(\t^{\sum_{j=0}^{i-1}R_j+N_1}\o\right).$$
 Let $N_3=N-\sum_{i=1}^sR_i$, then $N_3\le 2(N_2+2N_1)$. Therefore, by (\ref{E:SepNEstO}), we have that
 \begin{align*}
 &\frac{\log K(\o)}N\\
 =&\frac{\sum_{i=1}^s\log l\left(\t^{\sum_{j=0}^{i-1}R_j+N_1}\o\right)}N\\
 \ge&\frac{\left(\underline h(M\times\O|\O)-\frac12\g\right)\left(\sum_{i=1}^s\left(m\left(\t^{\sum_{j=0}^{i-1}R_j+N_1}\o\right)+2N_1\right)\right)}N\\
 =&\left(\underline h(M\times\O|\O)-\frac12\g\right)\frac{{\sum_{i=1}^sR_i}}N=\left(\underline h(M\times\O|\O)-\frac12\g\right)\frac{N-N_3}N\\
  \ge&\left(\underline h(M\times\O|\O)-\frac12\g\right)\left(1-\frac2L\right).
 \end{align*}
The largeness of $L$ is to ensure the  following
$$\frac{\log K(\o)}N>\underline h(M\times\O|\O)-\g.$$
For any given $\o\in \O$, let $\{y_j(\t^i\o)\}_{0\le i\le N-1}$, $1\le j\le K(\o)$, be the segments of $(\o,\frac23\a)$-pseudo orbits of $\phi$ defined above. By the uniform continuous of $\phi$, there exist $\d_1,\d_2>$ such that (\ref{E:CondDelt1}) and (\ref{E:CondDelt2}) hold, where $n$ in (\ref{E:CondDelt2}) is replaced by $N$. Let $\{\O_s'\}_{s\in S}$ be a finite Borel measurable partition of $\O$ with diameter less than $\d_2$. We fix an $\o'_s\in \O'_s$ for each $s\in S$, and let $K=\min\{K(\o'_s)\}_{s\in S}$. Note that we
$$\frac{\log K}N>\underline h(M\times\O|\O)-\g.$$
Define the simple functions as follows
$$g_{i,j}=\sum_{s\in S}y_j(\t^i\o_s')\chi_{\t^i(\O'_s)},\ 0\le i\le N-1,\ 1\le j\le K.$$
where $\chi_{\t^i(\O'_s)}$ is the characteristic function of $\t^i(\O'_s)$. It is not hard to see that (\ref{E:CondDelt1}) and (\ref{E:CondDelt2}) imply that
\begin{equation}\label{E:SimpFuncShd1}
d_{L^\infty(\O, M)}(\tilde\phi(g_{i,j}),g_{i+1,j})\le d_M(y_j(\t^i\o'_s),y_j(\t^{i+1}\o'_s)+\frac16\a<\a,\ 0\le i\le N-2,
\end{equation}
and
\begin{equation}\label{E:SimpFuncShd2}
d_{L^\infty(\O, M)}(\tilde\phi(g_{N-1,j}),g_{0,j'})\le \text{ diameter of }U_0+\frac16\a<\a,\ 1\le j,j'\le K.
\end{equation}
Also note that for $1\le j,j'\le K$ with $j\neq j'$ and each $\o\in\O$, $\{y_j(\t^i\o)\}_{0\le i\le N-1}$ and $\{y_{j'}(\t^i\o)\}_{0\le i\le N-1}$ are $3\b$-separated with respect to $d^\phi_{M\times \O,m(\o)}$, thus we have that
\begin{equation}\label{E:SHorsSepa}
\inf_{\o\in\O}\max_{0\le i\le N-1}d_M(g_{i,j}(\o), g_{i,j'}(\o))\ge 3\b.
\end{equation}

Now, we are ready to construct the strong full random horseshoe.

For $\hat a=(a_k)_{k\in\mb Z}\in \mc S_k$, let
\begin{equation}\label{E:SHorsPsOr}
 g^{\hat a}_i=g_{i-\left[\frac iN\right]N,a_{\left[\frac iN\right]}},\ i\in \mb Z,
\end{equation}
where $\left[\frac iN\right]$ is the largest integer less than $\frac iN$. By (\ref{E:SimpFuncShd1}) and (\ref{E:SimpFuncShd2}), we have that $\{ g^{\hat a}_i\}_{1\in\mb Z}$ forms an $\a$-pseudo orbit of $\tilde \phi$, thus Lemma \ref{L:LShad} implies a true orbit of $\tilde \phi$, $ \{\tilde\phi^i(g^{\hat a})\}_{i\in\mb Z}$, which is $\b$-shadowing $\{ g^{\hat a}_i\}_{1\in\mb Z}$.

Define $\Psi:\mc S_k\to L^\infty(\O, M)$ by letting
$$\Psi(\hat a)=g^{\hat a},$$where the continuity of $\Psi$ follows from  Lemma \ref{L:LContShad} immediately.
It is obvious that (\ref{E:SHorsSepa}) and the argument for (\ref{E:SHorsSeparation}) imply Condition i) of Definition \ref{D:SRandomHorseshoe} with $\Delta=\b L^{-N}$. ii) of Definition and the continuity of $\Psi^{-1}$ follow exactly from the same arguments as (\ref{E:PsiInverseContinuity}), thus we do not repeat it here.
This completes the proof of Part A. of Theorem \ref{T:TheoryAnosovMix2}. \end{proof}


Next, we  construct a full random horseshoe which can capture enough entropy that is close to the upper topological fiber entropy.
\begin{proof}[Proof of Part B. of Theorem \ref{T:TheoryAnosovMix2}]
For a given $\g>0$, by the definition of $\overline h(M\times \O|\O)$, we have that there exist $\o_0\in\O$ and $\e_0>0$ such that for any $0<\e<\e_0$, there exist infinitely many $n\in \mb N$ and corresponding $5\e$-separated set in $M_{\o_0}$ with respect to the Bowen metric $d_{M\times \O,n}^{\phi}$, whose cardinality is greater than $e^{n(\overline h(M\times\O|\O)-\frac13\g)}$.

Firstly, we fix an such $\e$ as  above. Without loss of generality, we assume that $\e<\frac12\b_0$, where $\b_0$ is as in Lemma \ref{L:Shadowing}. Letting $\b=\frac13\e$ and applying Lemma \ref{L:Shadowing}, we have that there exists $\a>0$ corresponding to $\b$, which is taken to be smaller than $\b$.  Then we fix a $\d$ such that $\d<\frac12\e$ and
\begin{equation}\label{E:JumpControl}
d_M(\pi_M\phi(x,\o_1),\pi_M\phi(x,\o_2))<\frac16 \a,\ \forall x\in M\text{ and }d_\O(\o_1,\o_2)<\d.
\end{equation}
Let $\{U_0,U_1,\cdot,U_p\}$ be a collection of  $\frac13\a$-balls in $M$ such that $M=\cup_{i=1}^pU_i$. By Condition H2), there exists $N_1\in \mb N$ such that the following holds
$$\pi_M\phi^n(U_i,\o)\cap U_j\neq \emptyset,\ \forall i,j\in \{0,1,\cdots,p\},\ n\ge N_1.$$
Now we fix an $N_2>> N_1$ such that there exists a $5\e$-separated set in $M_{\o_0}$ with respect to the Bowen metric $d_{M\times \O,N_2}^{\phi}$, whose cardinality is greater than $e^{N_2(\overline h(M\times\O|\O)-\frac13\g)}$.  By the continuity of $\phi$, one can easily obtain the following technical lemma:
\begin{lem}\label{L:SeparatedSetFiber}
There exist a nonempty open subset $O\subset \O$ and  $k$ points $y_i:\in M,1\le i\le k$ with $$k\ge  e^{n(\overline h(M\times\O|\O)-\frac13\g)}$$  such that for any $\o\in O$, $\{(y_i,\o)\}_{1\le i\le k}$ is a $3\e$-separated set with respect to $d^{\phi}_{M\times \O,N_2}$.
\end{lem}
Actually, $O$ can be taken to be a small open neighborhood of $\o_0$. For the sake  of convenience and without loss of generality, we assume that the following hold:
\begin{itemize}
\item[a1)] $O$ is in a  $\frac12\d$-ball with respect to $d^{\t}_{\O,N_2}$, where $\d$ is as in (\ref{E:JumpControl});
\item[a2)] $\frac{\log k}{N_2}> \overline h(M\times \O|\O)-\frac23\g$;
\item[a3)] $\{y_i\}_{1\le i\le k}\subset U_0$ (Note that  $M$ is compact).
\end{itemize}


Let $N=N_2+N_1$. Note that both $\O_1:=\bigcup_{i\in\mb Z}\t^{iN}(O)$ and $\O_2:=\O\setminus \O_1$ are $\t^N$-invariant sets since $\t$ is invertible. Clearly, $\O_1$ and $\O_2$ are both measurable ( actually $\O_1$ is open and $\O_2$ is closed). We will deal with $\O_1$ and $\O_2$ separately.

Following exactly the same procedure as the proof of Lemma \ref{L:Horseshoe}, we can construct a  finite Borel measurable partition of $\O$, which is denoted by $\xi:=\{\O'_t, t\in T\}$,  Borel measurable simple functions $\{g^{0,0},g^{i,0},\ 1\le i\le p\}$, and a finite set $\{\o_t\in \O'_t,t\in T\}$ satisfying the following properties:
\begin{itemize}
\item[b1)] For any $\o,\o'$ fall in an element of $\xi$, the following holds
$$d_\O(\t^{ n}\o,\t^{ n}\o')<\d,\ \forall 0\le n\le N-1,$$
where $\d$ satisfies (\ref{E:JumpControl});
\item[b2)] For $1\le i\le p$, $g^{i,0}=\sum_{t\in T}g^{i,0}_t\chi_{\O'_t}$, where
$$g^{i,0}_t\in U_i\text{ and }\pi_M\phi^{N_1}(g^{i,0}_t,\o_t)\in U_0;$$
\item[b3)] $g^{0,0}=\sum_{t\in T}g^{0,0}_t\chi_{\O'_t}$, where
$$g^{0,0}_t\in U_0\text{ and } \pi_M\phi^{N}(g^{0,0}_t,\o_t)\in U_0 \text{ for some }\o_t\in \O'_t.$$
\end{itemize}
For each $i\in[1,k]\cap \mb N$,  there exists a finite Borel partition of $O$, which is denoted by $\eta_i:=\left\{O^i_q|q\in\{0,1,\cdots, p\}\right\}$, such that the following holds:
\begin{itemize}
\item[c1)]
$\pi_M\phi^{N_2}(\{y_{i}\}\times O^i_{q})\subset U_q$, where $\{y_i\}_{1\le i\le k}$ is as in Lemma \ref{L:SeparatedSetFiber}.
\end{itemize}

 Now, we are ready to construct the full random horseshoe by constructing proper pseudo orbits.

\noindent
{\bf \em  Step 1.} For each $\o\in \O_1=\bigcup_{i\in\mb Z}\t^{iN}(O)$, we first collect all the hitting time of orbit $\{\t^{iN}\o\}_{i\in\mb Z}$ with set $O$, which is denoted by $\{n_i(\o)\}_{i\in I(\o)}\subset \mb Z\cup\{\pm \infty\}$. We arrange $I(\o)$ and $n_i(\o)$ in the following way:

 \noindent
 {\em Case (1.1).} $n_i(\o)$ has the same sign as $i$, and for any $i,j\in I (\o)$,
 $$i< j\text{ if and only if  }n_i(\o)< n_j(\o) ;$$

  \noindent
 {\em Case (1.2).} For the positive part of $I(\o)$, if $\{\t^{iN}\o\}_{i\in\mb N}$ hits $O$ infinitely many times, let
 $$I(\o)\cap \mb N=\mb N;$$
 otherwise, let
 $$I(\o)\cap \mb N=\{1,2,\cdots, m^+(\o)\}\text{ and }n_{m^+(\o)+1}=+\infty,$$
 where $m^+(\o)$ is the number of times  $\{\t^{iN}\o\}_{i\in\mb N}$ hits $O$;

  \noindent
 {\em Case (1.3).} Similarly, for the negative part of $I(\o)$,  if $\{\t^{-iN}\o\}_{i\in\mb N}$ hits $O$ infinitely many times, let $$I(\o)\cap -\mb N=-\mb N;$$
 otherwise, let
 $$I(\o)\cap -\mb N=\{-1,-2,\cdots, -m^-(\o)\}\text{ and }n_{m^-(\o)-1}=-\infty,$$ where
 $m^-(\o)$ is the number of how many times  $\{\t^{-iN}\o\}_{i\in\mb N}$ hitting $O$;

  \noindent
 {\em Case (1.4).} If $\o\in \t^{-N_1}O$, put $0$ in $I(\o)$ and let $n_0(\o)=0$; otherwise, $0\notin I(\o)$.\\

  We remark here that $I(\o)$ is always non-empty for $\o\in \O_1$, however, each of $I(\o)\cap \mb N$, $I(\o)\cap \{0\}$, and $I(\o)\cap -\mb N$ could be an empty set.
 \medskip

 \noindent
 {\bf \em Step 2.} For a given $\o\in\O_1$ and $\hat a=(a_n)\in \mc S^k$, we define a pseudo orbit as in the following:

  \noindent
 {\em Case 1.} For $n=n_i(\o)$, $i\in I(\o)$, suppose that
 $$\t^{nN}\o\in O^{a_n}_q \text{ and } \t^{nN+N_2}\o\in \O'_{t_1},$$
  then define
 \begin{equation}\label{E:POCase1}
 y'_{\hat a,\o}(nN+j)=
 \begin{cases}
 \pi_M\phi^{j}(y_{a_n},\t^{nN}\o),&0\le j\le N_2-1\\
 \pi_M\phi^{j-N_2}(g^{q,0}(\t^{nN+N_2}\o),\o_{t_1}),&N_2\le j\le N-1\\
 \end{cases}.
 \end{equation}

  \noindent
 {\it Case 2.} For $n\notin\{n_i(\o)\}_{i\in I(\o)}$, suppose that $\t^{nN}\o\in \O'_{t_2}$,
 then define
 \begin{equation}\label{E:POCase2}
 y'_{\hat a,\o}(nN+j)=\pi_M\phi^j(g^{0,0}(\t^{nN}\o),\o_{t_2}),\ 0\le j\le N-1.
 \end{equation}

 For Case 1., note that $\{(y'_{\hat a,\o}(nN+j),\t^{nN+j}\o\}_{0\le j\le N_2-1}$ is a true orbit, and $\t^{nN}\in O^{a_n}_q$, therefore
 $$\pi_M\phi(y'_{\hat a,\o}(nN+N_2-1),\t^{nN+N_2-1}\o),\ y'_{\hat a,\o}(nN+N_2)\in U_q,$$
 thus
 $$d_M\left(\pi_M\phi(y'_{\hat a,\o}(nN+N_2-1),\t^{nN+N_2-1}\o),\ y'_{\hat a,\o}(nN+N_2)\right)\le \frac23\a.$$


 By b1) and (\ref{E:JumpControl}), we have that for $N_2\le j\le N-1$,
 \begin{align*}
 &d_M\left(\pi_M\phi\left( y'_{\hat a,\o}(nN+j),\t^{nN+j}\o\right), y'_{\hat a,\o}(nN+j+1)\right)\\
 =&d_M\left(\pi_M\phi\left(\pi_M\phi^j(g^{q,0}(\t^{nN+N_2}\o),\o_{t_1}),\t^{nN+N_2+j}\o\right), \pi_M\phi^{j+1}(g^{q,0}(\t^{nN+N_2}\o),\o_{t_1})\right)\\
 \le &\frac16\a,
 \end{align*}
 which implies that $\{(y'_{\hat a,\o}(nN+j),\t^{nN+j}\o\}_{N_2\le j\le N-1}$ is a segment of an $(\o,\frac16\a)$-pseudo orbit.
 Hence $\{(y'_{\hat a,\o}(nN+j),\t^{nN_j}\o\}_{0\le j\le N-1}$ is a segment of an $(\o,\frac23\a)$-pseudo orbit with length $N$. \\

 For Case 2., by b3), (\ref{E:JumpControl}), and using the same argument as for $N_2\le j\le N-1$ in Case 1., we have that
$\{(y'_{\hat a,\o}(nN+j),\t^{nN_j}\o\}_{0\le j\le N-1}$ is a segment of an $(\o,\frac16\a)$-pseudo orbit with length $N$.

Also note that, all the segments of $(\o,\a)$-pseudo orbit defined in both Case 1. and Case 2. start and end in $U_0$ which has diameter $\frac23\a$. Hence, the pseudo orbit $\{(y'_{\hat a,\o}(j),\t^{j}\o\}_{j\in \mb Z}$ is an $(\o,\frac23\a)$-pseudo orbit.

It is clear that $\{(y'_{\hat a,\o}(j),\t^{j}\o\}_{j\in \mb Z},\ \o\in\O$ induces an $\a$-pseudo orbit of $\tilde\phi$ by letting
$$g'_{\hat a,j}(\o)=y'_{\hat a,\o}(j),\ j\in \mb Z,\ \o\in\O_1.$$
The Borel measurability of each $g'_{\hat a,j}$ follows from (\ref{E:POCase1}) and (\ref{E:POCase2}) immediately, since each $g'_{\hat a,j}$ is either a simple function or an image of a simple function under iterations of $\phi$ or $\phi^{-1}$. Therefore, by applying Theorem \ref{L:LShad} for $\tilde\phi|_{L^\infty(\O_1)}$, we have that   there exists a unique true orbit of $\tilde\phi|_{L^\infty(\O_1)}$, $\{(\tilde\phi|_{L^\infty(\O_1)})^n(g_{\hat a})\}_{n\in\mb Z}$, which is $\b$-shadowing $\{g'_{\hat a,j}\}_{j\in \mb Z}$.

\medskip

Now we define the desired horseshoe map by the letting
$$\Psi(\hat a)|_{\O_1}=g_{\hat a},\ \forall \hat a\in \mc S_k.$$
The continuity of $\Psi|_{\O_1}$ follows from  Lemma \ref{L:ContinousShadowing} straightforwardly. The next, we prove that   Condition i) of Definition \ref{D:RandomHorseshoe} for $\Psi|_{\O_1}$ holds.

For $\hat a'\neq\hat a'' (\in \mc S_k)$, let $s=\min\left\{|n|\big |\ a'_n\neq a''_n\right\}$ and $r$  be such that $|r|=s$ and $a'_r\neq a''_r$ .
By the choice of $O,y_i,U_i,\e,\a,\b$ at the beginning of this section,  we have that
 $$d^\phi_{M\times \O,N_1}\left(\phi^{rN}(\Psi(\hat a_1)(\o),\o),\phi^{rN}(\Psi(\hat a_2)(\o),\o)\right)>\e,\ \forall \o\in \t^{-rN}(O),$$
 where $d^\phi_{M\times \O,N_1}$ is the Bowen metric.
 Note that there is a constant $L>1$ such that for any $\o\in\O$ and $x,y\in M$ with $d(x,y)>0$, we have
 $$d(\phi^{\pm1}(x,\o),\phi^{\pm1}(y,\o))\le Ld_M(x,y).$$
 Therefore, for any $\o$ from the non-empty open set $\t^{-rN}(O)$, the following holds,
 \begin{align*}
& d_M\left(\pi_M\phi^{rN}\left(\Psi(\hat a_1)(\o),\o\right),\pi_M\phi^{rN}\left(\Psi(\hat a_2)(\o),\o\right)\right)\\
=&d\left(\phi^{rN}\left(\Psi(\hat a_1)(\o),\o\right),\phi^{rN}\left(\Psi(\hat a_2)(\o),\o\right)\right)\\
\ge& L^{-N_1}d^\phi_{M\times \O,N_1}\left(\phi^{rN}(\Psi(\hat a_1)(\o),\o),\phi^{rN}(\Psi(\hat a_2)(\o),\o)\right)\\
>& \e L^{-N},
 \end{align*}
 which implies Condition i) of  Definition \ref{D:RandomHorseshoe} by taking $\Delta=  \e L^{-N}$.
  The continuity of $\Psi^{-1}$ follows from the same argument as (\ref{E:PsiInverseContinuity}).

 \noindent {\bf \em  Final Step.} To complete the proof, we  need to address the following  issues:

 Firstly, Condition ii) in Definition \ref{D:RandomHorseshoe} for $\Psi|_{\O_1}$ simply follows from the uniqueness of $g_{\hat a}$ and the same argument used in the proof of Lemma \ref{L:Horseshoe}.

 Secondly,  we need to extend the definition of $g_{\hat a}$ properly to the whole $\O$. This can be done by Theorem \ref{T:TheoryAnosovMix}. Since $\O_1$ and $\O_2=\O\setminus \O_2$ are both $\t^N$-invariant, we can view $\O_2$ independently, and give up to capture entropy from the system $\phi^N|_{\O_2}$ (as the entropy captured from $\phi^N|_{\O_1}$ is already enough). Theorem \ref{T:TheoryAnosovMix} implies that for large $N$, there exists a Borel function $g:\O_2\to M$ such that
 $$\phi^N|_{\O_2}(graph(g))=graph(g).$$
 We simply extend $g_{\hat a}$ to $\tilde g_{\hat a}$ by letting
 $$\tilde g_{\hat a}|_{\O_1}=g_{\hat a}\text{ and }\tilde g_{\hat a}|_{\O_2}=g.$$
 It is not hard to see that such an extension is harmless for  Conditions i) and ii) of Definition \ref{D:RandomHorseshoe}. So the map $\Psi:\mc S_k\to L^\infty$ which takes $\hat a$ to $\tilde g_{\hat a}$ gives a full random horseshoe.

 Finally, as we required at the beginning of this section that $N_2>>N_1$, a2) yields
 $$\frac{\log k}N=\frac{\log k}{N_2+2N_1}>\frac{\log k}{N_2}-\frac13\g>\overline h(M\times\O|\O)-\g.$$
This completes  the proof of theorem.

\end{proof}


\section{Density of Periodic Random Invariant Measures}

In this section, we prove Theorem \ref{T:PeriodicApprox}, namely, any $\phi$-invariant measure with marginal $\mb P$ can be approximated by random periodic measures.

\vskip0.05in
\noindent
{\it Proof of Theorem \ref{T:PeriodicApprox}.} We first fix a $\t$-invariant probability measure $\mb P$ on $\O$. We will show that for any $\mu\in \mc I_{\mb P}(M\times\O)$, there is a sequence of periodic measures $\{\mu_i\}_{i=1,2,\cdots}$ such that for any continuous real function $h$ defined on $M\times\O$ the following holds
$$\int h d\mu=\lim_{i\to\infty} \int hd\mu_i.$$
Since $C(M\times\O)$, the space of continuous real functions defined on a compact set, is separable, it suffices to show that for any $\mu\in \mc I_{\mb P}(M\times\O)$, $\e>0$ and $h\in C(M\times \O)$, there is a periodic measure $\mu_{\e}^h$ such that
\begin{equation}\label{E:MeasAppr}
\left|\int hd\mu-\int hd\mu_{\e}^h\right|<\e.
\end{equation}
To the end, based on the above assertion, it is not hard to see that $\{\mu_{\frac1n}^{h_m}\}_{n,m\in \mb N}$ forms an countable dense subset of $\mc I_{\mb P}(M\times\O)$, where $\{h_m\}_{m\in\mb N}$ is a dense subset of $C(M\times\O)$.

Fixing $\mu\in \mc I_{\mb P}(M\times\O)$, $\e>0$ and $h\in C(M\times \O)$, we now construct $\mu_{\e}^h$.

Since $h$ is uniformly continuous on $M\times \O$ by the compactness of $M\times\O$,
there exists $\d_1>0$ such that the following holds
\begin{equation}\label{E:Delta1AndEpsilon}
|h(x_1,\o)-h(x_2,\o)|<\frac1{10}\e,\text{ provided }d_M(x_1,x_2)<\d_1.
\end{equation}

By Lemma \ref{L:LShad}, there exists $\d_2>0$ such that any $\d_2$-pseudo orbit  can be $\frac13\d_1$-shadowed by a unique true orbit of $\tilde\phi$.
By Lemma \ref{L:JoinSegments}, there exists an $N_1\in \mb N$ such that for any $g_1,g_2\in L^\infty(\O, M)$,  there exists a $g_3\in L^\infty(\O, M)$ such that
\begin{equation}\label{E:JoinGraph}
\max\left\{d_{L^\infty(\O, M)}(g_3,g_1),d_{L^\infty(\O, M)}(\tilde\phi^{N_1}(g_3),g_2)\right\}<\frac13\d_2.
\end{equation}
For $m\in\mb N$ and $s>0$, let $$A_{m,s}^h:=\left\{(x,\o)\in M\times\O\Big|\ \left|\int hd\mu-\frac1n\sum_{i=0}^{n-1}h\left(\phi^i(x,\o)\right)\right|<s,\ \forall n\ge m\right\}.$$
It is not hard to see that $A_{m,s}^h$ is a $G_{\d}$-set by the continuity of $h$ and $\phi$, and $\lim_{m\to\infty}\mu(A_{m,s}^h)=1$ by the Birkhoff ergodic theorem. Now fix an $N_2\in\mb N$ such that the following holds
\begin{equation}\label{E:LargePeriod}
\mb P \left(\pi_\O\left(A_{N_2,\frac15\e}^h\right)\right)>1-\frac1{4M}\e\text{ and }\frac{N_1}{N_2}<\frac1{10M}\e,
\end{equation}
where $M=\sup_{M\times\O}|h|$.

Note that $\mu$ is regular. Then, there exists a compact set $K\subset A_{N_2,\frac15\e}^h$ such that
$$\mb P \left(\pi_\O\left(K\right)\right)>1-\frac14\e.$$
Applying Theorem \ref{T:BMST}, we have that there is a Borel map $g:\O\to M$ such that the following holds
$$g(\o)\in K,\ \forall x\in \pi_\O\left(K\right).$$

Let $N=N_1+N_2$. By (\ref{E:JoinGraph}), we have that there exists an $\d_2$-pseudo orbit $\{(g_n'\}_{n\in \mb Z}$ satisfying the following conditions:
\begin{itemize}
\item[a)] $g'_{lN}=g,\ \forall l\in \mb Z$;
\item[b)] $g'_{lN+k}=\tilde\phi^{k}(g),\ \forall 0\le k\le N_2-1$.
 \end{itemize}
Therefore, by Lemma \ref{L:LShad}, there exists a true orbit $\{\tilde\phi^n(\tilde g_0)\}_{n\in\mb Z}$ of $\tilde\phi$,  which is $\frac13\d_1$-shadowing $\{g_n'\}_{n\in \mb Z}$.  By the uniqueness of the shadowing orbit, it is not hard to see that $\tilde g_0$ is a random periodic point with period $N$.

Next, we show that the periodic measure induced by $graph(\tilde g_0)$
$$\mu_\e^h:=\frac1N\sum_{i=0}^{N-1}\mu_{\tilde g_i},\text{ provided }graph(\tilde g_i)=\phi^i(graph(\tilde g_0)),\ 0\le i\le N-1,$$
satisfies (\ref{E:MeasAppr}).

Note that $\mb P$ can be decomposed into the following form
$$\mb P=\sum_{i\in I} a_i\nu_i,$$
where $card(I)\le N$, $a_i>0$, $\sum_{i\in I}a_i=1$, and $\nu_i$ is $\t^N$-ergodic.

By Birkhoff ergodic theorem, there exists an $\phi$-invariant Borel function $\tilde h:M\times \O$ such that
$$\tilde h(x,\o)=\lim_{n\to\infty}\frac1n\sum_{i=0}^{n-1}h\left(\phi^i(x,\o)\right),\ \mu_\e^h-a.e.,$$
and
$$\int \tilde hd\mu_\e^h=\int hd\mu_\e^h.$$
Suppose that $(\tilde g_0(\o),\o)$ is a regular point of $\mu_\e^h$ while $\o$ is a regular point of $\nu_i$ for some $i\in I$, which can be done because the set of regular points of an invariant measure occupies full probability. Denote that
$$\mc K(n):=card\left(\left\{k|\ 0\le k\le n-1,\t^{kN}\in \pi_\O(K)\right\}\right).$$ Then we have that,
\begin{align*}
&\left|\tilde h(\tilde g_0(\o),\o)-\int hd\mu\right|\\
=&\lim_{n\to\infty}\left|\frac1{nN}\sum_{i=0}^{nN-1}h\left(\phi^i(\tilde g_0(\o),\o)\right)-\int hd\mu\right|\\
\le&\lim_{n\to\infty}\frac1{nN}\Bigg\{\sum_{0\le k\le n-1,\t^{kN}\in \pi_\O(K)}\left|\sum_{i=kN}^{(k+1)N-1}h\left(\phi^i(\tilde g_0(\o),\o)\right)-N\int hd\mu\right|\\
&\quad\quad\quad\quad\quad+\sum_{0\le k\le n-1,\t^{kN}\notin \pi_\O(K)}\left|\sum_{i=kN}^{(k+1)N-1}h\left(\phi^i(\tilde g_0(\o),\o)\right)-N\int hd\mu\right|\Bigg\}\\
\le&\lim_{n\to\infty}\frac1{nN}\Bigg\{\sum_{0\le k\le n-1,\t^{kN}\in \pi_\O(K)}\sum_{i=kN}^{(k+1)N-1}\left|h\left(\phi^i(\tilde g_0(\o),\o)\right)-h\left(\phi^{i-kN}( g(\t^{kN}\o),\t^{kN}\o)\right)\right|\\
&\quad\quad+\sum_{0\le k\le n-1,\t^{kN}\in \pi_\O(K)}\left[\left|\sum_{i=kN}^{kN+N_2-1}h\left(\phi^{i-kN}( g(\t^{kN}\o),\t^{kN}\o)\right)-N_2\int hd\mu\right|+2N_1M\right]\\
&\quad\quad+\sum_{0\le k\le n-1,\t^{kN}\notin \pi_\O(K)}2NM\Bigg\}\\
\le &\lim_{n\to\infty}\frac1{nN}\left\{\mc K(n)\left(\frac1{10}N\e+\frac15N_2\e+2N_1M\right)+2NM(n-\mc K(n))\right\}\\
\le&\frac12\nu_i(K)\e+2M(1-\nu_i(\pi_\O (K)))\\
\le& \frac12\e+2M(1-\nu_i(\pi_\O(K))).
\end{align*}
In the above argument, condition b), (\ref{E:Delta1AndEpsilon}) and  $K\subset A_{N_2,\frac15\e}^h$ are used.

Therefore, as $\tilde h$ being $\phi$-invariant, we have that
\begin{align*}
&\left|\int \tilde hd\mu_\e^h-\int h d\mu\right|=\left|\int \tilde hd\mu_{\tilde g_0}-\int h d\mu\right|\\
\le&\int\left|\tilde h(\tilde g_0(\o),\o)-\int hd\mu\right|d\mu_{\tilde g_0}=\int\left|\tilde h(\tilde g_0(\o),\o)-\int hd\mu\right|d\mb P\\
=&\sum_{i\in I}a_i\int\left|\tilde h(\tilde g_0(\o),\o)-\int hd\mu\right|d\nu_i\le\sum_{i\in I}a_i\left(\frac12\e+2M(1-\nu_i(\pi_\O(K)))\right)\\
=&\frac12\e+2M(1-\mb P(\pi_\O(K)))<\e,
\end{align*}
which completes the proof. \qed


\section{Liv\v sic Theorem}
Throughout this section , we let  $\mb P$  be a fixed $\t$-invariant probability measure on $\O$.
Before proving the main result, we introduce several notions and derive two lemmas, namely Random specification and Bowen property.

\subsection{Random specification} The specification property plays an important role in the study of equilibrium state of uniformly hyperbolic systems such as Anosov diffeomorphisms or Axiom A systems. Under the setting of this paper, specification can be introduced for random systems in a natural way.  For the sake of simplicity, we first introduce some notions.

For $g\in L^{\infty}(\O)$ and $n\in \mb N$, define that
\begin{align*}
&(g,n)=\left\{(g_0,g_1,\cdots,g_{n-1})\in \left(L^\infty(\O, M)\right)^n\big|\ g_i=\tilde\phi^i\left(g\right),\ 0\le i\le n-1\right\},\\
&B_n(g,\e)=\left\{h\in L^\infty(\O, M)\big |\ d_{n}^{\tilde\phi}\left(g,h\right)<\e\right\},
\end{align*}
where $d_{n}^{\tilde\phi}$ is the Bowen metric induced by $\tilde\phi$ and $d_{L^\infty(\O, M)}$. We call $(g,n)$ a random orbit segment of $g$ with length $n$ (or $n$-segment of $g$ for short), and $B_n(g,\e)$ the Bowen ball of $g$ with length $n$ and radius $\e$ (or $(n,\e)$-Bowen ball of $g$ for short).
\begin{defn}\label{D:Specification}
A collection of random orbit segments $\mc G\subset L^\infty(\O, M)\times \mb N$ has specification  at scale $\e$ if there exists $\tau\in\mb N$ such that for every $\{(g_i,n_i)|\ 1\le i\le k\}\subset \mc G$, there exists a $h\in L^\infty(\O, M)$ satisfying that
$$h_i\in B_{n_i}(g_i,\e),\ 1\le i\le k$$
where $h_i\in L^\infty(\O, M)$ is given by
\begin{align*}
h_i=\tilde\phi^{\sum_{j=0}^{i-1}n_j+(i-1)\tau}h.
\end{align*}
\end{defn}

\begin{lem}\label{L:MixToSpecification}
For the systems satisfying Conditions H1) and H2), $L^\infty(\O, M)\times \mb N$ has specification at any scale.
\end{lem}
\begin{proof}
This result simply follows from Lemma \ref{L:JoinSegments} and the shadowing property (Lemma \ref{L:LShad}). We omit the detailed proof here.

We remark that the connecting time "$\tau$" in the definition of specification above only depends on the scale "$\e$" and the given system. More precisely, the connecting time is determined  by the accuracy of shadowing and the speed of mixing.
\end{proof}

\subsection{Bowen property}
For a potential $\psi\in L^1_{\mb P}(\O,\mc C^{0,\a}(M))$ and a collection of random orbit segments $\mc G\subset L^\infty(\O, M)\times \mb N$, we define the variation $V(\mc G,\psi,\mb P, \e)$ and the average variation $\bar V(\mc G,\psi,\mb P, \e)$ of $\psi$ on $\mc G$ with respect to $\mb P$ at scale $\e$ by the following
\begin{align}\begin{split}\label{D:VarP}
&V(\mc G,\psi,\mb P, \e):=\sup_{(g,n)\in \mc G,\ h\in B_n(g,\e)}\left\{\sup_{\mb P}\left|S_n\psi\left(g(\o),\o\right)-S_n\psi\left(h(\o),\o\right)\right|\right\},
\end{split}
\end{align}
\begin{align}\begin{split}\label{D:AveVarP}
&\bar V(\mc G,\psi,\mb P, \e):=\sup_{(g,n)\in \mc G,\ h\in B_n(g,\e)}\left\{\int_{\O}\left|S_n\psi\left(g(\o),\o\right)-S_n\psi\left(h(\o),\o\right)\right|d\mb P\right\},
\end{split}
\end{align}
where we set $S_n(\psi(x,\o))=\sum_{i=0}^{n-1}\psi(\phi^i(x,\o))$ for $(x,\o)\in M\times \O$.

\begin{defn}\label{D:BowenPro}
Given a collection of random orbit segments $\mc G\subset L^\infty(\O, M)\times \mb N$, a potential $\psi$ has the {\bf Bowen property} (or {\bf the average Bowen property}) on $\mc G$ with respect to  $\mb P$ at scale $\e$ if the following holds
$$V(\mc G,\psi,\mb P,\e)<\infty\ (\text{ or }\bar V(\mc G,\psi,\mb P,\e)<\infty).$$

\end{defn}
Let $L^\infty_{\mb P}(\O,\mc C^{0,\a}(M))\subset L^1_{\mb P}(\O,\mc C^{0,\a}(M))$ be the collection of functions $\psi$ with the following holds
$$\sup_{\mb P}\|\psi\|_{C^{0,\a}(M_\o)}<\infty.$$
We have the following results.

\begin{lem}\label{L:BowenPro}
Let $\phi$ satisfy Condition H1) and H2), then for $\mc G=L^\infty(\O, M)\times \mb N$, there exists $\e_0>0$ such that for any $\e\in(0,\e_0)$ each potential $\psi\in L^\infty_{\mb P}(\O,\mc C^{0,\a}(M))$ (or $\psi\in L^1_{\mb P}(\O,\mc C^{0,\a}(M))$ has the {\bf  Bowen property} ( or {\bf average Bowen property}) on $\mc G$ with respect to $\mb P$ at scale $\e$.
\end{lem}
\begin{proof}
By Lemma \ref{L:InvMani}, \ref{L:LocCoor}, and P2) in the proof of Lemma \ref{L:ClosingLemma}, for a fixed positive number $\l\in(0,\l_0)$, there exist $\e>0$, $\d\in(0,\e)$, $L>1$ such that for any $x,y\in M$ with $d_M(x,y)<\d$, the following hold
\begin{align*}
& [x,y]_\o=:W^s_{\e}(x,\o)\cap W^u_{\e}(y,\o),\text{ which is continuous on }x,y\text{ and }\o,\\
&\max\left\{d_M(x,[x,y]_\o),d_M(y,[x,y]_\o)\right\}\le Ld_M(x,y),\\
&d_M(\pi_M\phi^n(x,\o),\pi_M\phi^n(y,\o))\le e^{-n\l}d_M(x,y)\text{ for }(y,\o)\in W^s_{\e}(x,\o),\ n\ge 0,\\
&d_M(\pi_M\phi^{-n}(x,\o),\pi_M\phi^{-n}(y,\o))\le e^{-n\l}d_M(x,y)\text{ for }(y,\o)\in W^u_{\e}(x,\o),\ n\ge 0.
\end{align*}
For $(g,n)\in \mc G$ and $h\in B_n(g,\e)$, we have that for any $\o\in\O$ and $0\le i\le n-1$
\begin{align}\begin{split}\label{E:EstOnOrbit}
&d_M(\pi_M\phi^i(g(\o),\o),\pi_M\phi^i(h(\o),\o)))\\
\le &d_M\left(\pi_M\phi^i(g(\o),\o),\left[\pi_M\phi^i(g(\o),\o),\pi_M\phi^i(h(\o),\o))\right]_{\t^i\o}\right)\\
&+d_M\left(\pi_M\phi^i(h(\o),\o),\left[\pi_M\phi^i(g(\o),\o),\pi_M\phi^i(h(\o),\o))\right]_{\t^i\o}\right)\\
\le &e^{-i\l}d_M\left(g(\o),\left[g(\o),h(\o)\right]_{\o}\right)\\
&+e^{-(n-1-i)\l}d_M\left(\pi_M\phi^{n-1}(h(\o),\o),\left[\pi_M\phi^{n-1}(g(\o),\o),\pi_M\phi^{n-1}(h(\o),\o))\right]_{\t^{n-1}\o}\right)\\
\le&2Le^{-i\l}d_M(g(\o),h(\o))+2Le^{-(n-1-i)\l}d_M(\pi_M\phi^{n-1}(g(\o),\o),\pi_M\phi^{n-1}(h(\o),\o))\\
\le&4Le^{-\min\{i,n-1-i\}\l}\e.
\end{split}
\end{align}

Then
\begin{align}\begin{split}\label{E:VarEstBP}
&\sup_{\mb P}|S_n(\psi(g(\o),\o))-S_n(\psi(h(\o),\o))|\\
\le &\sup_{\mb P}\|\psi\|_{C^{0,\a}(M_\o)}(4L\e)^\a\sum_{i=0}^{n-1}e^{-\min\{i,n-1-i\}\a\l}\\<&\frac{2 \sup_{\mb P}\|\psi\|_{C^{0,\a}(M_\o)}(4L\e)^\a}{1-e^{-\a\l}},
\end{split}
\end{align}
which is independent on $n$, $g$ and $h$.

For the average Bowen property, we need only to replace the $\sup_{\mb P}$ in (\ref{E:VarEstBP}) by integration on $\O$ with respect to $\mb P$ and obtain the following
\begin{align}\begin{split}\label{E:VarEstABP}
&\int_{\O}|S_n(\psi(g(\o),\o))-S_n(\psi(h(\o),\o))|d\mb P\\
\le&\int_{\O}\sum_{i=0}^{n-1}\left|\psi\left(\phi^i(g(\o),\o)\right)-\psi\left(\phi^i(h(\o),\o)\right)\right|\\
\le &\int_{\O}(4L\e)^\a\sum_{i=0}^{n-1}\|\psi\|_{C^{0,\a}(M_{\t^i\o})}e^{-\min\{i,n-1-i\}\a\l}d\mb P\\
<&\frac{2(4L\e)^\a}{1-e^{-\a\l}}\int_{\O} \|\psi\|_{C^{0,\a}(M_\o)}d\mb P,
\end{split}
\end{align}
which is independent on $n$, $g$ and $h$.
This completes the proof of this lemma.
\end{proof}

\medskip


\subsection{Proof of the weak version of random Liv\v sic Theorem} This proof contains two main steps: in the first step, we show that there is an $\a$-H\"older functional $\hat\Psi$ defined on $C(\O,M)$ so that the cohomology  equation (\ref{E:CoBoundary}) holds for all $g\in C(\O,M)$, where $C(\O,M)$ is the space of continuous maps from $\O$ to $M$ endowed with the supreme metric; in the second step, we extend the domain of functional $\hat \Psi$ from $C(\O,M)$ to $L^\infty(\O, M)$ with (\ref{E:CoBoundary}) kept. Note that the second step  is not a straightforward consequence of the first step since $C(\O,M)$ is {\bf NOT }dense in $L^\infty(\O, M)$.


\begin{lem}\label{L:WLivThC}
For any $\Phi\in L^1_{\mb P}(\O,C^{0,\a}(M))$ satisfying (\ref{E:PeriodicNull}), there exists a bounded functional $\hat\Psi:C(\O,M)\to \mb R$ which is $\a$-H\"older continuous and satisfies the following cohomology equation
\begin{equation}\label{E:CoBoundary}
\tilde\Phi(h)=\hat \Psi(\tilde\phi(h))-\hat\Psi(h),\ \forall h\in C(\O,M).
\end{equation}
\end{lem}
\begin{proof}
At first, we prove a the following technical lemma:
\begin{lem}\label{L:TransitiveOnC(Ome,M)}
There exists an orbit of $\tilde\phi$, whose closure contains $C(\O,M)$.
\end{lem}
\begin{proof}
Note that $C(\O,M)$ is separable since $M$ and $\O$ are compact. Let $\{g_i\}_{i=1,2,\cdots}\subset C(\O,M)$ be a countable dense set of $C(\O,M)$. Denote $U_{ij}$ the $\frac1j$-ball centered at $g_i$ in $C(\O,M)$.
We relabel the countable set $\left\{U_{i,j}\right\}_{i,j=1,2,\cdots}$ and denote the relabeled set by $\{V_k\}_{k=1,2\cdots}$ which satisfies the following rule:
$$\text{ for }V_{k_1}=U_{i_1j_1}, V_{k_2}=U_{i_2j_2},\ k_1\le k_2\text{ iff } i_1+j_1\le i_2+j_2\text{ or } i_1+j_1= i_2+j_2,i_1\le i_2.$$
We also define a function $l:\mb N\to \mb N$ by letting $l(k)=i$ when $V_k=U_{ij}$.\\

Let $\a$ be as in Lemma \ref{L:LShad} for a fixed $\b\in(0,\b_0)$.
By Lemma \ref{L:JoinSegments}, for any $n_k,n_{k+1}\in\mb N$, there exist $\{g'_k\in L^\infty(\O, M)\}_{k=1,2,\cdots}$ and $m\in\mb N$ such that
$$\max\left\{d_{L^\infty(\O, M)}\left(\tilde\phi^{n_k}g_{l(k)},g'_k\right),\ d_{L^\infty(\O, M)}\left(\tilde\phi^{m}g'_k,\tilde\phi^{-n_{k+1}}g_{l(k+1)}\right)\right\}<\a.$$
Therefore, we have that $\mc G_{0},\mc G_1,\mc G_2,\cdots$ form an $\a$-pseudo orbit of $\tilde\phi$, where we let
\begin{align*}
 \mc G_{0}&=\{\tilde\phi^{-n}g_{l(1)}\}_{n\ge 1}\\
\mc G_k&=\left\{g_{l(k)},\cdots,\tilde\phi^{n_k-1}g_{l(k)},g'_k,\cdots,\tilde\phi^{m-1}g'_k,\tilde \phi^{-n_{k+1}}g_{l(k+1)},\cdots, \tilde\phi^{-1}g_{l(k+1)}\right\},\ k\in \mb N.
\end{align*}
  Thus,  Lemma \ref{L:LShad} implies that there exists a $g\in L^\infty(\O, M)$ such that
$$d_{L^\infty(\O, M)}\left(\tilde\phi^{\pm(n_1-1)}g,\tilde\phi^{\pm(n_1-1)}g_{l(1)}\right)<\b,$$
 and
 $$d_{L^\infty(\O, M)}\left(\tilde\phi^{\pm (n_{k+1}-1)}\left(\tilde\phi^{\sum_{i=1}^k(n_i+n_{i+1}+m)}g\right),\tilde\phi^{\pm n_{k+1}-1}g_{l(k+1)}\right)<\b,\ \forall k\ge 1.$$
By applying Lemma \ref{L:ContinousShadowing} and choosing $n_k$ large enough for each $k\ge 1$, we have the following:
\begin{equation}\label{E:ADenseOrbit}
g\subset V_1\text{ and }\tilde\phi^{\sum_{i=1}^k(n_i+n_{i+1}+m)}g\subset V_{k+1},\ \forall k\ge 1,
\end{equation}
which implies that
\begin{equation}\label{E:ADenseOrbit1}
C(\O,M)\subset \text{ closure of }\{\tilde \phi^n(g)\}_{n\ge 0} \text{ in }L^\infty(\O, M).
\end{equation}
\end{proof}
\medskip

Now we define a functional on $\{\tilde \phi^n(g)\}_{n\ge 0}$ by the following
\begin{equation}\label{E:PsiOnOrb}
\hat \Psi'(\tilde \phi^n(g))=\sum_{i=0}^n\tilde\Phi(\tilde \phi^i(g))
\end{equation}

From Lemma \ref{L:MixToSpecification} and the second part of Lemma \ref{L:BowenPro}, it follows  that $\hat \Psi ''$ is bounded on $\{\tilde \phi^n(g)\}_{n\ge 0}$. In fact, there exist $N\in\mb N$ and $\e\in (0,\e_0)$ (where $\e_0$ is as in Lemma \ref{L:BowenPro}) such that for any $n\in \mb N$, there exists a periodic point $h\in L^\infty(\O, M)$ such tha the following hold:
$$\tilde \phi^{n+N}(h)=h\text{ and }d_{L^\infty(\O, M)}(\tilde \phi^i(g),\tilde \phi^i(h))\le \e,\ i=0,1,\cdots, n.$$
Therefore,  by (\ref{E:VarEstABP}) we have that
\begin{align}\begin{split}\label{E:PsiSupNormEst}
\left|\hat \Psi'(\tilde \phi^n(g))\right|&\le \sum_{i=n+1}^{n+N-1}\left|\tilde \Phi(\tilde\phi^i(h))\right|+\sum_{i=0}^n\left|\tilde\Phi(\tilde\phi^i(h))-\tilde\Phi(\tilde\phi^i(g))\right|\\
&\le \left(N+\frac{2(4L\e)^\a}{1-e^{-\a\l}}\right)\int_{\O} \|\Phi\|_{C^{0,\a}(M_\o)}d\mb P,
\end{split}
\end{align}
which is independent on $n$.

Let $\a_0>0$ and $C$ be the constants as in the second part of Lemma \ref{L:LShad}. We have that for any $\a\in(0,\a_0)$ each $\a$-pseudo orbit of $\tilde\phi$ can be $C\a$-shadowed. Without loss of generality, we assume that $C\a_0<\e_0$, where $\e_0$ is as in Lemma \ref{L:BowenPro}. Suppose that for $n,k\in\mb N$ $d_{L^\infty(\O, M)}(\tilde \phi^{n+k}(g),\tilde \phi^{n}(g))<\a_0$. By applying Lemma \ref{L:LShad}, there exists a random periodic point $g_0\in L^\infty(\O, M)$ such that the following holds
$$\tilde \phi^k(g_0)=g_0\text{ and }d_{L^\infty(\O, M)}(\tilde \phi^i(g_0),\tilde\phi^{n+i}(g))<(C+1)d_{L^\infty(\O, M)}(\tilde \phi^{n+k}(g),\tilde \phi^{n}(g)), i=0,1,\cdots, k.$$
Therefore, we obtain that
\begin{align}\begin{split}\label{E:PsiHolderNormEst}
&\left|\hat\Psi'(\tilde \phi^{n+k}(g))-\hat\Psi'(\tilde \phi^{n}(g))\right|\\
=&\left|\sum_{i=n+1}^{n+k}\tilde \Phi(\tilde\phi^i(g))-\sum_{i=1}^{k}\tilde \Phi(\tilde\phi^i(g_0))\right|\ \left(\text{by }(\ref{E:PeriodicNull})\sum_{i=1}^{k}\tilde \Phi(\tilde\phi^i(g_0))=0\right)\\
\le&\sum_{i=1}^{k}\left|\tilde \Phi(\tilde\phi^{n+i}(g))-\tilde\Phi(\tilde\phi^i(g_0))\right|\\
\le&C_1\left(d_{L^\infty(\O, M)}(\tilde \phi^{n+k}(g),\tilde \phi^{n}(g))\right)^\a \ \left(\text{Applying }(\ref{E:VarEstABP})\right),
\end{split}
\end{align}
where $C_1=\frac{2\left(4L(C+1)\right)^\a}{1-e^{-\a\l}}\int_{\O}\|\Phi\|_{C^{0,\a}(M_\o)}d\mb P$ and $L$ is as in (\ref{E:VarEstABP}).

 So far, we have shown that $\hat \Psi '$ is $\a$-H\"older continuous and bounded on $\{\tilde \phi^n(g)\}_{n\ge 0}$. As $\{\tilde \phi^n(g)\}_{n\ge 0}$ being dense in $C(\O,M)$, $\hat \Psi '$ defined by (\ref{E:PsiOnOrb}) has a unique $\a$-H\"older extension $\hat\Psi$ on $C(\O,M)$,  which completes the proof.
  \end{proof}
%

To extend $\hat\Psi$ to the required $\tilde\Psi$ in Theorem \ref{T:WLivTh}, we need the following lemma which demonstrates a kind of "absolute continuity" of $\hat\Psi$.

\begin{lem}\label{L:PsiAbsCon}
There exists $L'>0$ such that  for any given $\Phi\in L^1_{\mb P}(\O,C^{0,\a}(M))$ satisfying (\ref{E:PeriodicNull}) and the associated $\hat\Psi$ as in Lemma \ref{L:WLivThC}  the following holds
\begin{equation}\label{E:PsiAbsCon}
\left|\hat \Psi(h_1)-\hat \Psi(h_2)\right|\le L'\overline\Phi\left(3\mb P\left(\O\setminus \O_{h_1,h_2,\r}\right)\right)+L'\r^\a\int_{\O}\|\Phi\|_{C^{0,\a}(M_\o)}d\mb P,\ \forall h_1,h_2\in C(\O,M),
\end{equation}
where $$\O_{h_1,h_2,\r}=\left\{\o\in\O\big|\ d_M(h_1(\o), h_2(\o))< \r\right\}$$ and $$\overline\Phi(s):=\sup_{A\subset \O:\mb P(A)\le s}\int_{A} \|\Phi\|_{C^{0,\a}(M_\o)}d\mb P,\  s\in (0,1).$$
\end{lem}
\begin{proof}
Since $\hat \Psi'$ is bounded, it suffices to consider the case when $\r$  and  $\mb P\left(\O\setminus \O_{h_1,h_2,\r}\right)$ are small.

Let $g\in L^\infty(\O, M)$ be the point that the closure of whose orbit contains $C(\O,M)$, and $\a_0$ be as in the second part of Lemma \ref{L:LShad}, and  $\a_1\in (0,\frac16\a_0)$. By Lemma \ref{L:JoinSegments}, there exists an $N_1\in \mb N$ such that
\begin{equation}\label{E:N-connect}
 \forall g_1,g_2\in L^\infty(\O, M),\ \exists g_3\in L^\infty(\O, M) \text{ s.t. }g_3\in B(g_1,\a_1)\text{ and }\tilde\phi^{N_1} (g_3)\in B(g_2,\a_1),
 \end{equation}
  where $B(h,r)$ is the $r$-ball centered at $h$ in $L^\infty(\O, M)$.


Given $\r\in (0,\a_1)$ and $h_1,h_2\in C(\O,M)$, we denote $\O_1=\O\setminus \O_{h_1,h_2,\r}$. For $\a_2\in(0,\frac12\r)$, let $n,k\in\mb N$ be such that $\tilde \phi^n(g)\in B(h_1,\a_2)$ and $\tilde \phi^{n+k}(g)\in B(h_2,\a_2)$. Without loss of generality, we assume that $n,k>>N_1$. By the choice of $\a_1$ and $N_1$, there exists $h_3\in L^\infty(\O, M)$ satisfying the following
$$h_3\in B(\tilde \phi^{n+k-N_1}(g),\a_1)\text{ and }\tilde\phi^{N_1}(h_3)\in B(h_1,\a_1).$$
Define a periodic sequence $\{g_i\}_{i\in \mb Z}\subset L^\infty(\O, M)$ with period $k$ by the following:
\begin{equation}\label{E:PeriodicPartialConnect}
g_i=\begin{cases}
\tilde \phi^{i+n}(g), &\text{when }0\le i\le k-N_1-1\\
\tilde \phi^{i+N_1}(g)\chi_{\O\setminus \t^{-N_1}\O_1}+h_3\chi_{ \t^{-N_1}\O_1},&\text{when }i=k-N_1\\
\tilde\phi^{i-k+N_1}g_{k-N_1}, &\text{when }k-N_1+1\le i\le k-1\\
g_{i\mod k}&\text{otherwise}
\end{cases},
\end{equation}
where $\chi_{O}$ is the characteristic function of $O\subset \O$.
It is clear that
$$\tilde\phi(g_i)\in B(g_{i+1},\r+2\a_1+2\a_2),\ \forall i\in\mb Z.$$
Thus, using $\a_1,\a_2, \r\in(0,\frac16\a_0)$ and Lemma \ref{L:LShad}, we obtain that there exists a random periodic point $g'$ satisfying that
\begin{equation}\label{E:PeriodicShadowg'}
\tilde\phi^i(g')\in B(g_i,C(\r+2\a_1+2\a_2)),\ i\in\mb Z,
\end{equation}
where $C$ is the one as in Lemma \ref{L:LShad}.

Note that for each $\o\in \O\setminus \O'$ where $\O'= \O_1\cup \t^{-k}\O_1\cup \t^{-2k}\O_1$, we have that
$$
g_i(\t^i\o)=\begin{cases}
\tilde\phi^{n+i-k}(g)(\t^i\o)&k\le i\le 2k-1\\
\tilde\phi^{n+i}(g)(\t^i\o)&0\le i\le k-1\\
\tilde\phi^{n+k+i}(g)(\t^i\o)&-k\le i\le-1\\
\end{cases}.
$$
Thus, we have  a $(\r+\a_1,\o)$-pseudo orbit $\{(y_i,\t^i\o)\}_{i\in \mb Z}$ by letting
$$(y_i,\t^i\o)=
\begin{cases}
\phi^{i-2k+1}\left(g_{2k-1}(\t^{2k-1}\o),\t^{2k-1}\o\right)&i\ge 2k\\
\left(g_i(\t^i\o),\t^i\o\right)&-k\le i\le 2k-1\\
\phi^{i+k}\left(g_{-k}(\t^{-k}\o),\t^{-k}\o\right)&i\le -k-1\\
\end{cases}.$$
Here $\{(y_i,\t^i\o)\}_{i\in\mb Z}$ is a $(\r+2\a_2,\o)$-pseudo orbit because $\{\left(g_i(\t^i\o),\t^i\o\right)\}_{-k\le i\le 2k-1}$ is a segment of a $(\r+2\a_2,\o)$-pseudo orbit. By Lemma \ref{L:ClosingLemma} and the choice of $\a_2$ and $\r$, we have that there exists a true orbit $\{\phi^i(y',\o)\}_{i\in\mb Z}$ which is $(C(\r+2\a_2),\o)$-shadowing $\{(y_i,\t^i\o)\}_{i\in\mb Z}$. Recall that $\{\phi^i(g'(\o),\o)\}_{-k\le i\le 2k}$ is $(C(\r+2\a_1+2\a_2),\o)$-shadowing $\{\left(g_i(\t^i\o),\t^i\o\right)\}_{-k\le i\le 2k-1}$, which means that
$$\max_{-k\le i\le 2k-1}d_M\left(\pi_M\phi^i(g'(\o),\o),\pi_M\phi^i(y',\o)\right)\le C(2\a_1+4\a_2+2\r),\ \forall \o\in\O\setminus \O'.$$
Note that $k$ can be chosen arbitrarily large. Hence Lemma \ref{L:ContinousShadowing} implies that one can make $d_M\left(\pi_M\phi^{k-1}(g'(\o),\o),\pi_M\phi^{k-1}(y',\o)\right)$ and $d_M\left(g'(\o),y'\right)$ arbitrarily small by taking $k$ large enough. It is worth to point out that the largeness of $k$ is depending on $C(2\a_1+4\a_2+2\r)$ only while it is nothing to do with  the choice on specific $\o$ as long as $\o\in \O\setminus \O'$. Now we  choose $k$ large enough so that the following holds
$$\max\left\{d_M\left(\pi_M\phi^{k}(g'(\o),\o),\pi_M\phi^{k}(y',\o)\right),d_M\left(g'(\o),y'\right)\right\}<C\r,\ \forall \o\in\O\setminus \O'.$$
Thus
\begin{equation}
\max\left\{d_M\left(\pi_M\phi^{k}(g'(\o),\o),\tilde \phi^{n+k}(g)(\t^{k}\o)\right),d_M\left(g'(\o),\tilde \phi^{n}(g)(\o)\right)\right\}<4C\r,\ \forall \o\in\O\setminus \O',
\end{equation}
here we used the fact that $C(2\r+2\a_2)\le 4C\r$.

Therefore, by using the same argument as in (\ref{E:EstOnOrbit}), we obtain that
\begin{equation}\label{E:EstOnOrbSeg}
d_M\left(\tilde\phi^i(g')(\t^i\o),\tilde \phi^{n+i}(g)(\t^i\o)\right)\le 16e^{-\min\{i,k-i\}\l}C\r,\ \forall \ \o\in \O\setminus \O', 0\le i\le k.
\end{equation}
Now we are ready to estimate $\left|\hat\Psi'(h_2)-\hat\Psi'(h_1)\right|$:
\begin{align*}
&\left|\hat\Psi(h_2)-\hat\Psi(h_1)\right|\\
\le&\left|\hat\Psi(h_2)-\sum_{i=0}^{n+k}\tilde\Phi(\tilde\phi^i(g))\right|+\left|\hat\Psi(h_1)-\sum_{i=0}^{n}\tilde\Phi(\tilde\phi^i(g))\right|+\left|\sum_{i=n+1}^{n+k}\tilde\Phi(\tilde\phi^i(g))\right|\\
\le &2C'\a_2^\a+\left|\sum_{i=n+1}^{n+k}\tilde\Phi(\tilde\phi^i(g))\right|\ \left(C'\text{ is the H\"older constant of functional }\hat \Psi'\right)\\
=&2C'\a_2^\a+\left|\sum_{i=n+1}^{n+k}\tilde\Phi(\tilde\phi^i(g))-\sum_{i=1}^{k}\tilde\Phi(\tilde\phi^i(g'))\right|\ \left(g' \text{ is a periodic point of } \tilde\phi\right)\\
\le& 2C'\a_2^\a+\sum_{i=1}^k\left|\tilde\Phi(\tilde\phi^{i+n}(g))-\tilde\Phi(\tilde\phi^i(g'))\right|\\
\le&2C'\a_2^\a+\sum_{i=1}^k\int_{\t^i(\O\setminus \O')}\left|\Phi\left(\tilde\phi^{i+n}(g)(\o),\o\right)-\Phi\left(\tilde\phi^i(g')(\o),\o\right)\right|d\mb P\\
&+\sum_{i=1}^{k-N_1}\int_{\t^i( \O')}\left|\Phi\left(\tilde\phi^{i+n}(g)(\o),\o\right)-\Phi\left(\tilde\phi^i(g')(\o),\o\right)\right|d\mb P\\
 &+\sum_{i=k-N_1+1}^{k}\int_{\t^i( \O')}\left|\Phi\left(\tilde\phi^{i+n}(g)(\o),\o\right)-\Phi\left(\tilde\phi^i(g')(\o),\o\right)\right|d\mb P\\
 \le& 2C'\a_2^\a+\sum_{i=1}^k\int_{\t^i(\O\setminus \O')}\|\Phi\|_{C^{0,\a}(M_\o)}\left(16LC\r\right)^\a e^{-\min\{i,k-i\}\a\l}d\mb P\ \left(\text{by }(\ref{E:EstOnOrbSeg})\right)\\
 &+\sum_{i=1}^{k-N_1}\int_{\t^i( \O')}\|\Phi\|_{C^{0,\a}(M_\o)}\left(20LC\a_1\right)^\a e^{-\min\{i,k-N_1-i\}\a\l}d\mb P\ \left(\text{by }(\ref{E:PeriodicShadowg'})\text{ and }(\ref{E:EstOnOrbit})\right)\\
 &+\sum_{i=k-N_1+1}^{k}\int_{\t^i( \O')}2\|\Phi\|_{C^{0,\a}(M_\o)}d\mb P\\
 \le&2C'\a_2^\a+\frac{2\left(16LC\r\right)^\a}{1-e^{-\a\l}}\int_{\O}\|\Phi\|_{C^{0,\a}(M_\o)}d\mb P+\left(\frac{2\left(20LC\a_1\right)^\a}{1-e^{-\a\l}}+2N_1\right)\overline \Phi(\mb P(\O')),
\end{align*}
where the fact that $0<\a_2<\r<\a_1$ is used.

Note that, in the last line of above estimate, $\a_1$ and $N_1$ are fixed while $\a_2$ can be chosen arbitrarily small which consequentially forces $n$ and $k$ tend to $+\infty$. Thus
\begin{align*}
\left|\hat\Psi(h_2)-\hat\Psi(h_1)\right|\le& L'\overline\Phi(\mb P(\O'))+L'\r^\a\int_{\O}\|\Phi\|_{C^{0,\a}(M_\o)}d\mb P\\
\le& L'\overline\Phi(3\mb P(\O_1))+L'\r^\a\int_{\O}\|\Phi\|_{C^{0,\a}(M_\o)}d\mb P,
\end{align*}
where $L'=\max\left\{\frac{2\left(16C\right)^\a}{1-e^{-\a\l}},\frac{2\left(20C\a_1\right)^\a}{1-e^{-\a\l}}+2N_1\right\}$. The proof is complete.
\end{proof}

Now we are ready to define the functional $\tilde \Psi:L^{\infty}(\O)\to \mb R$ based on $\hat \Psi'$ and Lemma \ref{L:PsiAbsCon}. For any $h\in L^\infty(\O, M)$ and $\e>0$, by Lusin's Theorem, there exists $h_\e\in C(\O,M)$ such that $\mb P(\{\o\in\O|\ h(\o)\neq h_\e(\o)\})<\e$. Let $\{\e_i\}_{i\ge 1}$ be a sequence of positive numbers, then by taking $\r$ arbitrarily small while applying Lemma \ref{L:PsiAbsCon}, we have that
$$\left|\hat \Psi(h_{\e_n})-\hat \Psi(h_{\e_k})\right|\le L'\overline \Phi\left(3(\e_n+\e_k)\right),\ \forall n,k\ge 1.$$
Since $\|\Phi\|_{C^{0,\a}(M_\cdot)}\in L^1_{\mb P}(\O)$, the absolute continuity of $\int \|\Phi\|_{C^{0,\a}(M_\o)}d\mb P$ implies that
\[\{\hat\Psi'(h_{\e_i})\}_{i\ge 1}\]
is a Cauchy sequence as long as $\lim_{i\to\infty}\e_i=0$. Thus
$\lim_{i\to\infty}\hat \Psi(h_{\e_i})$ exists when $\lim_{i\to\infty}\e_i=0$, and moreover such limit is independent on the choice of $\{h_{\e_i}\}_{i\ge 1}$.
Thereafter, define
\begin{equation}\label{E:DefTildePsi}
\tilde \Psi(h)=\lim_{\e\to 0^+}\hat \Psi(h_{\e}).
\end{equation}
It is clear that $\tilde \Psi$ is $\a$-H\"older continuous by (\ref{E:PsiAbsCon}), and the absolute continuity of
\[\int \|\Phi\|_{C^{0,\a}(M_\o)}d\mb P\]
together with (\ref{E:CoBoundary}) implies that
$$\tilde \Phi(h)=\tilde \Psi(\tilde \phi(h))-\tilde \Psi(h),\ \forall h\in L^{\infty}(\O).$$
It remains to show that such $\tilde \Psi$ is uniquely defined up to a constant functional. Suppose that both $\tilde \Psi_1$ and $\tilde \Psi_2$  satisfy (\ref{E:CoBoundary}). Then for any $h\in L^\infty(\O, M)$, we have that $\tilde\Psi_1-\tilde\Psi_2$ is constantly valued on $\{\tilde\phi^n(h)\}_{n\in\mb Z}$. Since  the closure of $\{\tilde\phi^n(h)\}_{n\in\mb Z}$ contains $C(\O,M)$ for some $h$, and by the H\"older continuity of $\tilde \Psi$, we have that $\tilde\Psi_1-\tilde\Psi_2$ is constantly valued on $C(\O,M)$, which means that $\hat \Psi_1-\hat\Psi_2$ is a constant functional on $C(\O,M)$. Thus, $\tilde \Psi_1$ and $\tilde \Psi_2$ defined by (\ref{E:DefTildePsi}) varies at a constant functional. This completes  the proof of Theorem \ref{T:WLivTh}.

 \qed


\section{Examples}\label{Examples}
In this section, we give concrete examples which can be adapted  into the theoretic framework in previous sections. The systems we consider mainly contains two types: a) Transitive fiber Anosov system driven by a quasi-periodic driven force; b) Fiber volume preserving Anosov system with a continuous driven force, of which the precise formulations will be given later. At the end of this section, we will make some further discussions on the random horseshoe with respect to a given marginal $\mb P$, and give some relevant results.

\subsection{Systems driven by a quasi-periodic  force}\label{S:Quasi}
In this section, we consider a type of  Anosov systems which is driven by a quasi-periodic force. This section contains two main parts: in the first part, we state and prove the general result in a general setting;  in the second part, we will investigate a concrete example which is generated by rotation on torus and affine Anosov maps on 2-d torus, and will show that such a system satisfies the Conditions 1A.-1C. below and illustrate some interesting phenomenas (see Remark \ref{R:NonContiPerOrb}).

\subsubsection{Setting and results}\label{S:S1)SetResult}
We assume that the systems  satisfy the following conditions:
\begin{itemize}
\item[1A.] $(\t,\O)$ is a minimal irrational rotation on a compact torus;
\item[1B.] $\phi$ is Anosov on fibers (Condition H1));
\item[1C.] $\phi$ is topological transitive on $M\times \O$.
 \end{itemize}
 For the sake of convenience, we call the systems satisfying Conditions 1A.-1C. the S1) systems.
 The main result we will prove is the following Theorem:
 \begin{thm}\label{T:Spectrum11}
 All S1) systems satisfy Condition H2) that is topological mixing on fibers.
 \end{thm}

  We first prove a weaker version of Theorem \ref{T:TheoryAnosovMix} for S1) systems, which is the following proposition:
 \begin{prop}\label{P:DensePeriodicOrb1}
 Let $\phi$ be an S1) system. Then for any $g\in C(\O,M)$ and $\e>0$, there exists a random periodic point $\tilde g\in L^{\infty}(\O)$ such that
 \begin{equation}\label{E:LInfiClose}
 d_{L^{\infty}(\O)}(g,\tilde g)\le \e.
 \end{equation}
 \end{prop}
\begin{proof}
Actually,  this proposition still holds when $\t$ is assumed to be equicontinuous only. We say $\t$ is {\bf equicontinuous}  if  there exist $c>1$ and $\d>0$ such that for any $\o_1,\o_2\in \O$ with $d_\O(\o_1,\o_2)<\d$,
 $$\frac1cd_\O(\o_1,\o_2)\le d_\O(\t^n\o_1,\t^n\o_2)\le cd_\O(\o_1,\o_2), \forall n\in \mb Z.$$

 When $\t$ is equicontinuous, if we set
 $$\widetilde{d}_\O(\o_1,\o_2)=\sup_{n\in \mathbb{Z}}d_{\O}(\t^n \o_1,\t^n\o_2),\, \forall \o_1,\o_2\in \O,$$
 then $\widetilde{d}_\O$ is  a metric on $\O$  satisfying that
 \begin{itemize}
 \item $\widetilde{d}_\O$ and $d_\O$ are equivalent;

 \item $\widetilde{d}_\O(\t^n\o_1,\t^n\o_2)=\widetilde{d}_\O(\o_1,\o_2)$ for any $\o_1,\o_2\in \O$
 and $n\in \mathbb{Z}$.
 \end{itemize}
 Therefore, without loss of generality, for an equicontinuous $\t$, we always assume that
 \begin{equation}\label{E:EquDist}
 d_\O(\t^n\o_1,\t^n\o_2)=d_\O(\o_1,\o_2),\ \forall \o_1,\o_2\in\O\text{ and }n\in \mb Z.
 \end{equation}

  For a given $g\in C(\O,M)$, we define a neighborhood of its graph as follows.  We consider a finite Borel partition $\xi=\{\xi_i\}_{1\le i\le p}$ of $\O$, of which each element contains an nonempty open subset. Such a partition exists because of compactness of $\O$. For any $\o\in\O$, we define $\xi(\o)$ the element in $\xi$ which contains $\o$, and also denote that
$$B_g(\xi_i,\d):=\{(x,\o)|\ \o\in \xi_i, d_M(x,g(\o))<\d\},$$ which induces
a $\d$-neighborhood of $graph(g)$: $\bigcup_{\xi_i\in\xi} B_g(\xi_i,\d)$. We denote it by $B_g(\d)$.


\begin{lem}\label{L:MReturn}
For any $\d>0$, 
there exists $m\in\mb N$ such that for any $\xi_i\in\xi$, there exists $(x_i,\o_i)\in B_g(\xi_i,\d)$ such that
$$\phi^m(x_i,\o_i)\in B_g(\xi_i,\d).$$
\end{lem}
\begin{proof}
Since $\phi$ is transitive, there exists $(x_0,\o_0)$ such that the orbit of $\{\phi^n(x_0,\o_0)\}_{n\in\mb Z}$ is dense in $M\times \O$.
 By the finiteness of the partition, we have that there exist a point $(x,\o)$ in the interior of $B_g(\xi_1,\d)$ and  an $m_0\in\mb N$ such that
$$\phi^{n_i}(x,\o)\in \text{interior of }B_g(\xi_i,\d)\text{ for some }n_i\in[1,m_0],\ \forall i\in[2,p].$$
By the continuity of $\phi$ and the finiteness of the partition, we have that there exists a small open neighbourhood of $(x,\o)$ in $M\times \O$, which is denoted by $U$, such that
$$U\subset B_g(\xi_1,\d)\text{ and }\phi^{n_i}(U)\subset \text{interior of }B_g(\xi_i,\d),\ \forall i\in[2,p].$$
 Then, by the transitivity of the system, there exist an $m\in \mb N$ and $(x',\o')\in U$ such that $\phi^m(x',\o')\in U$. Note that $$\phi^m(\phi^{n_i}(x',\o'))=\phi^{n_i}(\phi^m(x',\o'))\in \phi^{n_i}(U) \subset \text{interior of }B_g(\xi_i,\d),\ \forall i\in[2,p].$$
We complete the proof by letting $(x_1,\o_1)=(x',\o')$ and $(x_i,\o_i)=\phi^{n_i}(x',\o')$ for $i\in[2,p]$.

\end{proof}

\bigskip

Next, we will construct a $(\u,\d)$-{\bf pseudo orbit}  for $\u\in \O$, $g\in C(\O,M)$, and $\d>0$. Given $\d>0$, by Lemma \ref{L:MReturn}, there exist $m=m(\d,g,\phi) >0$ and $(x_i,\o_i)\in B_g(\xi_i,\d)$ for $1\le i\le p$ such that
\begin{equation}\label{E:MReturn}
\phi^m(x_i,  \o_i)\in  B_g(\xi_i,\d)\text{ for }1\le i\le p.
\end{equation}
For $l\in\mb Z$, suppose that $\t^{lm}\u\in \xi_{i_l}$, then we define
\begin{equation}\label{E:DefPseOrb}
(y_j,\u_j)=\begin{cases}
(g(\t^{lm}\u),\t^{lm}\u), &\text{ when }j=lm\\
(\pi_M\phi^{j-lm}(x_{i_l},\o_{i_l}),\t^j\u),&\text{ when }j\in[lm+1,(l+1)m-1]\\
\end{cases}.
\end{equation}
\begin{lem}\label{L:PseudoOrbit}
For any $\e'>0$, there exists $\d'(\phi,g,\e')>0$ such that if $\d$ and the  diameter of the partition $\xi$  are less than $\d'$,  the pseudo orbit defined by (\ref{E:DefPseOrb}) starting from $(y,\upsilon)\in graph(g)$  is a $(\u,\e')$-pseudo orbit.
\end{lem}
\begin{proof}
Noting that  $g$ is uniformly continuous since $\O$ being compact. Therefore, for any $\e_1>0$, there exists $\d_1>0$ such that for any $\o_1,\o_2\in \O$ if $d_\O(\o_1,\o_2)<\d_1$ then
\begin{equation}\label{E:POEst1}
  d_M(g(\o_1),g(\o_2))<\e_1.
\end{equation}
 Also note that both $\phi$ and $\phi^{-1}$ are uniformly continuous, then for any $\e_2>0$ there exists a $\d_2>0$ such that for any $(x_1,\o_1),(x_2,\o_2)\in M\times \O$ once $d((x_1,\o_1),(x_2,\o_2))<\d_2$ the following holds
 \begin{equation}\label{E:POEst2}
   \max\{d(\phi(x_1,\o_1),\phi(x_2,\o_2)),d(\phi^{-1}(x_1,\o_1),\phi^{-1}(x_2,\o_2))\}<\e_2.
 \end{equation}
Note that $d_\O(\t^j\u,\t^{j-lm}\o_{i_0})$ is less than the diameter of the partition $\xi$ (which is less than $\d'(\phi,g,\e')$) for all $j\in[lm,(l+1)m]$ and $l\in \mb Z$ by (\ref{E:EquDist}). For  a given $(y,\u)\in graph(g)$ with $\u\in \xi_{i_0}$,  the following holds:
\begin{enumerate}
 \item[a)] When $j\in [lm+2,(l+1)m-1]$,  taking $\d'<\d_2$,  then (\ref{E:POEst2}) implies that
  \begin{align*}
  d_M(y_j,\pi_M\phi(y_{j-1},\u_{j-1}))
  =&d_M(y_j,\pi_M\phi(y_{j-1},\t^{j-1}\u))\\
  \le &  d(\phi(y_{j-1},\t^{j-lm-1}\o_{i_l}),\phi(y_{j-1},\t^{j-lm-1}\u_{lm}))\\
  \le& \e_2,
  \end{align*}
  where the fact that $\u_{lm}\in\xi_{i_l}$ is used.
  \item[b)] When $j=lm+1$, taking $\d'<\frac12 \d_2$ (thus $\d<\frac12\d_2$), then
  $$d\left((x_{i_l},\o_{i_l}),(g(\t^{lm}\u),\t^{lm}\u)\right)<\d_2.$$ Hence, (\ref{E:POEst2}) implies that
   $$d_M\left(y_j,\pi_M\phi(y_{j-1},\u_{j-1})\right)=d_M\left(\pi_M\phi(x_{i_l},\o_{i_l}),\pi_M\phi(g(\t^{lm}\u),\t^{lm}\u)\right)<\e_2;$$

  \item[c)] For $j=(l+1)m$,  taking $\d'<\min\{\d_1,\d_2\}$, then (\ref{E:POEst1}), (\ref{E:MReturn}), and (\ref{E:POEst2}) implies that
  \begin{align*}
  & d_M(y_j,\pi_M\phi(y_{j-1},\u_{j-1}))\\
  =&d_M\left(g(\t^j\u),\pi_M\phi(y_{j-1},\u_{j-1})\right)\\
 \le & d_M\left(g(\t^m\u_{lm}),g(\t^m\o_{i_l})\right)+d_M\left(g(\t^m\o_{i_l}),\pi_M\phi^m(x_{i_l},\o_{i_l})\right)\\
    &+d_M(\pi_M\phi^m(x_{i_l},\o_{i_l}),\pi_M\phi(y_{j-1},\t^{m-1}\u_{lm}))\\
    \le& \e_1+\d+d_M(\pi_M\phi(y_{j-1},\t^{m-1}\o_{i_l}),\pi_M\phi(y_{j-1},\t^{m-1}\u_{lm}))\\
    \le&\e_1+\d+\e_2.
  \end{align*}
\end{enumerate}
Therefore, the required estimate can be achieved  if we choose $,\e_1+\d +\e_2<\e'$ and take a $\d'>0$  which satisfies all the proposed conditions  in a), b), and c) above. So for such $\d'>0$ if $\d$ and $diam(\xi)$ are less than $\d'$, the  pseudo-orbit defined by (\ref{E:DefPseOrb}) is a $(\u,\e')$-pseudo orbit.

\end{proof}

\begin{lem}\label{L:LipG}
For any $\e'>0$ and $L>0$, there exists $\d'(\phi,L,\e')>0$ such that for any $g\in C(\O,M)$ with $Lip\ g\le L$ if $\d$ and the  diameter of the partition $\xi$  are less than $\d'$,  then the pseudo orbit defined by (\ref{E:DefPseOrb}) starting from $(y,\upsilon)\in graph(g)$  is a $(\u,\e')$-pseudo orbit.
\end{lem}
\begin{proof}
In the proof of Lemma \ref{L:PseudoOrbit} above, the dependence of $\d'$ on $g$ is given by  (\ref{E:POEst1}) which appears in part c). So the exact property of ``$g$'' being involved in this lemma is how uniformly continuous it is. Therefore, under the setting of the current lemma, the exactly same argument can be applied here, in which ``$\d'$'' depends on the Lipchitz constant ``$L$'' rather than the choice of particular map ``$g$''. The detailed proof is omitted.
\end{proof}
Now we are ready to prove Proposition \ref{P:DensePeriodicOrb1}. For a given small $\e>0$ and $g\in C(\O,M)$, choose $0<\b<\min\{\frac12\e,\beta_0\}$. By Lemma \ref{L:LShad}, there exists a corresponding $\a>0$ for $\b$. Then, by applying Lemma \ref{L:MReturn} and \ref{L:PseudoOrbit}, for $\e'=\a$ there exists a $\d'>0$ such that for any $\u\in\O$ and $\d\in(0,\d')$, there exists $m\in \mb N$ which is depending on $\d,g,\phi$ only and a $(\u,\a)$-pseudo orbit defined by (\ref{E:DefPseOrb}), which is denoted by $\{(y_i,\t^n\u)\}_{n\in\mb Z}$, such that
$$(y_{km},\t^{km}\u)\in graph(g),\ \forall k\in\mb Z.$$
Define  $g_i:\O\to M$ for $i\in \mb Z$ as follows
$$g_j(\u)=\begin{cases}
g(\u), &\text{ when }j=lm\\
\pi_M\phi^{j-lm}(x_{i_l},\o_{i_l}),&\text{ when }j\in[lm+1,(l+1)m-1]\text{ and }\t^{-(j-lm)}\u\in \xi_{i_l}\\
\end{cases}.$$
It is clear that each $g_i$ is measurable since $g_i$ is either a continuous function or a pre-image of a simple function under continuous iterations. Also note that
$$(y_i,\t^n\u)\in graph(g_i),\ \forall\u\in\O\text{ and }i\in \mb Z.$$
Thus $\{g_i\}_{i\in\mb Z}$ is an $\a$-pseudo orbit of $\tilde\phi$. By applying Lemma \ref{L:LShad}, there is a true orbit $\{\tilde\phi(\tilde g)\}_{i\in\mb Z}$ which is $\b$-shadowing $\{g_i\}_{i\in\mb Z}$. The uniqueness of such a $\tilde g$ and the periodicity of $\{g_i\}_{i\in\mb Z}$ imply that $\tilde g$ is a periodic point of $\tilde \phi$.

This completes the proof of Proposition \ref{P:DensePeriodicOrb1}.
\end{proof}

\medskip

\noindent{\bf Proof of  Theorem \ref{T:Spectrum11}.}  Actually, what was required for the driven system is the following: $\O$ is connected, and for any $n\in\mb N$, $\t^n$ is uniquely ergodic and $\mb P$ is the ergodic probability measure.

By Proposition \ref{P:DensePeriodicOrb1}, we have that for any  $g\in C(\O,M)$ and $\e>0$, there exists a  random periodic point $\tilde g\in L^\infty(\O, M)$ satisfying (\ref{E:LInfiClose}). For this $\tilde g$,  let
\begin{equation}\label{E:XTg}
X_{\tilde g}= \bigcup_{\o\in\O}\{(x,\o)| x\in \overline{W^u(\tilde g(\o),\o)}\},
\end{equation}
where $\overline{W^u(\tilde g(\o),\o)}$ is the closure of the global unstable manifold of $(\tilde g(\o),\o)$ in $(M,\o)$.   We first require the following condition on the parameters:

{\bf C1:} By taking $\e$ small enough, we require that $\e<\frac12 \d$ where $\d$ is as in Lemma \ref{L:LocCoor} corresponding to $\frac12\e_0$ with $\e_0$ being a prefixed small number.

%

First, we show that
\begin{equation}\label{E:OpenClose}
X_{\tilde g}=M\times \O.
\end{equation}

The next lemma is the key result needed in the proof:
 \begin{lem}\label{L:ContiG}
  For any $(x,\o)\in M\times\O$, $g\in C(\O,M)$, and an open neighborhood $V\subset \O$ of $\o$,  there exist an nonempty open set $A\subset\O\setminus\{\o\}$ with $\O\setminus A\subset V$ and a function $g'\in C(\O,M)$ such that
 \begin{equation}\label{E:g'}
 g'(\o)=x\text{ and }g'|_A=g|_A.
 \end{equation}
 \end{lem}
 The proof of this lemma may be standard, however we give it here for the sake of completeness.
 \begin{proof}
 First, there exists a path $\g:[0,1]\to M$ such that $\g(0)=g(\o)$ and $\g(1)=x$. Note that there exist finitely many numbers $0=s_0<s_1<s_2<\cdots<s_n=1$ such that $\g(s_i)$ and $\g(s_{i+1})$ are in a same chart  of the atlas of $M$, say $(U_i,\psi_i)$, for all $i=0,\cdots, n-1$.
 Consider an open ball $B(\o,R_1)$ centered at $\o$ in $\O$,  with $R_1$ small enough so that $B(\o,R_1)\subset V$,  $\mb P(\overline{B(\o,R_1)})<1$ and $g(B(\o,R_1))\subset U_1$. Note that  $\{\o\}$ and $\O\setminus B(\o,\frac12 R_1)$ are disjoint closed sets. By Urysohn's lemma, there exists a continuous function $f_1:\O\to[0,1]$ such that
 $$f_1|_{\O\setminus B(\o,\frac12R_1)}=0\text{ and }f_1(\o)=1.$$
 Then, the following is well defined
 $$g_1(\u)=\begin{cases}g_0(\u),&\text{ if }\u\notin B(\o,\frac12R_1)\\
\psi_1^{-1}\left[ f_1(\u)\psi_1(\g(s_1))+(1-f_1(\u))\psi_1(g_0(\u))\right],&\text{ if }\u\in B(\o,\frac 12R_1)\end{cases},$$
where we set $g_0=g$ for the sake of convenience.

 It is obvious that $g_1$ is continuous and satisfies that
 $$g_1(\o)=\g(s_1)\text{ and }A_1:=\{\o|g_1(\o)=g(\o)\}\text{ contains an open set}.$$
 Inductively, we can construct $g_{i+1}$ based on $g_i$ by using   exactly the same argument above to ensure the following:
 $$g_{i+1}(\o)=\g(s_{i+1})\text{ and }A_{i+1}:=\{\o|g_1(\o)=g(\o)\}\text{ contains an open set}.$$
 To guarantee the existence of  $A$, it is sufficient to take $R_{i+1}<R_i$. In this case, we can take $A=\O\setminus \overline{B(\o,\frac12R_1)}$.

 \end{proof}

 \begin{rem}
 Actually, it is not hard to see that $\mb P(A)$ can be taken arbitrarily close to $1$ since $\mb P$ is not  atomic.
 \end{rem}

Now, we are ready to prove Theorem \ref{T:Spectrum11}.
For an arbitrary $(x,\o)\in  M\times\O$, and $g$ chosen at the beginning of the proof, by Lemma \ref{L:ContiG}, there exists $g'\in C(\O,M)$ satisfying (\ref{E:g'}).  For a small $\e'\in (0,\frac12\d)$, where $\d$ is as in C1). We apply Proposition \ref{P:DensePeriodicOrb1} to $g'$, and then derive that there exists a random periodic point $\tilde g'\in L^\infty(\O, M)$ such that
\begin{equation}\label{E:ClosePerOrb}
\sup_{\u\in A}\{d_M(\tilde g'(\u),\tilde g(\u))\}\le 2\e',
\end{equation}
where $\tilde g$ is the periodic point in (\ref{E:XTg}).
Let $m$ and $m'$ be the periods of $\{graph(\tilde g)\}$ and $\{graph(\tilde g')\}$ respectively. Note that $\t^{-mm'}$   is $\mb P$-ergodic since $\t$ is invertible and $\t^{mm'}$ is $\mb P$-ergodic, and each $\o'\in\O$ is a $\t^{mm'}$-minimal point. That means that there is an infinite increasing sequence of positive integers, $0<n_1<n_2<n_3<\cdots$, such that $\t^{-n_imm'}\o\in A$ for $i\ge 1$ as long as $A$ has nonempty interior.

For $i\ge 1$, let
\begin{align*}
x'_i&=[\pi_M\phi^{-n_imm'}(\tilde g'(\o),\o),\tilde g(\t^{-n_imm'}\o)]_{\t^{-n_imm'}\o}\\
&=W^s_{\frac{\e_0}2}(\phi^{-n_imm'}(\tilde g'(\o),\o))\bigcap W^u_{\frac{\e_0}2}(\tilde g(\t^{-n_imm'}\o),\t^{-n_imm'}\o),
\end{align*}
which is well defined by the choice of $\e_0,\e'$ and (\ref{E:ClosePerOrb}).

Note that, by Lemma \ref{L:InvMani}, we have that
$$\lim_{i\to\infty}\phi^{n_imm'}(x'_i,\t^{-n_imm'}\o)\to (\tilde g'(\o),\o)\text{ and }\phi^{n_imm'}(x'_i,\t^{-n_imm'}\o)\in W^u(\tilde g(\o),\o),$$
which implies that
$$(\tilde g'(\o),\o)\in \overline{W^u(\tilde g(\o),\o)}\subset X_{\tilde g}.$$
Since $d_M(g'(\o),x)<\e'$ and $\e'$ can be taken arbitrarily small, we have that $(x,\o)\in X_{\tilde g}$. Thus (\ref{E:OpenClose}) holds by the arbitrariness of $(x,\o)$. \\


Finally, we show that $\phi$ is {\bf topological mixing on fibers}.
Roughly speaking, the proof is based on the fact that the global stable manifold of a  random periodic point is dense in $M\times\O$.
The following lemma is the key  to derive the mixing property.
\begin{lem}\label{L:Converging}
For any $\e_1>0$ and $g\in C(\O,M)$, there exists $\d_1>0$ such that the following holds:\\
For any  random periodic point $graph(\tilde g)\in B_g(\d_1)$ (with period $l$), $\e'_1>0$ and $h\in C(\O,M)$,
 there exist $N\in\mb N$, and a measurable function $h_1:\O\to M$ satisfying that
$$graph(h_1)\subset B_{h}(\e'_1)\text{ and }\phi^{nl}(graph(h_1))\subset B_g(\e_1),\forall n>N.$$
\end{lem}
\begin{proof}
Here we only need to take $\d_1\in(0,\frac16\e_1)$ small enough such that any $2\d_1$-pseudo orbit of $\tilde\phi$ can be $\frac13\e_1$-shadowed by a unique true orbit. This can be done by Lemma \ref{L:LShad}. For the rest of the proof of this lemma, such a $\d_1$ is fixed.

 It is not hard to see that the existence of $\tilde g$ follows from Proposition \ref{P:DensePeriodicOrb1} directly. For an arbitrarily fixed $\o\in\O$, by Lemma \ref{L:ContiG}, there exist an nonempty  open set $A_\o\subset \O$ and $g'_\o\in C(\O,M)$  such that
 \begin{equation}\label{E:g'gAomega}
 g'_\o(\o)=h(\o)\text{ and }g'_\o|_{A_\o}=g|_{A_\o}.
 \end{equation}
  By applying Proposition \ref{P:DensePeriodicOrb1}, we have that there is a random periodic point $\tilde g'_\o\in L^\infty(\O, M)$ with period $l_\o$ such that
  \begin{equation}\label{E:TilG'omega}
  graph(\tilde g'_\o)\in B_{g'_\o}\left(\frac13\e'_1\right),
  \end{equation}
  where, without loss of generality, we assume that $\e'_1\in(0,\d_1)$.

  Let $B_\o\subset A_\o$ be an open ball  in $\O$ with radius $r_\o>0$, and $C_\o$ be an open neighborhood of $\o$  with radius $\frac12r_\o$. This $r_\o$ is taken additionally small to ensure the following inequality
  \begin{equation}\label{E:GomegaHmu}
  |g'_\o(\u)-h(\u)|<\frac13\e'_1,\ \forall \u\in C_\o.
  \end{equation}
  Since $\t^{l_\o l}$ is an isometric map on $\O$ and any nonempty open subset of $\O$ is not $\mb P$-null set, there is an infinite set $I_\o\subset\mb N$ such that the following holds
  $$\t^{nl_\o l}(C_\o)\subset B_\o,\ \forall n\in I_\o.$$
  For  $\u\in C_\o$ and $n\in I_\o$, define that
  $$(y_i,\u_i)_{n}=\begin{cases}
  \phi^i(\tilde g'_\o(\u),\u)&\text{ when }i< nl_\o l\\
  \phi^{i-nl_\o l}(\tilde g(\t^{nl_\o l}\u),\t^{nl_\o l}\u)&\text{ when }i\ge nl_\o l
  \end{cases}.$$
  It is clearly $\{(y_i,\u_i)_{n}\}_{i\in\mb Z}$ is an $(\u,\frac32\d_1)$-pseudo orbit because of (\ref{E:TilG'omega}) and (\ref{E:g'gAomega}). By the choice of $\e_1$ and $\d_1$ at the beginning of the proof, we have that there exists a unique true orbit $\{\phi^i(z_{\o,n}(\u),\u)\}_{i\in\mb Z}$ which is $(\u,\frac13\e_1)$-shadowing $\{(y_i,\u_i)_n\}_{i\in\mb Z}$. Note that $\d_1+\frac32\d_1+\frac13\e_1<\frac34\e_1$, then we have that
  \begin{equation}\label{E:1}
  \phi^{ml}(z_{\o,n}(\u),\u)\in B_g\left(\frac34\e_1\right),\ \forall m\ge nl_\o.
  \end{equation}
  By Lemma \ref{L:ContinousShadowing}, there exists $N_\o\in \mb N$ such that if we choose $nl_\o> N_\o$, then
  $$d_M\left(z_{\o,n}(\u),\tilde g'_\o(\u)\right)<\frac13\e'_1,$$
  which together with (\ref{E:TilG'omega}) and (\ref{E:GomegaHmu}) implies that
    \begin{equation}\label{E:2}
  d_M(z_{\o,n}(\u),h(\u))<\e'_1.
  \end{equation}
  So far, we have defined a measurable function $z_{\o,n}:C_\o\to M$ with an integer $N_\o\in \mb N$ ($n$ is required to be larger than $N_\o$)  satisfying (\ref{E:1}) and (\ref{E:2}) for all $\u\in C_\o$. We remark here that the measurability of $z_\o$ can be proved exactly in the  same way as the argument in the proof of Proposition \ref{P:DensePeriodicOrb1}, thus is omitted here.

  Note that $\{C_\o\}_{\o\in\O}$ forms an open cover of the compact set $\O$. Thus there is an finite open sub-cover which is denoted by $\{C_{\o_1},\cdots C_{\o_s}\}$. The proof is completed if we take $N=\max\{N_{\o_1},\cdots, N_{\o_s}\}$ and define
  $$h_1(\o)=\begin{cases}
  z_{\o_1}(\o)&\text{ when }\o\in C_{\o_1}\\
  z_{\o_i}(\o)&\text{ when }\o\in  C_{\o_i}\setminus \bigcup_{1\le j\le i-1}C_{\o_j},\ \forall 2\le i\le s
  \end{cases}.
  $$
\end{proof}

Let $U$ and $V$ be nonempty open sets in $M$ with an $x\in U$ and a $y\in V$ fixed. Then there exists $\e_2>0$ such that
$$B_M(x,\e_2)\subset U\text{ and }B_M(y,\e_2)\subset V.$$
Applying  Lemma \ref{L:Converging} by taking $\e_1=\e_2$ and $g=h_y$ where $h_y(\o)=y$ for all $\o\in\O$, we can derive a $\d_2>0$ correspondingly. By Proposition \ref{P:DensePeriodicOrb1}, there exists a  random periodic point $graph(\tilde h_y)\in B_{h_y}(\d_2)$ with period $p$.
Since $\phi$ is homeomorphism, for any $i\in\mb N$, $\phi^{i}(U\times \O)$ is open in $M\times\O$ while $\phi^{i}(graph(h_x))\subset \phi^{i}(U\times\O)$ is compact.
Therefore, there exists $\e_3\in(0,\e_2)$ such that
$$B_{g_i}(\e_3)\subset \phi^{i}(U\times \O), \ \forall 0\le i \le p-1,$$
where $g_i(:=\tilde\phi(h_x))\in C(\O, M)$. Applying Lemma \ref{L:Converging} $p$ times, we have that there exist $\{N_0,\cdots, N_{p-1}\} \subset \mb N$ and measurable functions $\{g'_i:\O\to M\}_{0\le i\le p-1}$ such that for all $0\le i\le p-1$ the following hold
\begin{align}\begin{split}\label{E:3}
&graph(g'_i)\subset B_{g_i}(\e_3)\subset \phi^{i}(U\times\O)\\
&\quad\quad\quad\quad\quad\quad\quad\text{and}\\
&\phi^{np}(graph(g'_i))\subset B_{h_y}(\e_2)\subset V\times \O,\forall n>N_i.
\end{split}
\end{align}
Obviously, (\ref{E:3}) holds if we replace $N_i$ by $N:=\max\{N_0,\cdots, N_{p-1}\}$. Therefore, for any $m\in\mb N$, if $m=np+j$ where $n>N$ and $j\in[0,p)$, then
$$\phi^{np}(graph(g'_i))\subset B_g(\e_2)\bigcap \phi^{m-j}(\phi^{j}(U\times \O))\subset V\times \O \bigcap \phi^{m}(U\times \O),$$
which implies that for any $\o\in\O$,
$$\phi^{m}(U\times \{\o\})\cap V\times \{\t^{m}\o\}\neq \emptyset.$$
  This completes the proof of Theorem \ref{T:Spectrum11}.

  \qed
 \subsubsection{Fiber Anosov maps on 2-d tori}\label{S:Example1}
In this section, we give an example of Fiber Anosov maps on 2-d tori.

Let $\phi:\mathbb{T}^2\times \mathbb{T}\rightarrow \mathbb{T}^2\times \mathbb{T}$ given by
$$\phi\left(\binom{x}{y},\omega\right)=\left( A\binom{x}{y}+h(\omega),\omega+\alpha\right),$$ where $\alpha\in \mathbb{R}\setminus \mathbb{Q}$,
$A=\left( \begin{array}{cc} 1 & 1\\ 2& 1\end{array}\right)$ and $h(\omega)=\binom{h_1(\omega)}{h_2(\omega)}$ is a continuous map from $\mathbb{T}$ to $\mathbb{T}^2$.
It is clear that $\phi$ is a quasi-periodic forced system and  $\phi$ induces a cocycle over $(\mathbb{T},\theta)$, where $\theta(\omega)=\omega+\alpha$.

\begin{lem} $\phi$ is Anosov on fibers and transitive.
\end{lem}
\begin{proof} Let $f_\omega \binom{x}{y}= A\binom{x}{y}+h(\omega)$ for $\omega\in \mathbb{T}$. Then $\phi\left(\binom{x}{y},\omega\right)=\left(f_\omega \binom{x}{y},\omega+\alpha\right)$.
Note that $Df_\omega \binom{x}{y}=A$ for any $(\binom{x}{y},\omega)\in \mathbb{T}^2\times \mathbb{T}$ and $A$ is a hyperbolic matrix. Hence $(\mathbb{T}^2\times \mathbb{T},\phi)$ is Anosov on fibers.

Let $\mu$ be the Lebesgue measure on $\mathbb{T}^2\times \mathbb{T}$. To show that $T$ is transitive, it is sufficient to show that $(\mathbb{T}^2\times \mathbb{T},\mathcal{B}_{\mathbb{T}^2\times \mathbb{T}}, \phi,\mu)$ is an ergodic measure-preserving system. First, it is clear that $(\mathbb{T}^2\times \mathbb{T},\mathcal{B}_{\mathbb{T}^2\times \mathbb{T}}, \phi,\mu)$ is a measure-preserving system. In the following we show that $(\mathbb{T}^2\times \mathbb{T},\mathcal{B}_{\mathbb{T}^2\times \mathbb{T}}, \phi,\mu)$ is ergodic.

Let $f\in L^2(\mu)$ with $f\circ \phi=f$ $\mu$-a.e.. Let
$$f\left(\binom{x}{y},\omega\right)=\sum_{(k,l,n)\in \mathbb{Z}^3} c_{k,l,n} e^{2\pi i \left(k\omega+(l,n)\binom{x}{y}\right)}$$
be the Fourier series of $f$ on $\mathbb{T}^2\times \mathbb{T}$. Then
$$\sum_{(k,l,n)\in \mathbb{Z}^3} |c_{k,l,n}|^2=\|f\|_{L^2(\mu)}<+\infty.$$
Fix $(n',l')\in \mathbb{Z}^2$. For any $g\in L^2(m_{\mathbb{T}})$, where $m_{\mathbb{T}}$ is the Lebesgue measure on $\mathbb{T}$,
we have
\begin{align*}
&\hskip0.5cm \int \left( \sum_{k\in \mathbb{Z}} c_{k,l',n'} e^{2\pi i k \omega}\right) g(\omega) d m_{\mathbb{T}}(\omega)\\
&=\int \left(\sum_{(k,l,n)\in \mathbb{Z}^3} c_{k,l,n} e^{2\pi i \left(k\omega+(l,n)\binom{x}{y}\right)}\right) g(\omega)e^{-2\pi i(l',n')\binom{x}{y}} d \mu\left(\binom{x}{y},\omega\right) \\
&=\int f\left(\binom{x}{y},\omega\right) g(\omega)e^{-2\pi i(l',n')\binom{x}{y}} d \mu\left(\binom{x}{y},\omega\right)\\
&=\int f\left(\phi^{-1}\left(\binom{x}{y},\omega\right)\right)  g(\omega)e^{-2\pi i(l',n')\binom{x}{y}} d \mu\left(\binom{x}{y},\omega\right)  \\
&=\int f\left(\binom{x}{y},\omega\right) g(\omega+\alpha)e^{-2\pi i(l',n')\left(A\binom{x}{y}+h(\omega)\right)} d \mu\left(\binom{x}{y},\omega\right)\\
&= \int \left(\sum_{(k,l,n)\in \mathbb{Z}^3} c_{k,l,n} e^{2\pi i \left(k\omega+(l,n)\binom{x}{y}\right)}\right) g(\omega+\alpha)e^{-2\pi i(l',n')\left(A\binom{x}{y}+h(\omega)\right)} d \mu\left(\binom{x}{y},\omega\right)\\
&=\int\left( \sum_{k\in \mathbb{Z}} c_{k,(l',n')A} e^{2\pi i k \omega}\right)e^{-2\pi i (l',n')h(\omega)} g(\omega+\alpha) d m_{\mathbb{T}}(\omega)\\
&=\int \left( \sum_{k\in \mathbb{Z}} c_{k,(l',n')A} e^{2\pi i k (\omega-\alpha)}\right)e^{-2\pi i (l',n')h(\omega-\alpha)} g(\omega) d m_{\mathbb{T}}(\omega).
\end{align*}
Since the above equality holds for any $g\in L^2(m_{\mathbb{T}})$, we have
\begin{equation}\label{eq-11}
 \sum_{k\in \mathbb{Z}} c_{k,l',n'} e^{2\pi i k \omega}= \left( \sum_{k\in \mathbb{Z}} c_{k,(l',n')A} e^{2\pi i k (\omega-\alpha)}\right)e^{-2\pi i (l',n')h(\omega-\alpha)}
\end{equation}
in $L^2(m_{\mathbb{T}})$.
Thus
$$\int \left| \sum_{k\in \mathbb{Z}} c_{k,l',n'} e^{2\pi i k \omega}\right|^2 d m_{\mathbb{T}}(\omega)=
\int \left|\left( \sum_{k\in \mathbb{Z}} c_{k,(l',n')A} e^{2\pi i k (\omega-\alpha)}\right)e^{-2\pi i (l',n')h(\omega-\alpha)}\right|^2 d m_{\mathbb{T}}(\omega).
$$
Hence
\begin{equation*}
\sum_{k\in \mathbb{Z}} |c_{k,l',n'}|^2=\sum_{k\in \mathbb{Z}} |c_{k,(l',n')A}|^2.
\end{equation*}
When $(l',n')\neq (0,0)$, $(l',n')A^m\neq (l',n')A^k$ for $m,k\in \mathbb{N}$ with $m\neq k$.
Hence if $(l',n')\neq (0,0)$ then for any $N\in \mathbb{N}$,
\begin{align*}
\sum_{k\in \mathbb{Z}} |c_{k,l',n'}|^2=\frac{1}{N}\sum_{i=0}^{N-1}\sum_{k\in \mathbb{Z}} |c_{k,(l',n')A^i}|^2 \le \frac{1}{N} \sum_{(k,l,n)\in \mathbb{Z}^3} |c_{k,l,n}|^2\le \frac{1}{N} \|f\|_{L^2(\mu)}.
\end{align*}
Let $N\nearrow +\infty$, we get $\sum_{k\in \mathbb{Z}} |c_{k,l',n'}|^2=0$ for $(l',n')\neq (0,0)$.
When $(l',n')= (0,0)$, by \eqref{eq-11} we have
$$\sum_{k\in \mathbb{Z}} c_{k,0,0} e^{2\pi i k \omega}=\sum_{k\in \mathbb{Z}} c_{k,0,0}e^{-2\pi i k\alpha} e^{2\pi i k \omega}$$
in $L^2(m_{\mathbb{T}})$.
Thus
$c_{k,0,0}=c_{k,0,0}e^{-2\pi i k\alpha}$ for $k\in \mathbb{Z}$. When $k\neq 0$, we have $c_{k,0,0}=0$ as $\alpha \in \mathbb{R}\setminus \mathbb{Q}$.
Summing over we get $f=c_{0,0,0}$. This shows that $(\mathbb{T}^2\times \mathbb{T},\mathcal{B}_{\mathbb{T}^2\times \mathbb{T}}, \phi,\mu)$ is ergodic. The proof of this lemma is complete.

\end{proof}

For a continuous function $p$ from $\mathbb{T}$ to $\mathbb{T}$, we denote the degree of $p$ by $deg(p)$.

\begin{lem} \label{ex-lem-2} Let $deg(h_i)=n_i$ for $i=1,2$. If $n_2$ is odd number, then there are no $g\in C(\mathbb{T},\mathbb{T}^2)$  and positive integer $n$ such that
$\phi^n (graph(g))=graph(g)$, that is, $(\mathbb{T}^2\times \mathbb{T},\phi)$ has no random periodic point whose graph is continuous.
\end{lem}
\begin{proof} If this is not true, then there are  $g(\omega)=\binom{g_1(\omega)}{g_2(\omega)}\in C(\mathbb{T},\mathbb{T}^2)$  and positive integer $n$ such that
$\phi^n (graph(g))=graph(g)$, i.e.
$$ \phi^n(g(\omega),\omega)=(g(\omega+n\alpha),\omega+n\alpha)$$
for $\omega\in \mathbb{T}$.
Let $deg(g_i)=k_i$ for $i=1,2$. Note that
$$\phi^n(g(\omega),\omega)=\left(A^n\binom{g_1(\omega)}{g_2(\omega)}+A^{n-1}\binom{h_1(\omega)}{h_2(\omega)}+\cdots+\binom{h_1(\omega+(n-1)\alpha)}{h_2(\omega+(n-1)\alpha)},\omega+n\alpha\right)$$
Hence
\begin{align*}
&\binom{g_1(\omega+n\alpha)}{g_2(\omega+n\alpha)}\\=&A^n\binom{g_1(\omega)}{g_2(\omega)}
+A^{n-1}\binom{h_1(\omega)}{h_2(\omega)}+A^{n-2}\binom{h_1(\omega+\alpha)}{h_2(\omega+\alpha)}+\cdots+\binom{h_1(\omega+(n-1)\alpha)}{h_2(\omega+(n-1)\alpha)}.
\end{align*}
Compare the degree of continuous functions from $\mathbb{T}$ to $\mathbb{T}$ appearing in the above equation,  we have
$$\binom{k_1}{k_2}=A^n\binom{k_1}{k_2}
+A^{n-1}\binom{n_1}{n_2}+A^{n-2}\binom{n_1}{n_2}+\cdots+\binom{n_1}{n_2}.$$
That is
$$(I-A^n)\binom{k_1}{k_2}=(I+A+\cdots+A^{n-1})\binom{n_1}{n_2}.$$
Since $\det(I-A^n)\neq 0$ and $(I-A)(I+A+\cdots+A^{n-1})=I-A^n$,  we have
$$(I-A)\binom{k_1}{k_2}=\binom{n_1}{n_2}.$$
This implies $k_1=-\frac{n_2}{2}\not \in \mathbb{Z}$, a contradiction. The proof is complete.
\end{proof}
\begin{rem}\label{R:NonContiPerOrb}
Note that Theorem \ref{T:TheoryAnosovMix} implies that the random periodic orbits are dense in $L^{\infty}(\O, M)$.  Questions such as "Is there any continuous random periodic orbit?" or    "When do there exist continuous periodic orbits?" raise naturally. Lemma \ref{ex-lem-2} somehow answers the first question in a negative way; while the   Theorem 1.1  from \cite{Liu} give the positive answer for the first question when the system is generated by i.i.d. small random perturbations of an Axiom A system. In another aspect, Lemma \ref{ex-lem-2} tells that the systems we considered here beyonds the framework of the study on structure stability.
\end{rem}

\subsection{Systems preserving a volume on fibers driven by a force}\label{S:VolumePreserving}
In this section, we consider a type of  Anosov systems which are volume-preserving on fiber. We first state and prove the mixing property in a relatively general setting. Then,  we  investigate a concrete example which is generated by the random composition of matrices on 2-d torus, whose entries are positive integers, and  show such a system satisfies the Conditions 2A.-2C. below.

\subsubsection{Setting and results}\label{S:S1)SetResult}
 We assume that the systems  satisfy the following conditions:
\begin{itemize}
\item[2A.] $(\t,\O)$ is a heomomorphism  on a compact metric space;
\item[2B.] $\phi$ is Anosov on fibers (Condition H1));
\item[2C.] There exists an $f$-invariant Borel probability measure $\nu$ with full support  (i.e. $supp\nu=M$).
\end{itemize}
Here by a Borel probability measure $\nu$ on $M$ being {\bf $f$-invariant}  we mean that for any Borel measurable set $A\subset M$ and any $\o\in\O$, $\nu(f_\o^{-1}(A))=\nu (A)$.

 For the sake of convenience, we call the systems satisfying Conditions 2A.-2C. the S2) systems.
 The main result we will prove is the following Theorem:
 \begin{thm}\label{T:Spectrum12}
 S2) systems satisfity Condition H2) that is topological mixing on fibers.
 \end{thm}
\begin{proof}
Unlike the case of S1) systems, the proof for S2) is based on the measure preserving property of $f_\o$, which leads to a much more direct approach.

 For a given $\o\in \O$, an open subset $V\subset M$, and a  real number $r>0$, we inductively define a sequence of subsets $V^r_n(\o)\subset M$, $n=0,1,2,\cdots$ as follows:
 \begin{equation}\label{E:V_n}
 V^r_n(\o)=\begin{cases}
 V,&n=0\\
 \cup_{x\in f_{\t^{n-1}\o}(V^r_{n-1}(\o))}B_M(x,r),&n>1
 \end{cases},
 \end{equation}
 where $B_M(x,r)$ is the ball in $M$ with radius $r$ centered at $x$. It is not hard to see that each $V^r_n(\o)$ is open and $\{\nu(V^r_n(\o))\}_{n=0,1,\cdots}$ is non-decreasing  on both $n$ and $r$ since $\nu$ is $f$-invariant and $f$ is invertible.  Thus $\lim_{n\to\infty}\nu(V^r_n(\o))$ exists.

 Next, we show that
 \begin{equation}\label{E:LimVn}
  \forall\ r>0, \exists\ N(\o)>0 \text{ such that }\forall\ n\ge N(\o), V^r_n(\o)=M.
 \end{equation}
 Since $supp\ \nu=M$, we know that any open subset of $M$ is not a $\nu$-null set. Additionally, because $M$ is compact, for any $R>0$, $M$ can be  covered by finitely many balls with radius $\frac13R$, each of which has a positive $\nu$-measure. Thus, \begin{equation}\label{E:NuInf}
 \forall\ R>0, \exists\ \e>0, \text{ such that }\forall\ x\in M,\ \nu(B_M(x,R))>\e.
 \end{equation}
 For a given $n\ge 1$, let $t^r_n(\o)=\min\{\sup\{d_H(x,f_{\t^{n-1}\o}(V^r_{n-1}(\o))|\ x\in M\},r\}$, where for $x\in M$ and $A\subset M$  $$d_H(x,A):=\inf\{d_M(x,y)|\ y\in A\}$$ is the Hausdoff distance. Since $d_H(x,V^r_n(\o))$ is continuous on $x$ and $M$ is connected, we have that there exists $x'\in M$ such that $$\frac12t^r_n(\o)< d_H(x',f_{\t^{n-1}\o}(V^r_{n-1}(\o)))<r.$$
 Therefore, $B\left(x',\frac12t^r_n(\o)\right)\subset V^r_n(\o)\setminus f_{\t^{n-1}\o}(V^r_{n-1}(\o))$ which implies that
 $$\nu \left(B_M\left(x',\frac12t^r_n(\o)\right)\right)\le \nu(V^r_n(\o))-\nu(f_{\t^{n-1}\o}(V^r_{n-1}(\o)))=\nu(V^r_n(\o))-\nu(V^r_{n-1}(\o)(\o)).$$
 Because of the existence of  $\lim_{n\to\infty} \nu (V^r_n(\o))$ and (\ref{E:NuInf}), we have that
 $$\lim_{n\to\infty} t_n^r=0.$$
 By the definition of $V_n^r(\o)$ as in (\ref{E:V_n}), once $t_n^r<r$, $V_n^r(\o)=M$, thus (\ref{E:LimVn}) holds. \\

 Now, we are ready to show that S2)-type systems are topological mixing on fibers.

 Let $U_1, U_2\subset M$ be two nonempty open sets. Choose two open balls $B_M(x_i,3R)\subset U_i$, $i=1,2$, with $R>0$. Applying Lemma \ref{L:Shadowing} by setting $\b=R$, we have an $\a>0$, such that any $(\o,\a)$-pseudo orbit can be $(\o,\b)$-shadowed by a true orbit.

  For a given $\o\in\O$, take $V=B_M(x_1,R)$ and $r=\a$, and define $V^r_n(\o)$ as in (\ref{E:V_n}). By (\ref{E:LimVn}), there exists an $N_1(\o)$ such that for any $n\ge N_1(\o)$, $V_n^r(\o)=M$. By continuity of $\phi$, we have that for any $\o\in\O$, there exists $\d(\o)>0$ such that the following holds:
 $$\forall\ \o'\in B_\O(\o,\d(\o))\text{ and }n\ge N_1(\o)+1, \ V^r_n(\o') =M,$$
 where $B_\O(\o,\d(\o))$ is the ball in $\O$ with radius $\d(\o)$ centered at $\o$.

 Since $\O$ is compact, and $\{B_\O(\o,\d(\o))\}_{\o\in\O}$ forms an open cover on $M$, there exists a finite sub-cover, $\{B_\O(\o_j,\d(\o_j)\}_{1\le j\le q},\ q\in\mb N$. Take $N=\max\{N_1(\o_1),\cdots, N_1(\o_q)\}+1$, then the following holds
 \begin{equation}\label{E:UniformN}
  \forall\ \o\in\O, \text{ if } n\ge N, \text{ then } V^r_n(\o)=M.
  \end{equation}
 Next, we show that
 \begin{equation}\label{E:TMix}
\forall\ \o\in\O\text{ and }n\ge N,\  \phi(n,U_1\times\{\o\})\cap U_2\times\{\t^n\o\}\neq \emptyset.
 \end{equation}
 This can be done by constructing an $(\o,r)$-pseudo orbit, $\{(y'_{n,i},\t^i\o)\}_{i\in\mb Z}$, connecting two balls $B_M(x_1,R)\times \{\o\}$ and $B_M(x_2,R)\times \{\t^n\o\}$, which is $(\o,R)$-shadowed by a true orbit \[\{(y_{n,i},\t^i\o)\}_{i\in\mb Z}.\] Thus this true orbit,  $\{(y_{n,i},\t^i\o)\}_{i\in\mb Z}$, will hit both $B_M(x_1,3R)\times \{\o\}$ and $B_M(x_2,3R)\times \{\t^n\o\}$. Hence (\ref{E:TMix}) holds. Actually, the existence of such an $(\o,r)$-pseudo orbits ($\forall\ \o\in \O$), $\{(y'_{n,i},\t^i\o)\}_{i\in\mb Z}$, follows from (\ref{E:UniformN}), which can be constructed in the following:

 For any $n\ge N$ and $\o\in\O$, since $V_N^r(\o)=M$, we have that there exists $z_{N-1}(\o)\in V_{N-1}^r(\o)$ such that
 $$\pi_M\phi^{-(n-N)}(x_2,\t^n\o)\in B_M(f_{\t^{N-1}\o}z_{N-1}(\o),r).$$
 Inductively, we have that for any $i\in\{0,\cdots, N-2\}$, there exists $z_i(\o)\in V_i^r(\o)$ such that
 $$z_{i+1}(\o)\in B_M(f_{\t^i\o}z_i(\o),r).$$
 Also note that $V_0^r(\o)=B_M(x_1,R)$, then we have that $\{z_i(\o)\}_{1 \le i\le N-1}$ satisfies the following properties:
 \begin{align}\begin{split}\label{E:ConnectX1X2}
  &d_M(z_0(\o),x_1)<r;\\
  &d_M(z_{i+1}(\o), f_{\t^i\o}z_i(\o))<r,\ 0\le i\le N-2;\\
  &d_M\left(\pi_M\phi^{-(n-N)}(x_2,\t^n\o),f_{\t^{N-1}\o}z_{N-1}(\o)\right)<r.
  \end{split}
   \end{align}
We define
$$
y'_{n,i}=\begin{cases}
\pi_M\phi^i(x_1,\o),&i\le -1\\
z_i(\o),&0\le i\le N-1\\
\pi_M\phi^{i-n}(x_2,\t^n\o), &i\ge N
\end{cases}.
$$
By (\ref{E:ConnectX1X2}), it is not hard to see that the above $y'_{n,i}$'s  fit all the requirements. The proof is  complete.

\end{proof}


\subsubsection{Random composition of $2\times 2$ area-preserving positive matrices } \label{S:Example2}
The main purpose of this section is to provide an example of  S2)-type systems, which is generated by the random composition of a class of area-preserving positive $2\times 2$ matrices. Then, we apply the results for S2)-type systems to investigate the dynamical behavior of the given example.

Let $$\left\{A_i=\left(\begin{matrix}a_i&b_i\\c_i&d_i \end{matrix}\right)\right\}_{1\le i\le p}$$ be $2\times 2$ matrices with positive integer entries and
$|\det A_i|=1,\ \forall i\in\{1,\cdots,p\}$, and $\mc S_p:=\{1,\cdots, p\}^{\mb Z}$ with the left shift operator $\s$ be the  symbolic dynamical system with $p$ symbols.

Define a map $f:\mc S_k\to \{A_1\cdots,A_p\}$ by letting
$$f(\hat a)=A_{\hat a(0)},\ \forall \hat a\in \mc S_k.$$

Define that, for $x\in [0,+\infty)\times [0,+\infty)$ and $\hat a\in \mc S_k$,
$$\phi^n_{\hat a}=\begin{cases}
I_2&\text{ when }n=0\\
f(\s^{n-1}\hat a)\circ\cdots\circ f(\hat a)&\text{ when }n>0\\
f^{-1}(\s^{n}\hat a)\circ\cdots\circ f^{-1}(\s^{-1}\hat a)&\text{ when }n<0\\
\end{cases},
$$
where $I_2=\left(\begin{matrix}1&0\\0&1\end{matrix}\right)$.
Note that $\phi$ induces a random dynamical system on 2 dimensional torus  $\mb T^2$, which is denoted by $\hat \phi$. Let \begin{equation}\label{E:AssumpExpan}
\kappa:=\min_{1\le i\le p}\min\left\{\sqrt{a_i^2+c_i^2},\sqrt{b_i^2+d_i^2}\right\},
\end{equation}
which is obviously greater than or equal to $\sqrt 2$.

\begin{prop}\label{P:2DTorusAnosov}
$\hat \phi:\mb T^2\times \mc S_k\to \mb T^2\times \mc S_k$ is an S2)-type system by letting $M=\mb T^2$ and $(\O,\t)=(\mc S_p,\s)$.
\end{prop}
\begin{proof}
Denote $\mb R^{2+}:=\left\{(x,y)\in \mb R^2\setminus \{(0,0)\}\Big|\ x\ge 0, y\ge 0\right\}$. For any $v_i=\left(\begin{matrix}x_i\\y_i\end{matrix}\right)\in \mb R^{2+}$, $i=1,2$, by a straightforward computation, we have that
 $$\sin\angle(v_1,v_2)=\left| \frac{y_2x_1-y_1x_2}{\sqrt{x_1^2+y_1^2}\sqrt{x_2^2+y_2^2}}\right|,$$
 where $ \angle(v_1,v_2)$ is the angle between $v_1$ and $v_2$. \\

 By a straightforward computation, we have that, for any $j\in\{1,\cdots,p\}$, $i=1,2$,
 \begin{align}\begin{split}\label{E:Expan&Contrac}
& |A_jv_i|=\sqrt{(a_j^2+c_j^2)x_i^2+2(a_jb_j+c_jd_j)x_iy_i+(b_j^2+d_j^2)y_i^2}\ge \kappa |v_i|\\
&\sin\angle(A_jv_1,A_jv_2)=\left|\frac{(a_jd_j-b_jc_j)(y_2x_1-y_1x_2)}{|A_jv_1||A_jv_2|}\right|\le \kappa^{-2}\sin\angle(v_1,v_2).
\end{split}
 \end{align}
Let $Proj^+:=\left\{\spa\{v\}|\ v\in \mb R^{2+}\right\}$ which is a compact subset of the projective space $\mb P^2$.  It is not hard to see that the second inequality of (\ref{E:Expan&Contrac}) yields that for any $\hat a\in \mc S_k$ and $n\in\mb N$ the map $\phi(\cdot,\hat a)$ induces a  contraction map on $Proj^+$ with a uniform contracting rate. By the contracting mapping theorem, we have that there is a unique $1$-dimensional space $E_{\hat a}$ such that
$$E_{\hat a}\in \bigcap_{n\ge 0}\bar\phi^n(Proj^+,\s^{-n}\hat a),$$
where we denote $\bar \phi(\cdot,\hat a)$ the induced map on $\mb P^2$ by  $\phi(\cdot,\hat a)$. By the uniqueness, we have that
$$\phi(E_{\hat a},\hat a)=E_{\s\hat a}.$$
The first inequality of (\ref{E:Expan&Contrac}) yields that $\phi(|_{E_{\hat a}},\hat a)$ is expanding with rate at least $\kappa$. The continuity of $E_{\hat a}$ on $\hat a$ simply follows from the second inequality of (\ref{E:Expan&Contrac}).\\

By applying exactly the  same argument, we can obtain that for any $\hat a\in \mc S_k$ there exists a unique $1$-dimensional space $F_{\hat a}\notin Proj^+$ such that $F_{\hat a}$ varies continuously on $\hat a$ and the following holds
$$\phi(F_{\hat a},\hat a)=F_{\s\hat a}\text{ and }\phi^{-1}(|_{F_{\hat a}},\hat a) \text {is expanding with rate at least }\kappa.$$
Note that $\phi$ itself is linear on $\mb T^2$, therefore the linearization of $\hat \phi^n$, $D_{\mb T^2}\hat \phi^n$, can be identified by $\phi$. Hence, $\hat \phi$ is an S2)-type system with continuous co-invariant splitting $\mb R^2=E_{\hat a}\oplus F_{\hat a}$ for $\hat a\in\mc S_k$.
\end{proof}
The following corollaries are the consequences of Proposition \ref{P:2DTorusAnosov} and Theorem \ref{T:TheoryAnosovMix}.
\begin{cor}\label{C:FrequentReturn}
For any $\e>0$, there exists $N\in\mb N$ such that for any $x\in \mb T^2$, $n\ge N$, and $\hat a\in \mc S_k$, there exists $x'\in \mb T^2$ satisfying
$$x',\ \pi_{\mb T^2}\hat \phi^n(x',\hat a)\in B_{\mb T^2}(x,\e).$$
\end{cor}

\begin{cor}\label{C:T^2Horseshoe}
For any $\g>0$ there exist $N,k$ such that the following hold:

\begin{itemize}
\item[i)] $\frac1N\log k>\log \kappa-\g$;
\item[ii)] $\hat \phi^N$ has a strong full random horseshoe with $k$ symbols.
\end{itemize}
Here $\kappa$ is given by (\ref{E:AssumpExpan}).
\end{cor}

\subsection{Further discussions on Random Horseshoes with marginal $\mb P$}\label{S:RandomHorseshoe}
In this section, we further discuss  the random horseshoes for systems with a given marginal $\mb P$ which is an $\t$-invariant probability measure on $\O$. The standard fiber topological entropy of $\phi$ with respect to $\mb P$ is defined by the following
\begin{equation}\label{E:FibTopEntrP}
h_{top}(M\times\O|\ \mb P)=\int_{\O}h_{top}(\phi|_{M_\o})d\mb P(\o).
\end{equation}
Since $h_{top}(\phi|_{M_\o})=h_{top}(\phi|_{M_{\t\o}})$ for $\o\in\O$ by the definition, we have that when $\mb P$ is ergodic the following holds
$$h_{top}(M\times\O|\mb P)\equiv h_{top}(\phi|_{M_{\o}})\text{ for } \mb P-a.e.\ \o.$$
For the random dynamical system $\phi$ over a metric dynamical system $(\O,\mb P,\t)$ mentioned in Remark \ref{R:RandomDynSys}, it is natural to modify the definitions of random horseshoes correspondingly.

Let $L^\infty_{\mb P}(\O)$ be the space of Borel measurable maps from $\O$ to $M$ endowed with the    metric $d_{L^\infty_\mb P(\O)}$:
$$d_{L^\infty_\mb P(\O)}(g_1,g_2)=\sup_{\mb P}\left\{d_M(g_1(\o),g_2(\o)),\o\in\O\right\},\ \forall g_1,g_2\in L^\infty(\O, M),$$
where $\sup_{\mb P}$ is the essential supreme with respect to $\mb P$. Similarly, $\phi:M\times\O$ induces a homeomorphism on $L^\infty_{\mb P}(\O)$, which is denoted by $\tilde\phi_\mb P$.
Analogically, we define two functions on $L^\infty_\mb P(\O)\times  L^\infty_\mb P(\O)$ to measure separations of elements of $L^\infty_\mb P(\O)$ on different levels:
\begin{align}\begin{split}\label{E:PSeparationFunc}
\overline d_{L^\infty_\mb P(\O)}(g_1,g_2)&=\sup_{\O'\subset \O,\ \mb P(\O')>0} \inf_{\o\in \O'}d_M(g_1(\o),g_2(\o)),\ \forall g_1,g_2\in L^\infty_\mb P(\O),\\
\underline d_{L^\infty_\mb P(\O)}(g_1,g_2)&=\inf_{\mb P}\left\{d_M(g_1(\o),g_2(\o)),\ \o\in\O\right\},\ \forall g_1,g_2\in L^\infty_\mb P(\O),
\end{split}
\end{align}
where $\inf_\mb P$ is the essential infimum  with respect to $\mb P$. First, we claim that
$$\overline d_{L^\infty_\mb P(\O)}=d_{L^\infty_\mb P(\O)}.$$
For the sake of completeness, we give the proof here. Recall that
$$d_{L^\infty_\mb P(\O)}(g_1,g_2)=\inf_{\O'\subset \O,\ \mb P(\O')=1} \sup_{\o\in \O'}d_M(g_1(\o),g_2(\o)),\ \forall g_1,g_2\in L^\infty_\mb P(\O).$$
For a given $\e>0$, let
$$\O'_\e(g_1,g_2)=\left\{\o\in\O|\ d_M(g_1(\o),g_2(\o))>d_{L^\infty_\mb P(\O)}(g_1,g_2)-\e\right\}.$$
By the definition of $d_{L^\infty_\mb P(\O)}(g_1,g_2)$, we have that $\mb P(\O'_\e(g_1,g_2))>0$. Thus
$$\overline d_{L^\infty_\mb P(\O)}(g_1,g_2)\ge \inf_{\o\in \O'_\e(g_1,g_2)}d_M(g_1(\o),g_2(\o))>d_{L^\infty_\mb P(\O)}(g_1,g_2)-\e.$$
Let $$\O''_\e(g_1,g_2)=\left\{\o\in\O|\ d_M(g_1(\o),g_2(\o))>d_{L^\infty_\mb P(\O)}(g_1,g_2)+\e\right\}.$$
By the definition of $d_{L^\infty_\mb P(\O)}(g_1,g_2)$, we have that $\mb P(\O''_\e(g_1,g_2))=0$. Thus
$$\mb P\left(\O'\cap (\O\setminus \O''_\e(g_1,g_2))\right)>0,\ \forall\ \O'\subset \O,\ \mb P(\O')>0,$$
which implies that
$$\overline d_{L^\infty_\mb P(\O)}(g_1,g_2)\le  d_{L^\infty_\mb P(\O)}(g_1,g_2)+\e.$$
By the arbitrariness of $\e$, we complete the proof of the claim. In the rest of this section, we will only use $d_{L^\infty_\mb P(\O)}(g_1,g_2)$ instead of $\overline d_{L^\infty_\mb P(\O)}(g_1,g_2)$.

Consequently, we can define the random horseshoes and strong random horseshoes with marginal $\mb P$.

 \begin{defn}\label{D:PRandomHorseshoe}
 We  call a continuous embedding $\Psi:\mc S_k\to L^\infty_\mb P(\O)$ a {\bf \em  random horseshoe on marginal $\mb P$} with $k$-symbols of $\phi$ if
the following hold:
\begin{itemize}
     \item[i)] There exists $\Delta>0$ such that for any $l\in\mb Z$ and $\hat a_1,\hat a_2\in \mc S_k\text{ with }\hat a_1(l)\neq \hat a_2(l)$
     \begin{equation}\label{E:RandomSepaImgSymb}
      d^{\tilde \phi_\mb P}_{l}\left(\Psi({\hat a_1}),\Psi({\hat a_2})\right)>\Delta,
     \end{equation}
     where for $l\ge 0$
     $$ d^{\tilde \phi_\mb P}_{l}(g_1,g_2):=\max_{0\le i\le l}\left\{ d_{L^\infty_\mb P(\O)}\left(\tilde \phi_\mb P^i(g_1),\tilde \phi^i_\mb P(g_2)\right)\right\},\ g_1,g_2\in L^\infty_\mb P(\O)$$  and for $l<0$ we let $ d^{\tilde\phi_\mb P}_{l}= d^{\tilde\phi^{-1}_\mb P}_{|l|};$
\item[ii)]  For all $\hat a\in\mc S_k$ $$\tilde \phi_\mb P\left(\Psi({\hat a})\right)=\Psi({\s \hat a}).$$
\end{itemize}
\end{defn}

By using an stronger separation function $\underline d_{L^\infty_\mb P(\O)}$ as in (\ref{E:PSeparationFunc}), we define a stronger random horseshoe as follows:
 \begin{defn}\label{D:PSRandomHorseshoe}
 We  call a continuous embedding $\Psi:\mc S_k\to L^\infty_\mb P(\O)$ a {\bf \em  strong  random horseshoe on marginal $\mb P$} with $k$-symbols of $\phi$ if
the following hold:
\begin{itemize}
     \item[i)] There exists $\Delta>0$ such that for any $l\in\mb Z$ and $\hat a_1,\hat a_2\in \mc S_k\text{ with }\hat a_1(l)\neq \hat a_2(l)$
     \begin{equation}\label{E:SRandomSepaImgSymb}
     \underline d^{\tilde \phi_\mb P}_{l}\left(\Psi({\hat a_1}),\Psi({\hat a_2})\right)>\Delta,
     \end{equation}
     where for $l\ge 0$
     $$\underline d^{\tilde \phi_\mb P}_{l}(g_1,g_2):=\max_{0\le i\le l}\left\{\underline d_{L^\infty_\mb P(\O)}\left(\tilde \phi^i_\mb P(g_1),\tilde \phi^i_\mb P(g_2)\right)\right\},\ g_1,g_2\in L^\infty(\O, M)$$  and for $l<0$ we let $\underline d^{\tilde\phi_\mb P}_{l}=\underline d^{\tilde\phi^{-1}_\mb P}_{|l|};$
\item[ii)]  For all $\hat a\in\mc S_k$ $$\tilde \phi_\mb P\left(\Psi({\hat a})\right)=\Psi({\s \hat a}).$$
\end{itemize}
\end{defn}
For the random horseshoe on marginal $\mb P$, we have the following result.
\begin{thm}\label{T:PHorseshoe}
 Let $\phi$ satisfy H1) and H2) and $\mb P$ be a $\t$-ergodic probability measure. Then for any $\g>0$,  there exist $N,k\in \mb N$ such that the following hold
   \begin{enumerate}
     \item[i)] $\frac1N\ln k\ge h_{top}(M\times\O|\mb P)-\g$;
     \item[ii)] $\phi^N$ has a $k$-symbol random horseshoe on marginal $\mb P$.
   \end{enumerate}
 \end{thm}
\begin{proof}
The proof of Theorem \ref{T:PHorseshoe} can be derived by employing the exact same paradigm used in the proof of Part B. of Theorem \ref{T:TheoryAnosovMix2}. To adapt the proof of Part B. of Theorem \ref{T:TheoryAnosovMix2} into the current circumstance, one  needs only to replace $\overline h(M\times\O|\O)$ by $h_{top}(M\times\O|\mb P)$, and make the "$O$" in Lemma \ref{L:SeparatedSetFiber} satisfy $\mb P(O)>0$ rather than "nonempty open subset", which is guaranteed by the definition of $h_{top}(M\times\O|\mb P)$. We omit the detailed proof here.
\end{proof}
For the strong random horseshoe on marginal $\mb P$, only partial result is obtained which is applicable for S1) systems. Note that for S1) systems, $\t:\O\to\O$ is uniquely ergodic, thus $\mb P$ can only be the unique Haar measure of $\t$.
\begin{thm}\label{T:SPHorseshoe}
 Let $\phi$ be an S1) system and $\mb P$ be the unique $\t$-ergodic probability measure. Then for any $\g>0$,  there exist $N,k\in \mb N$ such that the following hold
   \begin{enumerate}
     \item[i)] $\frac1N\ln k\ge h_{top}(M\times\O|\mb P)-\g$;
     \item[ii)] $\phi^N$ has a $k$-symbol strong random horseshoe on marginal $\mb P$.
   \end{enumerate}
 \end{thm}
 \begin{proof}
 This theorem follows from Part A. of Theorem \ref{T:TheoryAnosovMix2} and  the following lemma:
 \begin{lem}\label{L:S1EntropyEquality}
 Let $\phi$ be an S1) system. Then the following holds
 $$\underline h(M\times\O|\O)=\overline h(M\times\O|\O)=h_{top}(\phi),$$
 where $h_{top}(\phi)$ is the topological entropy of $\phi$.
 \end{lem}

 \begin{proof}[Proof of Lemma \ref{L:S1EntropyEquality}]
It is straightforward to observe that $\overline h(M\times\O|\O)\le h_{top}(\phi)$ by the definition. For the sake of simplicity, we will prove this lemma in the case of $h_{top}(\phi)<\infty$ only, as for the other case the  same argument is still applicable.
By the Theorem 17 in \cite{Bow3}, we have  that
$$h_{top}(\phi)\le \overline h(M\times\O|\O)+h_{top}(\t).$$
Since $\t$ is equicontinuous, we have $h_{top}(\t)=0$, thus
$$h_{top}(\phi)= \overline h(M\times\O|\O).$$
Therefore, it remains to show that $\underline h(M\times \O|\O)=\overline h(M\times \O|\O)$. It is sufficient to prove the following
\begin{equation}\label{E:FibEntrEqu}
h_{top}(\phi|_{M_{\o_1}})=h_{top}(\phi|_{M_{\o_2}}),\ \forall \o_1,\o_1\in\O.
\end{equation}
Recall that
$$h_{top}(\phi|_{M_\o})=\lim_{\e\to 0^+}\limsup_{n\to\infty}\frac1n\log N(\o,\e,n),$$
where
$$N(\o,\e,n):=\max\  \left\{card(E)\big|\ E\subset M_{\o} \text{ is an }\e\text{ separated set with respect to }d^{\phi}_{M\times \O,n}\right\},$$
and $d^\phi_{M\times \O,n}$ is the Bowen metric. For a given small positive number $\g>0$ and $\o_1\in\O$, there exists $\e_0>0$ such that for any $\e\in(0,\e_0)$ and  $N_0\in\mb N$, there is an $n\ge N_0$ so that the following holds
\begin{equation}\label{E:NumSepaSet}
N(\o_1,3\e,n)\ge e^{n( h_{top}(\phi|_{M_{\o_1}})-\frac12\g)}.
\end{equation}
Without loss of generality, we assume that $\e<\frac12\b_0$ and fix such an $\e$, where $\b_0$ is as in Lemma \ref{L:Shadowing}. Letting $\b=\frac12\e$ and applying Lemma \ref{L:Shadowing}, there exists $\a>0$ corresponding to $\b$ as in Lemma \ref{L:Shadowing}. Again, for the sake  of convenience, we assume that $\a<\b$. By the uniform continuity of $\phi$, there exists $\d>0$ such that for any $\o',\o''\in\O$ with $d_\O(\o',\o'')<\d$ and any $x\in M$ the following holds
\begin{equation}\label{E:VarEst}
d_M(\pi_M\phi(x,\o'),\pi_M\phi(x,\o''))<\a.
\end{equation}
We note that this $\d$ depends only on $\a$.

Let $E_{\o_1}\subset M_{\o_1}$ be an $3\e$-separated set with respect to $d^{\phi}_{M\times \O,n}$ satisfying $card(E_{\o_1})=N(\o_1,3\e,n)$. For $(x_j,\o_1)\in E_{\o_1}(=\{(x_j,\o_1)\}_{1\le j\le N(\o_1,3\e,n)})$ and $\o\in B_\O(\o_1,\d)$, define a pseudo orbit $\{(y_i'(x_j,\o),\o)\}_{i\in\mb Z}$ by letting
\begin{equation}\label{E:PsOrbSepa}
y_i'(x_j,\o)=\begin{cases}
\pi_M(\phi^i(x_j,\o_1)), &\text{ when }i\ge 0\\
\pi_M(\phi^i(x_j,\o)), &\text{ when }i< 0
\end{cases}.
\end{equation}
By the choice of $\d$ and  (\ref{E:VarEst}), it is not hard to see that $\{(y_i'(x_j,\o),\o)\}_{i\in\mb Z}$ is an $(\o,\a)$-pseudo orbit. Thus, by Lemma \ref{L:Shadowing}, there exists a true orbit, denoted by $\{\phi^i(y_j,\o)\}_{i\in\mb Z}$, which is $(\o,\b)$-shadowing $\{(y_i'(x_j,\o),\o)\}_{i\in\mb Z}$. Since $E_{\o_1}$ is a $3\e$-separated set with respect to $d^{\phi}_{M\times \O,n}$ and $\b=\frac12\e$, we have that $\{(y_j,\o)\}_{1\le j\le N(\o_1,3\e,n)}\subset M_\o$ is a $2\e$-separated set with respect to $d^{\phi}_{M\times \O,n}$.

Also note that, $(\O,\t)$ is uniquely ergodic, thus there is an $N_1\in\mb N$ depending on $B_\O(\o_1,\d)$ only, which satisfies that for any $\o\in\O$ there exists a positive integer $m(\o)\le N_1$ such that $\t^{m(\o)}(\o)\in B_\O(\o_1,\d)$. This implies that for any $\o\in\O$, the following holds
$$N(\o,2\e,n+N_1)\ge card (\phi^{-m(\o)}E_{\o_1})=N(\o_1,3\e,n)\ge e^{n( h_{top}\left(\phi|_{M_{\o_1}})-\frac12\g\right)}.$$
Therefore, we have that
\begin{align*}
h_{top}(\phi|_{M_\o})&=\lim_{\e'\to 0^+}\limsup_{m\to\infty}\frac1m\log N(\o,\e',m)\\
&\ge \frac1{n+N_1}\log N(\o,2\e,n+N_1)\\
&= \frac n{n+N_1}\left( h_{top}(\phi|_{M_{\o_1}})-\frac12\g\right).
\end{align*}
Since  such $n$ can be arbitrarily large and $\g$ can be arbitrarily small, we have that
$$ h_{top}(\phi|_{M_\o})\ge h_{top}(\phi|_{M_{\o_1}}).$$
Also note that $\o_1$ is arbitrarily chosen, thus the proof of (\ref{E:FibEntrEqu}) is complete.
 \end{proof}
 To end the proof, it suffices to note that
 $$h_{top}(M\times\O|\mb P)=\overline h(M\times\O|\O)=\underline h(M\times\O|\O)=h_{top}(\phi)$$
 and
 $$\underline d_{L^\infty(\O, M)}\le \underline d_{L^\infty_\mb P(\O)}.$$
 \end{proof}

 At the end of this section, we make the following remark on the random SRB measures. We summarize the notions and results of the random SRB measures in Section \ref{S:RandomSets} for the sake of completeness, which are mainly borrowed from \cite{GK}.
  \begin{rem}\label{R:SRBMeasure}
Under the setting of this paper, S1) and S2) systems satisfy the conditions of Theorem \ref{T:SRBMeasure}, thus, roughly speaking, we have both existence and uniqueness of SRB measure for a given $\t$-invariant ergodic probability measure $\mb P$.

For S1) systems, since the driven system $(\O,\t)$ is uniquely ergodic, each S1) system has a unique SRB measure $\mu$. By the variational principle, we have that $$h_\mu(\phi)\le h_{top}(\phi).$$
Therefore, by Theorem \ref{T:SPHorseshoe}, there exist strong random horseshoes on marginal $\mb P$ whose entropy can approach  $h_\mu(\phi)$.

For S2) systems, for any $\t$-ergodic measure $\mb P$, there exists a unique SRB measure $\mu$ whose marginal is $\mb P$. By the variational principle, we have that
$$h_{\phi}(\mu|\mb P)\le  h_{top}(M\times\O|\ \mb P), $$
$h_\phi(\mu|\mb P)$ is called the relative entropy of $\mu$ with respect to $\mb P$.
Therefore, by Theorem \ref{T:PHorseshoe}, there exist  random horseshoes on marginal $\mb P$ whose entropy can approach $h_\phi(\mu|\mb P)$.
\end{rem}


\appendix

\section{Random Sinai-Ruelle-Bowen measures}\label{S:RandomSets}
In this section, some notions and results will be given for random dynamical systems, especially the results of Random SRB measures, most of which are borrowed from \cite{GK}.

Let $(\O,\mc B,\mb P)$ be a probability space, $\mb E$ be a locally compact Hausdorff second countable topological space. Let $\mc F, \mc K,\mc G$ denote respectively the family of all closed, compact and nonempty open sets of $\mb E$.
\begin{defn}\label{D:RandomSet}
A map $X:\O\to\mc F$ is called a {\em random closed set} if, for every compact set $K\subset \mb E$, $$\{\o\in\O|X(\o)\cap K\neq \emptyset\}\in \mc B.$$

A map $Y:\O\to\mc G$ is called a {\em random open set} if its complement $X=Y^c$ (by this, we mean $X(\o)=E\setminus Y(\o)$ for $\mb P$-a.e. $\o\in\O$) is a random closed set. Since $\{\o\in\O|Y^c(\o)\cap F=\emptyset\}=\{\o\in\O|F\subset Y(\o)\}$, $Y$ is random open set if and only if $\{\o\in\O|F\subset Y(\o)\}$ is a measurable event for every $F\in \mc F$.
\end{defn}
Let $\t:(\O,\mc B,\mb P)\to(\O,\mc B,\mb P)$ be an invertible ergodic metric dynamical system, and $F:\mb E\times \O\to \mb E\times \O$ be a continuous random dynamical system over this metric dynamical system.

\begin{defn}\label{D:RTran}
 $F$ is called {\em random topological transitive} if for any random open sets $U,V$ with $U(\o),V(\o)\neq \emptyset$ for all $\o\in\O$, there exists a random variable $n$ taking values in $\mb Z$ such that the intersection $F^{n(\o)}(\t^{-n(\o)}\o,U(\t^{-n(\o)}\o))\cap V(\o)$ is non-empty $\mb P$-a.s..
 \end{defn}

The following is a technical lemma whose proof is included for the sake  of completeness.
\begin{lem}\label{L:MixToTran}
  $F$ is topological mixing on fibers (see Section \ref{S:Setting} for definition) implies that $F$ is random topological transitive.
\end{lem}
\begin{proof}
Let $U$ and $V$ be random open sets with $U(\o),V(\o)\neq \emptyset$ for $\mb P$-a.e. $\o\in\O$. Let $\{x_i\}_{ i\in\mb N}$ be a countable dense set of $\mb E$. For $i,j\in\mb N$, define
$$\O^U_{i,j}=\left\{\o|\ \overline{B\left(x_i,\frac1j\right)}\subset U(\o)\right\}\text{ and }\O^V_{i,j}=\left\{\o|\ \overline{B\left(x_i,\frac1j\right)}\subset V(\o)\right\}.$$
It is not hard to see that each $\O^U_{i,j}$ or $\O^V_{i,j}$ is measurable, and $\mb P(\cup_{i,j}\O^U_{i,j})=\mb P(\cup_{i,j}\O^V_{i,j})=1$. Without loss of generality, we suppose that $\mb P(\O^U_{i_0,j_0})>0$ for some $i_0,j_0\in\mb N$. Since $F$ is topological mixing on fibers, for any $i,j\in \mb N$, there exists $N_{i,j}\in\mb N$ such that for any $k>N_{i,j}$
\begin{equation}\label{E:Hitting}
F^k\left(B\left(x_{i_0},\frac1{j_0}\right)\times \{\t^{-k}\o\}\right)\cap B\left(x_i,\frac1j\right)\times\{\o\} \neq \emptyset,\ \forall \o\in\O.
\end{equation}
By the ergodicity of $\t^{-1}$ and Poincare Recurrence Theorem, we have that for $\mb P$-a.s. $\o\in\O$, $\{\t^{-k}\o\}_{k\in\mb N}$ will visit $\O^U_{i_0,j_0}$ infinitely many  times. Given $i,j\in\mb N$, for $\o\in \O^V_{i,j}$, define that
$$n(\o)=\min\{k|\ k>N_{i,j},\ \t^{-k}\o\in \O^U_{i_0,j_0}\}.$$
It is not hard to see that $n:\O_{i,j}^V\to\mb N$ is measurable and $n(\o)<\infty$ for $\mb P$-a.s. $\o$.  Additionally, (\ref{E:Hitting}) implies that for $\mb P$-a.s. $\o\in \O^V_{i,j}$, the following holds
\begin{align}\begin{split}\label{E:TTran}
&F\left(n(\o),U(\o)\times\left\{\t^{-n(\o)}\o\right\}\right)\cap V(\o)\times\{\o\}\\
\supset&F\left(n(\o),B\left(x_{i_0},\frac1{j_0}\right)\times\left\{\t^{-n(\o)}\o\right\}\right)\cap B\left(x_i,\frac1j\right)\times\{\o\}\\
\neq &\emptyset.
\end{split}
\end{align}
Although we only give the definition of the integer valued function $n$ on $\O_{i,j}^V$, it is not hard to extend $n$ onto $\O$ for which (\ref{E:TTran}) still holds for $\mb P$-a.e. $\o\in\O$, since there are only countable many $\O^V_{i,j}$s.
Thus $F$ is random topological transitive.
\end{proof}

The following results are mainly from \cite{GK}  whose proofs are omitted.

Let $\L$ be a random hyperbolic attractor of $F$, where the $F(\o)$ is uniformly hyperbolic on $\L(\o)$ $\mb P$-a.s.. For definition of random hyperbolic attractors, we refer the reader to \cite{GK}; while for more details about random attractors, we refer the reader to  \cite{A}.
\begin{defn}\label{D:RPO}
Let $\d$ be a strictly positive random variable on $(\O,\mc B,\mb P)$. Then for any $\o\in\O$ a sequence $\{y_n\}_{n\in \mb Z}$ in $\mb E$ is called an $(\o,\d)$ pseudo-orbit of $F$ if
$$d(y_{n+1},F(y_n,\t^n\o))\le \d(\t^{n+1}\o)\quad\text{ for all }n\in \mb Z.$$
For a strictly positive random variable $\e$ and any $\o\in \O$ a point $x\in\mb E$ is said to $(\o,\e)$-shadow the $(\o,\d)$ pseudo-orbit $\{y_n\}_{n\in\mb Z}$ if
$$d(F^n(x,\o),y_n)\le \e(\t^n\o)\quad\text{ for all }n\in\mb Z.$$
\end{defn}
The following lemma is the Proposition 3.6 of \cite{GK}, which is called {\bf Random Shadowing Lemma}.
\begin{lem}\label{L:RShL}
Let the random hyperbolic set $\L$ have local product structure. Then for every tempered random variable $\e>0$ there exists a tempered random variable $\b>0$ such that $\mb P$-a.s. every $(\o,\b)$ pseudo-orbit $\{y_n\}_{n\in\mb Z}$ with $y_n\in \L(\t^n\o)$ can be $(\o,\e)$-shadowed by a point $x\in \L(\o)$. If $2\e$ is chosen as an expansivity characteristic, then the shadowing point $x$ is unique. Moreover, if the $y_n$ are chosen to be random variables such that for $\mb P$-almost all $\o\in\O$ the sequence $\{y_n(\o)\}_{n\in\mb Z}$ is an $(\o,\b)$ pseudo -orbit, then the starting point $x(\o)$ of the corresponding $(\o,\e)$-shadowing orbit depends measurably on $\o$.
\end{lem}

Let $V^s, V^u$ be the local invariant manifolds of $F$ respectively, moreover take them of random size $\eta$ being an expansivity characteristic. Take a small enough strictly positive tempered random variable $\varpi$ the smallness of which is only depending on $\b$, $\eta$ and the hyperbolicity of the system. Then one can find a random variable $k:\O\to\mb N$ such that $\L(\o)$ can be covered by $k(\o)$ open balls of radius less than $\varpi(\o)$ and centres $p_i(\o)$, $i=1,\cdots,k(\o)$ with $p_i:\O\to\mb E$ measurable. Actually $\varpi$ can always be chosen log-integrable, and $k$ can be chosen log-integrable if $F$ is of tempered continuity.

Denote these balls by $B_{\varpi}(p_i,\o)$. Then $A$ is called a random matrix if, for each $\o\in \O$, $A(\o)\in \mb R^{k(\o)\times k(\t\o)}$ such that
$$A(\o)_{i,j}=\begin{cases}1\quad &\text{ if }F(p_i(\o),\o)\in B_{\varpi}(p_j,\t\o)\\ 0 &\text{ otherwise.}\end{cases}$$

\begin{defn}\label{D:RSD}
Let $k:\O\to\mb N^+$ be a random variable, $A$ a corresponding random transition matrix, and define for $\o\in\O$
$$\Sigma_k(\o):=\prod_{i=-\infty}^{\infty}\{1,\cdots,k(\t^i\o)\},$$
$$\Sigma_A(\o):=\{\bar x=(x_i)\in \Sigma_k(\o):A_{x_i,x_{i+1}}(\t^i\o)=1\text{ for all }i\in\mb Z\}.$$
Let $\s$ be the standard (left-) shift. The families $\{\s:\Sigma_k(\o)\to \Sigma_k(\t\o)\}$ and $\{\s:\Sigma_A(\o)\to \Sigma_A(\t\o)\}$ are called random $k$-shift and random subshift of finite type, respectively. Moreover, define $\Sigma_A:=\{(\bar x,\o): \bar x\in \Sigma_A(\o),\o\in\O\}$, which is a measurable bundle over $\O$, and also denote the respective skew-product transformation on $\Sigma_A$ by $\s$.
\end{defn}
One can also define one-sided versions of random $k$-shifts and subshifts of finite type on
$\Sigma_k^+:=\prod_{i=0}^\infty\{1,\cdots, k(\t^i\o)\}\text{ and corresponding }\Sigma_A^+(\o), \o\in\O$ respectively.
Let $C(\Sigma_A^+(\o))$ denote the space of random continuous functions on $\Sigma_A^+$ which are measurable on $\o$ and continuous on $x\in\mb E$ for fixed $\o$. If a random continuous function is $\mb P$-a.s. H\"older continuous with uniform exponent, such function is called random H\"older continuous function. For a random continuous function $\varphi$ and $\o\in\O$  random transfer operators $\mc L_{\varphi}(\o):C(\Sigma_A^+(\o))\to C(\Sigma_A^+(\t\o))$ is defined by
$$(\mc L_{\varphi}(\o)h)(x)=\sum_{y\in\Sigma_A^+(\o):\s y=x}\exp(\varphi(\o,y))h(y)$$
for $h\in C(\Sigma_A^+(\o)),x\in C(\Sigma_A^+(\t\o))$. $\mc L_{\varphi}^*(\o)$ denotes the random dual operator mapping finite signed measures on $C(\Sigma_A^+(\t\o))$ to those on $C(\Sigma_A^+(\o))$ by
$$\int hd\mc L_{\varphi}^*(\o)m=\int \mc L_{\varphi}hdm\text{ for all }h\in C(\Sigma_A^+(\o))$$
for a finite signed measure $m$ on $C(\Sigma_A^+(\t\o))$.

The following theorem is the Theorem 4.3 of \cite{GK}, which is the main result about SRB measures.
\begin{thm}\label{T:SRBMeasure}
Let $F$ be a $C^{1+\a}$ random dynamical system with a random topological transitive hyperbolic attractor $\L(\subset \mb E\times \O)$. Then there exists a unique $F$-invariant measure (SRB measure) $\nu$ supported by $\L$ and characterized by each of the following:
\begin{itemize}
\item[(i)] $h_\nu(F)=\int\sum\l_i^+d\nu$ where $\l_i$ are the Lyapunov exponents corresponding to $\nu$;
\item[(ii)] $\mb P$-a.s. the conditional measures of $\nu_\o$ on the unstable manifolds are absolutely continuous with respect to the Riemannian volume on these submanifolds;
\item[(iii)] $h_\nu(F)+\int fd\nu=\sup_{F-\text{invariant measure } m}\{h_m(F)+\int fd\nu\}$ and the latter is the topological pressure $\pi_F(f)$ of $f$ which satisfies $\pi_{F}(f)=0$;
\item[(iv)] $\nu=\psi \tilde \mu$ where $\tilde \mu$ is the equilibrium state for the two-sided shift $\s$ on $\Sigma_A$ and the function $f\circ \psi$. The measure $\tilde \mu$ can be obtained as a natural extension of the probability measure $\mu$ which is invariant with respect to the one-sided shift $\s$ on $\Sigma^+_A$ and such that $\mc L^*_\eta(\o)\mu_{\t\o}=\mu_\o$ $\mb P$-a.s. where $\eta-f\circ \psi=h-h\circ (\t\times \s)$ for some random H\:older continuous function $h$;
\item[(v)] $\nu$ can be obtained as a weak limit $\nu_\o=\lim_{n\to\infty}F(n,\t^{-n}\o)m_{\t^{-n}\o}$ $\mb P$-a.s. for any measure $m_\o$ absolutely continuous with respect to the Riemannian volume such that sup $m_\o\subset U(\o)$.
\end{itemize}
\end{thm}
Here $f(x,\o)=-\log\| \det D_xF(x,\o)|_{E^u(x,\o)}\|$, where $E^u(x,\o)$ is the invariant unstable Oseledets subspace of $D_xF$ on $(x,\o)$; and $\psi$ is the conjugation between $F$ on $\L$ and $\s$ on $\Sigma_A$ which is constructed based on Markov partitions, then reduce $\s$ on $\Sigma_A$ to the image of an unstable manifold to obtain $\s$ on $\Sigma_A^+$. The H\"older continuity follows from the $C^{1+\a}$ness of $F$, and the existence of Markov partitions are given by Theorem 3.9 of \cite{GK}.


\section{Measurable Selection}\label{S:MST}
In this section, we state a version of measurable selection theorems, which is taken from  \cite{BP}.

\begin{thm}\label{T:BMST}
Let $U,V$ be complete separable metric spaces and $E\subset U\times V$ be a Borel set. If, for each $u\in E$, the section $E_u$ is $\s$-compact there is a Borel selection, $S$, of $E$. Further $proj(E)$ is a Borel set and $\rho_S$ is a Borel measurable function defined on $proj(E)$.
\end{thm}
Here $S$ is said to be a Borel selection of $E$ provided
\begin{itemize}
\item[i)] $S$ is a Borel set;
\item[ii)] $S\subset E$;
\item[iii)] For each $u\in U$, the section $S_u=\{v\in V|\ (u,v)\in S\}$ contains at most one point:
\item[iv)] $proj(S)=proj(E)$.
\end{itemize}
And $\rho_S$ is the function induced by $S$, which assigns to each $u\in proj(E)$ the second coordinate of the unique member of $S$ with first coordinate $u$.


\bibliographystyle{plain}

\end{document}